\documentclass{bull-l}
\usepackage{amssymb,amscd,url}
\usepackage[all]{xy}
\usepackage{pict2e}  % allows lines and vectors of arbitrary slopes in pictures

\begin{document}

\allowdisplaybreaks

%%%%%%%%%%%%%%%%%%%%%%%%%%%%%%%%%%%%%%%%%%%%%%%%%%%%%%%%%%%%%%%%%%%%%%
%% Title and Author Information

\title[Trends and Problems in Arithmetic Dynamics]
      {Current Trends and Open Problems in Arithmetic Dynamics}
\date{\today}

\author[Benedetto, et al.]{Robert Benedetto}
\email{rlbenedetto@amherst.edu}
\address{Department of Mathematics and Statistics, Amherst College,
  Amherst, MA 01002 USA}

%% \author[]{Laura DeMarco}
%% \email{demarco@northwestern.edu}
%% \address{Department of Mathematics, Northwestern University, 2033
%%   Sheridan Road Evanston, IL 60208-2730, USA.  OrcID
%%   0000-0003-3934-8832}

\author[]{Patrick Ingram}
\email{pingram@yorku.ca}
\address{Department of Mathematics and Statistics, York University,
  N520 Ross, 4700 Keele Street, Toronto, ON M3J 1P3, Canada}

\author[]{Rafe Jones}
\email{rfjones@carleton.edu}
\address{Carleton College, Department of Mathematics, Northfield, MN
  55057}

\author[]{Michelle Manes}
\email{mmanes@math.hawaii.edu}
\address{Department of Mathematics, University of Hawaii, 2565
  McCarthy Mall, Honolulu, HI 96822, USA}

\author[]{Joseph H. Silverman}
\email{jhs@math.brown.edu}
\address{Mathematics Department, Box 1917,
         Brown University, Providence, RI 02912 USA. OrcID 0000-0003-3887-3248}

\author[]{Thomas J. Tucker}
\email{thomas.tucker@rochester.edu}
\address{Mathematics Department, 915 Hylan Building, University of
  Rochester, Rochester, NY 14627, USA}

\thanks{Benedetto's research supported by NSF Grant DMS-1501766;
  Manes's research  supported by Simons Collaboration Grant~\#359721;
  Silverman's research supported by Simons Collaboration Grant~\#241309,
Tucker's research supported by NSF Grant DMS-0101636.}


%%%%%%%%%%%%%%%%%%%%%%%%%%%%%%%%%%%%%%%%%%%%%%%%%%%%%%%%%%%%%%%%%%%%%%

% \allowdisplaybreaks

\hyphenation{ca-non-i-cal semi-abel-ian}

%%%%%%%%%%%%%%%%%%%%%%%%%%%%%%%%%%%%%%%%%%%%%%%%%%%%%%%%%%%%%%%%%%%%%%
% Theorem environments

\newtheorem{theorem}{Theorem}[section]
\newtheorem{lemma}[theorem]{Lemma}
\newtheorem{sublemma}[theorem]{Sublemma}
\newtheorem{conjecture}[theorem]{Conjecture}
\newtheorem{question}[theorem]{Question}
\newtheorem{proposition}[theorem]{Proposition}
\newtheorem{corollary}[theorem]{Corollary}
\newtheorem*{claim}{Claim}

\theoremstyle{definition}
% The * surpresses numbering
\newtheorem{definition}[theorem]{Definition}
\newtheorem{example}[theorem]{Example}

\theoremstyle{remark}
\newtheorem{remark}[theorem]{Remark}
\newtheorem*{acknowledgement}{Acknowledgements}

\numberwithin{equation}{section}

%%%%%%%%%%%%%%%%%%%%%%%%%%%%%%%%%%%%%%%%%%%%%%%%%%%%%%%%%%%%%%%%%%%%%%

%%%%%%%% Set Up Environment for Notation %%%%%%%%%%%%%%
% This is currently set to allow quite wide items to be defined
\newenvironment{notation}[0]{%
  \begin{list}%
    {}%
    {\setlength{\itemindent}{0pt}
     \setlength{\labelwidth}{4\parindent}
     \setlength{\labelsep}{\parindent}
     \setlength{\leftmargin}{5\parindent}
     \setlength{\itemsep}{0pt}
     }%
   }%
  {\end{list}}

%%%%%%%% Set Up Environment for Parts in Theorems %%%%%%%%%%%%%%
\newenvironment{parts}[0]{%
  \begin{list}{}%
    {\setlength{\itemindent}{0pt}
     \setlength{\labelwidth}{1.5\parindent}
     \setlength{\labelsep}{.5\parindent}
     \setlength{\leftmargin}{2\parindent}
     \setlength{\itemsep}{0pt}
     }%
   }%
  {\end{list}}
% Use \Part{(a)}, instead of \item[(a)], to ensure upright font
\newcommand{\Part}[1]{\item[\upshape#1]}

%%%%%%%% Set Up Macro for Cases %%%%%%%%%%%%%%
\def\Case#1#2{%
\paragraph{\textbf{\boldmath Case #1: #2.}}\hfil\break\ignorespaces}

%%%%%%%%%%%%%%%%%%
% Greek Alphabet %
%%%%%%%%%%%%%%%%%%
\renewcommand{\a}{\alpha}
\renewcommand{\b}{\beta}
\newcommand{\g}{\gamma}
\renewcommand{\d}{\delta}
\newcommand{\e}{\epsilon}
\newcommand{\f}{\varphi}
\newcommand{\bfphi}{{\boldsymbol{\f}}}
\renewcommand{\l}{\lambda}
\renewcommand{\k}{\kappa}
\newcommand{\lhat}{\hat\lambda}
\newcommand{\m}{\mu}
\newcommand{\bfmu}{{\boldsymbol{\mu}}}
\renewcommand{\o}{\omega}
\renewcommand{\r}{\rho}
\newcommand{\rbar}{{\bar\rho}}
\newcommand{\s}{\sigma}
\newcommand{\sbar}{{\bar\sigma}}
\renewcommand{\t}{\tau}
\newcommand{\z}{\zeta}

\newcommand{\D}{\Delta}
\newcommand{\G}{\Gamma}
\newcommand{\F}{\Phi}
\renewcommand{\L}{\Lambda}

%%%%%%%%%%%%%%%%%%%%
% Fraktur Alphabet %
%%%%%%%%%%%%%%%%%%%%
\newcommand{\ga}{{\mathfrak{a}}}
\newcommand{\gb}{{\mathfrak{b}}}
\newcommand{\gn}{{\mathfrak{n}}}
\newcommand{\gp}{{\mathfrak{p}}}
\newcommand{\gP}{{\mathfrak{P}}}
\newcommand{\gq}{{\mathfrak{q}}}

%%%%%%%%%%%%%%%%%%%
% Barred Alphabet %
%%%%%%%%%%%%%%%%%%%
\newcommand{\Abar}{{\bar A}}
\newcommand{\Ebar}{{\bar E}}
\newcommand{\kbar}{{\bar k}}
\newcommand{\Kbar}{{\bar K}}
\newcommand{\Pbar}{{\bar P}}
\newcommand{\Sbar}{{\bar S}}
\newcommand{\Tbar}{{\bar T}}

%%%%%%%%%%%%%%%%%%%%%%%%%
% Calligraphic Alphabet %
%%%%%%%%%%%%%%%%%%%%%%%%%
\newcommand{\Acal}{{\mathcal A}}
\newcommand{\Bcal}{{\mathcal B}}
\newcommand{\Ccal}{{\mathcal C}}
\newcommand{\Dcal}{{\mathcal D}}
\newcommand{\Ecal}{{\mathcal E}}
\newcommand{\Fcal}{{\mathcal F}}
\newcommand{\Gcal}{{\mathcal G}}
\newcommand{\Hcal}{{\mathcal H}}
\newcommand{\Ical}{{\mathcal I}}
\newcommand{\Jcal}{{\mathcal J}}
\newcommand{\Kcal}{{\mathcal K}}
\newcommand{\Lcal}{{\mathcal L}}
\newcommand{\Mcal}{{\mathcal M}}
\newcommand{\Ncal}{{\mathcal N}}
\newcommand{\Ocal}{{\mathcal O}}
\newcommand{\Pcal}{{\mathcal P}}
\newcommand{\Qcal}{{\mathcal Q}}
\newcommand{\Rcal}{{\mathcal R}}
\newcommand{\Scal}{{\mathcal S}}
\newcommand{\Tcal}{{\mathcal T}}
\newcommand{\Ucal}{{\mathcal U}}
\newcommand{\Vcal}{{\mathcal V}}
\newcommand{\Wcal}{{\mathcal W}}
\newcommand{\Xcal}{{\mathcal X}}
\newcommand{\Ycal}{{\mathcal Y}}
\newcommand{\Zcal}{{\mathcal Z}}

%%%%%%%%%%%%%%%%%%%%%%%%%%%%
% Blackboard Bold Alphabet %
%%%%%%%%%%%%%%%%%%%%%%%%%%%%
\renewcommand{\AA}{\mathbb{A}}
\newcommand{\BB}{\mathbb{B}}
\newcommand{\CC}{\mathbb{C}}
\newcommand{\FF}{\mathbb{F}}
\newcommand{\GG}{\mathbb{G}}
\newcommand{\NN}{\mathbb{N}}
\newcommand{\PP}{\mathbb{P}}
\newcommand{\QQ}{\mathbb{Q}}
\newcommand{\RR}{\mathbb{R}}
\newcommand{\ZZ}{\mathbb{Z}}

%%%%%%%%%%%%%%%%%%%%%%%%%%
% Boldface Math Alphabet %
%%%%%%%%%%%%%%%%%%%%%%%%%%
\newcommand{\bfa}{{\boldsymbol a}}
\newcommand{\bfb}{{\boldsymbol b}}
\newcommand{\bfc}{{\boldsymbol c}}
\newcommand{\bfd}{{\boldsymbol d}}
\newcommand{\bfe}{{\boldsymbol e}}
\newcommand{\bff}{{\boldsymbol f}}
\newcommand{\bfg}{{\boldsymbol g}}
\newcommand{\bfi}{{\boldsymbol i}}
\newcommand{\bfj}{{\boldsymbol j}}
\newcommand{\bfp}{{\boldsymbol p}}
\newcommand{\bfr}{{\boldsymbol r}}
\newcommand{\bfs}{{\boldsymbol s}}
\newcommand{\bft}{{\boldsymbol t}}
\newcommand{\bfu}{{\boldsymbol u}}
\newcommand{\bfv}{{\boldsymbol v}}
\newcommand{\bfw}{{\boldsymbol w}}
\newcommand{\bfx}{{\boldsymbol x}}
\newcommand{\bfy}{{\boldsymbol y}}
\newcommand{\bfz}{{\boldsymbol z}}
\newcommand{\bfA}{{\boldsymbol A}}
\newcommand{\bfF}{{\boldsymbol F}}
\newcommand{\bfB}{{\boldsymbol B}}
\newcommand{\bfD}{{\boldsymbol D}}
\newcommand{\bfG}{{\boldsymbol G}}
\newcommand{\bfI}{{\boldsymbol I}}
\newcommand{\bfM}{{\boldsymbol M}}
\newcommand{\bfP}{{\boldsymbol P}}
\newcommand{\bfzero}{{\boldsymbol{0}}}
\newcommand{\bfone}{{\boldsymbol{1}}}

%%%%%%%%%%%%%%%%%%%%%%%%%%%%%%
% Miscellaneous New Commands %
%%%%%%%%%%%%%%%%%%%%%%%%%%%%%%
\newcommand{\arb}{{\textup{arb}}}
\newcommand{\Aut}{\operatorname{Aut}}
\newcommand{\codim}{\operatorname{codim}}
\newcommand{\Comp}{\operatorname{Comp}}
\newcommand{\Compavg}{\operatorname{Comp}_{\textup{avg}}}
\newcommand{\Crit}{\operatorname{Crit}}
\newcommand{\crit}{{\textup{crit}}}
\newcommand{\Density}{\operatorname{\textup{\textsf{Density}}}}
\newcommand{\Disc}{\operatorname{Disc}}
\newcommand{\Div}{\operatorname{Div}}
\newcommand{\dyn}{{\textup{dyn}}}
\newcommand{\Dom}{\operatorname{Dom}}
\newcommand{\End}{\operatorname{End}}
\newcommand{\Fbar}{{\bar{F}}}
\newcommand{\Fix}{\operatorname{Fix}}
\newcommand{\Gal}{\operatorname{Gal}}
\newcommand{\GL}{\operatorname{GL}}
\newcommand{\GR}{\operatorname{\mathcal{G\!R}}}
\newcommand{\Hom}{\operatorname{Hom}}
\newcommand{\Image}{\operatorname{Image}}
\newcommand{\Isom}{\operatorname{Isom}}
\newcommand{\hhat}{{\hat h}}
\newcommand{\Ker}{{\operatorname{ker}}}
\newcommand{\limstar}{\lim\nolimits^*}
\newcommand{\limstarn}{\lim_{\hidewidth n\to\infty\hidewidth}{\!}^*{\,}}
\newcommand{\Mat}{\operatorname{Mat}}
\newcommand{\maxplus}{\operatornamewithlimits{\textup{max}^{\scriptscriptstyle+}}}
\newcommand{\MOD}[1]{~(\textup{mod}~#1)}
\newcommand{\Mor}{\operatorname{Mor}}
\newcommand{\Moduli}{\mathcal{M}}
\newcommand{\Norm}{{\operatorname{\mathsf{N}}}}
\newcommand{\notdivide}{\nmid}
\newcommand{\normalsubgroup}{\triangleleft}
\newcommand{\NS}{\operatorname{NS}}
\newcommand{\onto}{\twoheadrightarrow}
\newcommand{\ord}{\operatorname{ord}}
\newcommand{\Orbit}{\mathcal{O}}
\newcommand{\Per}{\operatorname{Per}}
\newcommand{\PrePer}{\operatorname{PrePer}}
\newcommand{\PGL}{\operatorname{PGL}}
\newcommand{\Pic}{\operatorname{Pic}}
\newcommand{\Poly}{\operatorname{Poly}}
\newcommand{\Prob}{\operatorname{Prob}}
\newcommand{\Proj}{\operatorname{Proj}}
\newcommand{\PBerk}{\mathsf{P}_{\textup{Berk}}^1}
\newcommand{\Qbar}{{\bar{\QQ}}}
\newcommand{\rank}{\operatorname{rank}}
\newcommand{\Rat}{\operatorname{Rat}}
\newcommand{\Resultant}{\operatorname{Res}}
\newcommand{\Set}{{\mathsf{Set}}}
\renewcommand{\setminus}{\smallsetminus}
\newcommand{\sgn}{\operatorname{sgn}} 
\newcommand{\shafdim}{\operatorname{ShafDim}}
\newcommand{\SL}{\operatorname{SL}}
\newcommand{\Span}{\operatorname{Span}}
\newcommand{\Spec}{\operatorname{Spec}}
\renewcommand{\ss}{\textup{ss}}
\newcommand{\stab}{\textup{stab}}
\newcommand{\evenstab}{{\textup{e-stable}}}
\newcommand{\Support}{\operatorname{\textup{\textsf{Support}}}}
\newcommand{\tors}{{\textup{tors}}}
\newcommand{\tr}{{\textup{tr}}} 
\newcommand{\Trace}{\operatorname{Trace}}
\newcommand{\trianglebin}{\mathbin{\triangle}} % symmetric set difference
\newcommand{\UHP}{{\mathfrak{h}}}    % Upper half plane
\newcommand{\Zsig}{\operatorname{\textup{\textsf{Zsig}}}}
\newcommand{\<}{\langle}
\renewcommand{\>}{\rangle}

\newcommand{\pmodintext}[1]{~\textup{(mod}~#1\textup{)}}
\newcommand{\ds}{\displaystyle}
\newcommand{\longhookrightarrow}{\lhook\joinrel\longrightarrow}
\newcommand{\longonto}{\relbar\joinrel\twoheadrightarrow}
\newcommand{\SmallMatrix}[1]{%
  \left(\begin{smallmatrix} #1 \end{smallmatrix}\right)}

%% This creates an \xrightarrow that's dashed. This is used to create a long dashed arrow.
\makeatletter
\newcommand{\xdashrightarrow}[2][]{\ext@arrow 0359\rightarrowfill@@{#1}{#2}}
\def\rightarrowfill@@{\arrowfill@@\relax\relbar\rightarrow}
\def\arrowfill@@#1#2#3#4{%
  $\m@th\thickmuskip0mu\medmuskip\thickmuskip\thinmuskip\thickmuskip
   \relax#4#1
   \xleaders\hbox{$#4#2$}\hfill
   #3$%
}
\newcommand{\longdashrightarrow}{\xdashrightarrow{\hspace{2em}}}
\makeatother

\newcommand{\WhoWrite}[1]{\noindent \framebox{This section written by: #1}\par}
\def\Who#1{{ [\textit{#1}]}}
\def\WhoDone#1{{ [\textit{#1}: Done]}}
%% Uncomment the next line to remove this information
\def\WhoWrite#1{}  \def\Who{} \def\WhoDone#1{}

%%%%%%%%%%%%%%%%%%%%%%%%%%%%%%%%%%%%%%%%%%%%%%%%%%%%%%%%%%%%%%%%%%%%%%

\begin{abstract}
Arithmetic dynamics is the study of number theoretic properties of
dynamical systems.  A relatively new field, it draws inspiration
partly from dynamical analogues of theorems and conjectures in
classical arithmetic geometry, and partly from $p$-adic analogues of
theorems and conjectures in classical complex dynamics.  In this
article we survey some of the motivating problems and some of the
recent progress in the field of arithmetic dynamics.
\end{abstract}

\maketitle

\setcounter{tocdepth}{1}
\tableofcontents

%%%%%%%%%%%%%%%%%%%%%%%%%%%%%%%%%%%%%%%%%%%%%%%%%%%%%%%%%%%%%%%%%%%%%%
\section{Introduction}
\label{section:introduction}
%%%%%%%%%%%%%%%%%%%%%%%%%%%%%%%%%%%%%%%%%%%%%%%%%%%%%%%%%%%%%%%%%%%%%%

In classical real and complex dynamics, one studies topological and
analytic properties of orbits of points under iteration of self-maps
of~$\RR$ or~$\CC$, or more generally self-maps of real or complex
manifolds.  In \emph{arithmetic dynamics}, a more recent subject, one
likewise studies properties of orbits of self-maps, but with a number
theoretic flavor.  In this article we describe some of the motivating
problems and recent advances in this relatively new field of
arithmetic dynamics.

Many of the motivating problems in arithmetic dynamics come via
analogy with classical problems in arithmetic geometry. A large part
of arithmetic geometry involves studying rational and integral points
on algebraic varieties, and it is fair to say that many deep results
in the subject involve abelian varieties, which are algebraic
varieties admitting a group structure. Table~\ref{table:dictionary}
provides a dictionary that suggests ways in which classical theorems
and conjectures in arithmetic geometry lead naturally to dynamical
problems of an arithmetic nature.

\begin{table}[ht]
\begin{center}
\begin{picture}(300,115)(-20,-35)

\put(0,65){\makebox(100,12)[c]{\textbf{Arithmetic Geometry}}}
\put(147,65){\makebox(120,12)[c]{\textbf{Dynamical Systems}}}

\put(-5,35){\framebox(110,25)[l]{%
    \begin{tabular}{l}
       rational and integral\\
       points on varieties\\
    \end{tabular}
}}

\put(105,48){\vector(1,0){40}}
\put(145,48){\vector(-1,0){40}}

\put(145,35){\framebox(125,25)[l]{%
    \begin{tabular}{l}
       rational and integral\\
       points in orbits\\
    \end{tabular}
}}

\put(-5,0){\framebox(110,25)[l]{%
    \begin{tabular}{l}
       torsion points on\\
       abelian varieties\\
    \end{tabular}
}}

\put(105,13){\vector(1,0){40}}
\put(145,13){\vector(-1,0){40}}

\put(145,0){\framebox(125,25)[l]{%
    \begin{tabular}{l}
       periodic and preperiodic\\
       points of rational maps\\
    \end{tabular}
}}
%%%%%%%%%%%%%%%%%%%%%%%%%%%%%%%%%%%%%%%%%%%%%%%%%%%%%%%%%%%%
\put(-5,-35){\framebox(110,25)[l]{%
    \begin{tabular}{l}
       abelian varieties with\\
       complex multiplication\\
    \end{tabular}
}}
\put(105,-22){\vector(1,0){40}}
\put(145,-22){\vector(-1,0){40}}

\put(145,-35){\framebox(125,25)[l]{%
    \begin{tabular}{l}
       post-critically finite\\
       rational maps\\
    \end{tabular}
}}
\end{picture}
\end{center}
\caption{An Arithmetical/Dynamical Dictionary \cite[\S6.5]{MR2884382}}
\label{table:dictionary}
\end{table}

A second source of inspiration for arithmetic dynamics comes from a
basic principle in number theory: before studying a problem over the
field of rational numbers~$\QQ$, one should first study it over all
completions of~$\QQ$, including in particular $p$-adic fields.  The
field~$\QQ_p$ of $p$-adic rationals comes equipped with an absolute
value, so many concepts from complex dynamics can be translated to the
$p$-adic setting, and one can try to prove $p$-adic analogues of
complex dynamical theorems.  However, since the $p$-adic absolute
value satisfies the nonarchimedean triangle inequality
$\|x+y\|\le\max\bigl\{\|x\|,\|y\|\bigr\}$, there are striking
differences between the $p$-adic and complex theories, many of them
traceable to two facts: (1) $p$-adic fields are totally disconnected;
(2) the smallest complete algebriacally closed field
containing~$\QQ_p$ is not locally compact. Despite, or possibly
because of, these differences, the study of $p$-adic dynamics is a
thriving branch of arithmetic dynamics.

%% *** Should we mention dynamics over finite fields here in the introduction

The goal of this article is to describe motivating problems and
conjectures in the field of arithmetic dynamics and to give a brief
description of recent progress. We start in
Sections~\ref{section:abstractDS} and~\ref{section:background} with an
overview of the terminology, notation, and background on dynamical
systems, number theory, and algebraic geometry that we will need. As
noted at the end of Section~\ref{section:background}, we have
attempted to make this article accessible by generally restricting
attention to maps of projective space, thereby obviating the need for
more advanced material from algebraic geometry. Our survey of
arithmetic dynamics then commences in Section~\ref{section:ub} and
concludes in Section~\ref{section:lgqd}. These nineteen sections are mostly
independent of one another, albeit liberally sprinkled with
instructive cross-references.

\subsection*{A Brief Timeline of the Early Days}
The study of Galois groups of iterated polynomials was initiated by
Odoni in a series of papers~\cite{MR805714,MR813379,MR962740} during
the mid-1980s. Beyond this, it appears that the first paper to
describe dynamical analogues for a wide range of classical problems in
arithmetic geometry did so not for polynomial maps or projective
space, but rather for~K3 surfaces admitting an automorphism of
infinite order.  The~1991 paper of
Silverman~\cite{silverman:K3heights} discusses dynamical analogues of
various problems, including: (1)~Uniform boundedness of periodic
points; (2)~Periodic points on subvarieties; (3)~Lower bounds for
canonical heights; (4)~Open image of Galois groups; (5)~Integral
points in orbits.  The rest of the 1990s saw an explosion of papers in
which numerous researchers began to study a wide variety of problems
in arithmetic dynamics.

The serious study of $p$-adic dynamics starts with a 1983 paper of
Herman and Yoccoz~\cite{MR730280}, two of the world's leading complex
dynamicists, in which they investigate a $p$-adic analogue of a
classical problem in complex dynamics. Later in the 1980s and 90s, the
physics literature includes several
papers~\cite{arrowsmith:padichorseshoe,arrowsmith:padicsiegeldisk,unpubBenMenahem,thiran:padicdynamics}
that discuss $p$-adic dynamics for polynomial maps, but the deeper
study of $p$-adic dynamics really took off with the Ph.D.\ theses of
Benedetto~\cite{benedetto:thesis,MR1781331} in~1998 and
Rivera-Letelier~\cite{riverathesis,MR2040006} in~2000.  The dynamics
of iterated $p$-adic power series was investigated by
Lubin~\cite{MR1310863} in a~1994 paper. We also mention that there is
an extensive literature on ergodic theory over~$p$-adic fields which,
although not within the purview of this survey, dates back to at
least~1975~\cite{MR0369301}.

%%%%%%%%%%%%%%%%%%%%%%%%%%%%%%%%%%%%%%%%%%%%%%%%%%%%%%%%%%%%%%%%%%%%%%
\section{Abstract Dynamical Systems\WhoDone{Joe}}
\label{section:abstractDS}
\WhoWrite{Joe}
%%%%%%%%%%%%%%%%%%%%%%%%%%%%%%%%%%%%%%%%%%%%%%%%%%%%%%%%%%%%%%%%%%%%%%

A \emph{discrete dynamical system} consists of a set~$S$ and a
self-map
\[
  f:S\longrightarrow S.
\]
We write
\[
  f^n := \underbrace{f\circ f\circ\cdots\circ f}_{\text{$n$ copies}}
\]
for the $n$th iterate of~$f$, with the convention that $f^0$ is the
identity map, and we write
\[
  \Orbit_f(\a) := \bigl\{ f^n(\a) : n\ge0 \bigr\}
\]
for the \emph{forward orbit} of an element~$\a\in S$.  A fundamental
problem of discrete dynamics is to classify the points in~$S$
according to the behavior of their orbits. Frequently this behavior is
described according to some additional properties of the set~$S$,
which for example might be a topological space or a metric space or an
analytic space. But even in the setting of abstract sets, we can make
interesting distinctions among orbits.

\begin{definition}  
Let $f:S\to S$ be a dynamical system, and let $\a\in S$. We say
that~$\a$ is
\begin{center}
  \begin{tabular}{rl}
  \emph{Wandering}\footnotemark  &\quad\text{if $\Orbit_f(\a)$ is infinite,} \\
  \emph{Preperiodic}&\quad\text{if $\Orbit_f(\a)$ is finite,} \\
  \emph{Periodic}&\quad\text{if $f^n(\a)=\a$ for some $n\ge1$.} \\
  \end{tabular}
\end{center}
\footnotetext{We mention that our definition of \emph{wandering point}
  is different from the classical definition for dynamical systems on
  topological, respectively measure, spaces, where~$\a$ is said to be wandering if
  there is a neighborhood~$U$ of~$\a$ such that~$f^n(U)\cap U$ is
  empty, respectively has measure~$0$, for all sufficiently large~$n$.}
If~$\a$ is periodic, the smallest~$n\ge1$ satisfying~$f^n(\a)=\a$ is
called the (\emph{exact}) \emph{period of~$\a$}.  We set the useful notation
\begin{align*}
  \PrePer(f,S)  &:= \bigl\{ \a\in S : \text{$f^n(\a)=f^m(\a)$ for some $n>m\ge0$} \bigr\}, \\
  \Per(f,S)     &:= \bigl\{ \a\in S : \text{$f^n(\a)=\a$ for some $n\ge1$} \bigr\}, \\
  \Per_n(f,S)   &:= \bigl\{ \a\in S : f^n(\a)=\a \bigr\}, \\
  \Per_n^*(f,S) &:= \bigl\{ \a\in S : \text{$f^n(\a)=\a$ and $f^i(\a)\ne\a$ for all $1\le i<n$} \bigr\}.
\end{align*}
Thus, $\PrePer(f,S)$ is the set of preperiodic points, $\Per(f,S)$ is
the set of periodic points, $\Per_n(f,S)$ is the set of periodic
points of period dividing~$n$, and $\Per_n^*(f,S)$ is the set of
periodic points of period exactly equal to~$n$.
\end{definition}

In this way the map~$f$ partitions the points of~$S$ into two subsets depending on
whether a point's orbit is finite or infinite, and it further
partitions the former into two subsets depending on whether a point's
orbit is strictly periodic.  More generally, if~$S$ has more
structure, for example if~$S$ is a metric space and~$f$ is continuous,
one classically partitions~$S$ into the Fatou set, consisting of
points where the iterates of~$f$ behave ``nicely'', and the
complementary Julia set, consisting of points where the iterates
of~$f$ behave ``chaotically.''

\begin{example}
\label{example:Gprepertorsion}
Let $G$ be a group, let $d\ge2$, and let $f:G\to G$ be the $d$th
power map $f(\g)=\g^d$. Then it is a fun exercise to verify that
\[
%%  \{\g\in G : \text{$\g$ is preperiodic for $f$}\}
  \PrePer(f,G) =
  \{\g\in G : \text{$\g$ has finite order}\}.
\]
\end{example}

Moving to the next level, dynamicists often study collections of
maps. Suppose that~$\Fcal$ is a collection of maps $f:S\to S$. Then a
basic meta-problem is to classify the maps~$f$ in~$\Fcal$ according to
the dynamical behavior of~$f$. Here is a typical example.

\begin{question}
Let~$S$ be a set, and let $T\subseteq S$ be a subset.  Classify those maps
$f:S\to S$ with the property that for every $\a\in S$, the set
\[
  \bigl\{ n \ge 0 : f^n(\a) \in T \bigr\}
\]
has some nice structure, for example, it is the union of finitely many
one-sided arithmetic progressions.
\end{question}

%%%%%%%%%%%%%%%%%%%%%%%%%%%%%%%%%%%%%%%%%%%%%%%%%%%%%%%%%%%%%%%%%%%%%%
\section{Background: Number Theory and Algebraic Geometry\WhoDone{Joe}}
\label{section:background}
\WhoWrite{Joe}
%%%%%%%%%%%%%%%%%%%%%%%%%%%%%%%%%%%%%%%%%%%%%%%%%%%%%%%%%%%%%%%%%%%%%%

In this paper we use $\ZZ,\QQ,\RR,\CC$ as usual to denote the
integers, rational numbers, real numbers, and complex numbers. We
write~$\FF_q$ for a finite field with~$q$ elements. We fix an
algebraic closure~$\Qbar$ for~$\QQ$. Elements of~$\Qbar$ are called
\emph{algebraic numbers}. A \emph{number field} is a finite
extension~$K$ of~$\QQ$.

The completions of~$\ZZ$ and~$\QQ$ for their $p$-adic valuations are
denoted~$\ZZ_p$ and~$\QQ_p$, respectively, and a completion of an
algebraic closure of~$\QQ_p$ is denoted~$\CC_p$.

At times we will want a way to measure the size of an infinite
collection of primes.

\begin{definition}
\label{definition:density} 
The (\emph{natural}) \emph{density} of a set of
primes~$\Pcal\subset\ZZ$ is the following quantity, provided that the
limit exists:
\[
  \Density(\Pcal) = \lim_{X\to\infty}
  \frac{\#\{p\in\Pcal : p\le X\}}{\#\{\text{primes $p$} : p\le X\}}.
\]
\end{definition}

Let $N\ge1$. Fora given field~$K$, we write~$\Kbar$ to denote a
separable closure of~$K$.\footnote{Usually we will take $K$ to either
  have characteristic~$0$ or be a finite field, in which case~$\Kbar$
  is an algebraic closure of~$K$.} \emph{Projective $N$-space}
(\emph{over~$K$}), denoted~$\PP^N$, is the quotient space
\[
  \PP^N := ( \Kbar^{N+1} \setminus \bfzero ) / \Kbar^*,
\]
where the group~$\Kbar^*$ acts on $(N+1)$-tuples via scalar
multiplication.  If we want to specify that the coordinates of
the~$(N+1)$-tuples are in the field~$K$, we write
\[
  \PP^N(K) := ( K^{N+1} \setminus \bfzero ) / K^*.
\]
We write
\[
  [a_0,a_1,\ldots,a_N] \in \PP^N
\]
to denote the point in~$\PP^N$ corresponding to the $(N+1)$-tuple
$(a_0,\ldots,a_N)$.  The \emph{field of definition}~$K(P)$ of a point
$P=[a_0,\ldots,a_N]\in\PP^N$ is defined by choosing any index~$i$ such
that~$a_i\ne0$ and setting
\[
  K(P) := K(a_0/a_i,a_1/a_i,\ldots,a_N/a_i).
\]

The Galois group~$\Gal(\Kbar/K)$ acts on~$\PP^N$ via its action
on the coordinates of an~$N$-tuple, and it is a consequence of Hilbert's Theorem~90
that
\[
  \PP^N(K) = \bigl\{ P\in\PP^N : \text{$\s(P)=P$ for all $\s\in\Gal(\Kbar/K)$} \bigr\}.
\]
More generally, the field of definition of~$P$ may be defined Galois-theoretically as
\[
  K(P) := \text{Subfield of $\Kbar$ fixed by
    $\bigl\{\s\in\Gal(\Kbar/K) : \s(P)=P\bigr\}$.}
\]

We next discuss maps from~$\PP^N$ to itself.  A \emph{rational map}
$f:\PP^N\dashrightarrow\PP^N$ is specified by an $(N+1)$-tuple of
homogeneous polynomials $f_0,\ldots,f_N\in\Kbar[X_0,\ldots,X_N]$ of
the same degree,
\[
  f = [f_0,\ldots,f_N] : \PP^N \longdashrightarrow \PP^N.
\]
Since the polynomial ring~$\Kbar[X_0,\ldots,K_N]$ is a unique
factorization domain, we always cancel any factors common to all
of~$f_0,\ldots,f_N$. With this proviso, the \emph{degree of~$f$} is
defined to be the common degree of the polynomials~$f_0,\ldots,f_N$.
The map~$f$ is said to be \emph{dominant} if the
polynomials~$f_0,\ldots,f_N$ do not themselves satisfy a non-trivial
polynomial relation.

If~$P=[a_0,\ldots,a_N]\in\PP^N$, then we define $f(P)$ to be
\[
  f(P) := \bigl[ f_0(a_0,\ldots,a_N),\ldots,f_N(a_0,\ldots,a_N) \bigr],
\]
provided at least one of the coordinates is non-zero. Note that if we
use a different representative for~$P$, say~$P=[ca_0,\ldots,ca_N]$ for
some $c\in\Kbar^*$, then the coordinates of~$f(P)$ are all multiplied
by~$c^{\deg(f)}$, so we get the same point in~$\PP^N$.

For $N=1$, a point~$P\in\PP^1$ is a \emph{critical point} for the
rational map $f:\PP^1\to\PP^1$ if~$f$ fails to be locally one-to-one
in a neighborhood of~$P$. Writing~$f$ as a rational
function~$f(x)\in\Kbar(x)$, the critical points of~$f$ are the points
where the derivative~$f'(x)$ vanishes, with some adjustment required
if~$P$ or~$f(P)$ is the point at infinity. Maps whose critical points
are all preperiodic are said to be \emph{post-critically finite}, or
PCF for short. The collection of PCF maps plays an important role in
the study of $1$-dimensional dynamics.

Returning now to arbitrary $N\ge1$, the \emph{indeterminacy locus
  of~$f$}, denoted~$I(f)$, is the set of points at which~$f$ is not
well-defined, i.e.,
\[
  I(f) = \bigl\{ [a_0,\ldots,a_N]\in\PP^N : 
   f_0(a_0,\ldots,a_N)=\cdots=f_N(a_0,\ldots,a_N)=0 \bigr\}.
\]
The map~$f$ is said to be a \emph{morphism} if~$I(f)=\emptyset$, in
which case~$f$ is well-defined at every point of~$\PP^N$, and we
write~$f$ using a solid arrow $f:\PP^N\to\PP^N$. For $N=1$, it turns
out that every rational map is a morphism, but this is no longer true
for~$N\ge2$. We set the notation
\[
  \End_d^N(K) :=
  \left\{\begin{tabular}{@{}l@{}}
  morphisms $f:\PP^N\to\PP^N$ of degree $d$ given by\\
  polynomials with coefficients in the field $K$\\
  \end{tabular}
  \right\}.
\]

The fact that a rational map $f:\PP^N\dashrightarrow\PP^N$ is not
well-defined at every point is potentially a big problem if we want to
study orbits of points under iteration of~$f$. This prompts the following
notation for the set of points having well-defined orbits:
\[
  \PP^N_f := \bigl\{ P \in \PP^N : \text{$f^n(P)\notin I(f)$ for all $n\ge0$} \bigr\}.
\]
If we further want the coordinates of~$P$ to be in the field~$K$, we
write~$\PP^N(K)_f$.

The zero set of a non-zero homogeneous polynomial
$g\in\Kbar[X_0,\ldots,X_N]$ is a well-defined subset of~$\PP^N$. The
complement of~$\{g=0\}$ is called a \emph{Zariski open subset
  of~$\PP^N$}, and these form a base of open sets for the
\emph{Zariski topology of~$\PP^N$}. For our purposes, the most
important concept is that of density for the Zariski topology.

\begin{definition}
A subset $S\subset\PP^N$ is \emph{Zariski dense} if
\[
  \bigl\{ g\in\Kbar[X_0,\ldots,X_N] :
  \text{$g(P)=0$ for every $P\in S$} \bigr\} = \{0\},
\]
i.e., there are no non-zero homogeneous polynomials
$g\in\Kbar[X_0,\ldots,X_N]$ that vanish at every point of~$S$.
\end{definition}

\subsection*{Note to the Reader}
To make this article accessble to a wide audience, we have mostly
restricted attention to dynamics on projective space~$\PP^N$ over
number fields, and frequently we restrict to maps of~$\PP^1$ defined
over~$\QQ$.  Some sections conclude with a short subsection describing
how the material generalizes to other algebraic varieties.  In these
addenda we have assumed that the reader is familiar with some basic
algebraic geometry terminology and constructions, including, but not
limited to, the following concepts.

A (\emph{projective}) \emph{algebraic set}~$V$ is a subset of~$\PP^N$ defined by
the set of solutions in~$\PP^N$ of a system of equations
\[
  g_1(X_0,\ldots,X_N) = \cdots = g_r(X_0,\ldots,X_N) = 0
\]
given by a collection of homogeneous polynomials
$g_1,\ldots,g_r\in\Kbar[X_0,\ldots,X_N]$. If~$g_1,\ldots,g_r$ have coefficients in~$K$,
then we say that~$V$ is \emph{defined over~$K$}, and then
\[
  V(K) = V \cap \PP^N(K).
\]
A \emph{rational map} $F:V\to V$ is the restriction to~$V$ of a
collection of rational maps $f_1,\ldots,f_s:\PP^N\to\PP^N$ having the
property that
\[
  f_i\bigl(V\setminus I(f_i)\bigr)\subseteq V
  \quad\text{and}\quad
  \text{$f_i(P) = f_j(P)$ for all $P\in V\setminus\bigl(I(f_i)\cup I(f_j)\bigr)$.}
\]
In other words, we create~$F$ in the standard way by gluing the~$f_i$
along the subsets of~$V$ on which they agree.  The map~$F$ is
well-defined at all points of~$V$ outside the \emph{indeterminacy
  locus}~$I(F):=\bigcap_{i=1}^s I(f_i)$. Thus if~$P\notin I(f_j)$ for
some~$j$, then $F(P)=f_j(P)$.  If~$I(F)=\emptyset$, then we say
that~$F$ is a \emph{morphism}, or an \emph{endomorphism of~$V$}.

%%%%%%%%%%%%%%%%%%%%%%%%%%%%%%%%%%%%%%%%%%%%%%%%%%%%%%%%%%%%%%%%%%%%%%
\section{Uniform Boundedness of (Pre)Periodic Points\WhoDone{Joe}} 
\label{section:ub}
\WhoWrite{Joe}
%%%%%%%%%%%%%%%%%%%%%%%%%%%%%%%%%%%%%%%%%%%%%%%%%%%%%%%%%%%%%%%%%%%%%%

A rational function $f(x)\in K(x)$ induces a map on the $K$-rational
points of the projective line,
\[
  f : \PP^1(K)\longrightarrow\PP^1(K).
\]
When $K=\RR$ or $K=\CC$, the distribution of the (pre)periodic points
is intimately connected to the dynamics of~$f$. For fields of
arithmetic interest, there is a famous, frequently rediscovered,
finiteness theorem of Northcott.

\begin{theorem}[Northcott, 1950, \cite{northcott:periodicpoints}]
\label{theorem:northcott}  
Let $N\ge1$ and $d\ge2$, and let~$K/\QQ$ be a number field. Then for
all degree~$d$ morphisms $f:\PP^N\to\PP^N$ defined over~$K$, the set
of $K$-rational preperiodic points
\[
  \PrePer\bigl(f,\PP^N(K)\bigr)~\text{is finite.}
\]
\end{theorem}

In view of Northcott's theorem, it is natural to inquire as to the
possible size of the finite set~$\PrePer\bigl(f,\PP^N(K)\bigr)$, which
leads to the following fundamental conjecture.

\begin{conjecture}[Dynamical Uniform Boundedness Conjecture:~Mor\-ton--Sil\-ver\-man, 1994,
\cite{mortonsilverman:rationalperiodicpoints}]
\label{conjecture:MSUBC}
Fix integers $N\ge1$, $d\ge2$, and $D\ge1$. There is a constant
$C(N,d,D)$ such that for all degree~$d$ morphisms $f:\PP^N\to\PP^N$
defined over a number field~$K$ of degree~$[K:\QQ]=D$, the number
of~$K$-rational preperiodic points is uniformly bounded,
\[
  \#\PrePer\bigl(f,\PP^N(K)\bigr) \le C(N,d,D).
\]
\end{conjecture}

Very little is known about this conjecture, although there are weaker
bounds in which the constant~$C(N,d,D)$ is replaced by a constant that
depends in various weak ways on the number of primes of bad reduction
of~$f$; see for example~\cite{MR1813109,MR2339471}, as well as the
discussion in Section~\ref{section:goodred}. Indeed, even for the
family of polynomials~$x^2+c$, we have only fragmentary knowledge, as
in the following result.

\begin{theorem}
  \label{theorem:x2cperpt}
  For $c\in\QQ$, let $f_c(x)=x^2+c$. Let $N\ge1$.
  \begin{parts}
    \Part{(a)}
    For $N\in\{1,2,3\}$, there are infinitely many $c\in\QQ$ such
    that~$f_c(x)$ has a periodic point $\a\in\QQ$ of exact
    period~$N$.
    \Part{(b)}
    \textup{(Morton, 1992, \cite{morton:quadraticmaps};
      Flynn--Poonen--Schaefer, 1997, \cite{flynnetal:cyclesofquads};
      Stoll, 2008, \cite{arxiv:0803.2836})}
    For $N\in\{4,5,6\}$, there are no values of $c\in\QQ$ such
    that~$f_c(x)$ has a periodic point $\a\in\QQ$ of exact
    period~$N$, with the caveat that the proof for $N=6$ is conditional
    on the conjecture of Birch and Swinnerton-Dyer.
  \end{parts}
\end{theorem}

The proofs of the three parts of Theorem~\ref{theorem:x2cperpt}(b) are
geometric and global in nature. See~\cite{MR3562026} for a local proof
of~$N=4$ and, with an additional assumption,~\text{$N=5$}.

\begin{conjecture}
\label{conjecture:perle3}
\textup{(Poonen, 1998, \cite{poonen:uniformboundrefined})}
For all $c\in\QQ$, the set of periodic points \text{$\Per(x^2+c,\QQ)$} contains
no points of period strictly greater than~$3$.
\end{conjecture}

\begin{theorem}
\textup{(Poonen, 1998, \cite{poonen:uniformboundrefined})}
If Conjecture~$\ref{conjecture:perle3}$ is true,
then
\[
  \sup_{c\in\QQ} \#\PrePer\bigl(x^2+c,\PP^1(\QQ)\bigr) = 9,
\]
with the maximum occurring for infinitely many values of~$c$,
including for example $c=-\frac{29}{16}$ and $c=-\frac{21}{16}$.
\textup(All possible preperiodic structures are described
in~\cite{poonen:uniformboundrefined}.\textup)
\end{theorem}

There are many other interesting families of rational maps that have
been studied. For example, rational maps of degree~$2$ that
satisfy~$f(-x)=-f(x)$ form the $1$-dimensional family\footnote{Despite
  the two parameters~$b$ and~$c$ appearing in $bx+cx^{-1}$,
  geometrically this is actually a 1-dimensional family, because over
  an algebraically closed field, the change of variable
  $c^{-1/2}{f(c^{1/2}x)}=bx+x^{-1}$ eliminates~$c$. In general, we
  call~$c$ a twist parameter, since the dynamical properties of the
  map depend only on~$c$ modulo squares.}  of maps $f(x)=bx+cx^{-1}$.

\begin{conjecture}
[Manes, 2008, \cite{MR2407816}]
\label{conjecture:perbdaut2}
A $\QQ$-rational periodic point of a map of the form $bx+cx^{-1}$ with
$b,c\in\QQ^*$ has period~$1$,~$2$, or~$4$.  Further,
\[
  \sup_{b,c\in\QQ^*} \#\PrePer\bigl(bx+cx^{-1},\PP^1(\QQ)\bigr) = 12.
\]  
\end{conjecture}

More generally, we may consider the $2$-dimensional family
consisting of all degree~$2$ rational maps of~$\PP^1$ defined
over~$\QQ$.

\begin{conjecture}
[Benedetto, Chen, Hyde, Kovacheva, White, 2014, {\cite{MR3277939}}; Blanc--Canci--Elkies, 2015, \cite{MR3431627}]  
\label{question:deg2numpreper}
Let $f\in\End_2^1(\QQ)$ be a rational map of degree~$2$ defined over~$\QQ$.
\begin{parts}
\Part{(a)}
Let $P\in\PP^1(\QQ)$ be periodic for~$f$. Then~$P$ has period at most~$7$.  
\Part{(b)}
Let $P\in\PP^1(\QQ)$ be preperiodic for~$f$. Then~$\#\Orbit_f(P)\le9$.
\Part{(c)}
The full set of rational preperiodic points satisfies $\#\PrePer\bigl(f,\PP^1(\QQ)\bigr) \le 14$.
\end{parts}
\end{conjecture}

We note that there are examples in~\cite{MR3277939,MR3431627} showing
that the values~$7$,~$9$, and~$14$ in
Conjecture~\ref{question:deg2numpreper} are optimal.

\begin{remark}
Fakhruddin~\cite{MR1836365} has shown that
Conjecture~\ref{conjecture:MSUBC}, which is a uniform boundedness
conjecture for projective space, implies an analogous folklore uniform
boundedness conjecture for torsion points on abelian varieties. This
latter conjecture generalizes theorems for elliptic curves due to
Mazur~\cite{mazur:uniformbound}, Kamienny~\cite{MR1172689}, and
Merel~\cite{merel:uniformbound}.
\end{remark}

\subsection{Generalization to Algebraic Varieties}

Poonen has asked whether a uniform bound, analogous to
Conjecture~\ref{conjecture:MSUBC}, is true for preperiodic points in
arbitrary families of maps. We restate a version of Poonen's question
as a conjecture.

\begin{conjecture}
[Poonen, 2012, {\cite[Question 3.2.]{arxiv1206.7104}}]  
\label{question:poonenunifbd}
Let~$K$ be a number field, let~$X$ and~$S$ be varieties defined
over~$K$, let $\pi:X\to S$ be a morphism, and let $f:X\to X$ be an
$S$-morphism, which means that $\pi\circ f=\pi$. For each point~$s\in
S(\Kbar)$, the map~$f$ restricted to the fiber $X_s:=\pi^{-1}(s)$
gives a morphism $f_s:X_s\to X_s$, and we view~$f:X\to X$ as a family
of maps parameterized by the points of~$S$. Then for all~$D\ge1$ we
have
\[
  \sup_{\substack{L/K\\ [L:K]=D\\}} \sup_{s\in S(L)}
  \#\PrePer\bigl(f_s,X_s(L)\bigr) < \infty.
\]
\end{conjecture}

Taking $X\to S$ to be the universal family of degree~$d$ endomorphisms
of~$\PP^N$, cf.\ Section~\ref{section:dms}, one sees that
Conjecture~\ref{question:poonenunifbd} implies the Uniform Boundedness
Conjecture (Conjecture~\ref{conjecture:MSUBC}).

%%%%%%%%%%%%%%%%%%%%%%%%%%%%%%%%%%%%%%%%%%%%%%%%%%%%%%%%%%%%%%%%%%%%%%
\section{Arboreal Representations\WhoDone{Rafe+Joe}}
\label{section:ar}
\WhoWrite{Rafe+Joe}
%%%%%%%%%%%%%%%%%%%%%%%%%%%%%%%%%%%%%%%%%%%%%%%%%%%%%%%%%%%%%%%%%%%%%%

Let~$K$ be a field of characteristic~$0$, let~$\Kbar$ be an algebraic
closure of~$K$, and let $f:\PP^N\to\PP^N$ be a morphism of
degree~$d\ge2$ defined over~$K$.\footnote{More generally, we could
  take~$f:V\to V$ to be a self-morphism of a quasi-projective variety,
  for example a polynomial map $f:\AA^N\to\AA^N$, subject to suitable
  assumptions.}  Rather than looking at forward orbits of points, in
this section we start with an initial point~$P \in \PP^N(K)$ and look
at the points in its \emph{backward orbit}
\[
  \Orbit_f^{\scriptscriptstyle-}(P) := \bigl\{ Q \in \PP^N(\Kbar): \text{$f^k(Q)=P$ for some $k\ge0$} \bigr\}.
\]
Writing
\[
  f^{-n}(P) := \{Q \in \PP^N(\Kbar): f^n(Q)=P \},
\]
we see that  the disjoint union
\[
T_{f,\infty}(P) :=  \{P\} \sqcup f^{-1}(P) \sqcup f^{-2}(P) \sqcup f^{-3}(P) \sqcup \cdots
\]
acquires the structure of a rooted tree, with root vertex $P$, where
we assign edges according to the action of~$f$.

\begin{definition}
We say that the tree $T_{f,\infty}(P)$ is \emph{complete} if
$\#f^{-n}(P)=d^n$ for every $n \geq 0$, i.e., if $T_{f,\infty}(P)$ is a
complete rooted $d$-ary tree.
\end{definition}

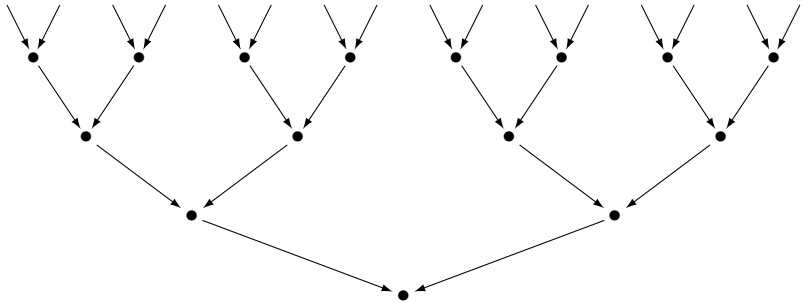
\begin{figure}[ht]
\begin{center}
\begin{picture}(260,120)(0,0)
  \multiput(0,90)(40,0){8}{\makebox(0,0)[c]{\textbullet}}   % top row of dots
  \multiput(20,60)(80,0){4}{\makebox(0,0)[c]{\textbullet}}  % second row of dots
  \multiput(60,30)(160,0){2}{\makebox(0,0)[c]{\textbullet}} % third row of dots
  \put(140,0){\makebox(0,0)[c]{\textbullet}}                % fourth row of dots (one dot)
  \multiput(2,87)(80,0){4}{\vector(2,-3){16}}               % top row of down-right arrows
  \multiput(38,87)(80,0){4}{\vector(-2,-3){16}}             % top row down-left arrows
  \multiput(24,57)(160,0){2}{\vector(4,-3){32}}             % second row down-right arrows
  \multiput(96,57)(160,0){2}{\vector(-4,-3){32}}            % second row down-left arrows
  \put(64,28.5){\vector(8,-3){72}}                          % bottom row down-right arrow
  \put(216,28.5){\vector(-8,-3){72}}                        % bottom row down-left arrow
  \multiput(-10,110)(40,0){8}{\vector(1,-2){8.5}}           % extra row of down-right arrows at top
  \multiput(10,110)(40,0){8}{\vector(-1,-2){8.5}}           % extra row of down-left arrows at top
\end{picture}
\end{center}
\caption{A complete binary inverse image tree }
\label{figure:2branchtree} 
\end{figure}

For example, if $N=1$, then $T_{f,\infty}(P)$ is complete if and only
if for inverse image $f^{-n}(P)$ contains no critical points of~$f$
for all $n\ge1$.  Figure~\ref{figure:2branchtree} illustrates a complete
inverse image tree for a degree~$2$ map~$f:\PP^1\to\PP^1$.

The coordinates of the points in~$\Orbit_f^{\scriptscriptstyle-}(P)$ generate finite
extension fields of~$K$. For each~$n\ge1$, we let
\[
  K_{f,n}(P) := \text{the extension of $K$ generated by the coordinates of
    points in $f^{-n}(P)$,}
\]
and we let
\[
  K_{f,\infty}(P) := \bigcup_{n\ge0} K_{f,n}(P)
\]
be the field generated by the coordinates of all of the points in the
full backward orbit of~$P$. In the case that~$K$ is a number field, we
can use the theory of height functions as described in
Section~\ref{section:ado} to show that~$T_{f,\infty}(P)$ is a set of 
points of bounded height, which in turn implies that
$\bigl[K_{f,\infty}(P):K\bigr]=\infty$ provided that the tree has
infinitely many distinct points.

For the rest of this section, we make the following assumption:
\[
  \framebox{The tree  $T_{f,\infty}(P)$ is complete.}
\]

Then the Galois group $G_K:=\Gal(\Kbar/K)$ acts on $f^{-n}(P)$ through
its action on the coordinates of each point. Moreover, this action
commutes with $f$, since~$f$ is defined over $K$. We let\footnote{We
  remark that the automorphism group of an $n$-level $d$-branched tree
  may be identified with the $n$-fold wreath product of the symmetric
  product~$\Scal_d$, so $\Aut\bigl(T_\infty(P)\bigr)$ is isomorphic to
  the inverse limit of $\Scal_d\wr\Scal_d\wr\cdots\wr\Scal_d$ ($n$
  copies) as $n\to\infty$.}
\[
  \Aut\bigl(T_{f,\infty}(P)\bigr) = \text{the group of tree automorphisms of $T_{f,\infty}(P)$}.
\]
The action of~$G_K$ on the points in~$T_{f,\infty}(P)$ yields a
homomorphism
\[
  \rho^\arb_{f,P}: G_K \to \Aut\bigl(T_{f,\infty}(P)\bigr).
\]
The map $\rho^\arb_{f,P}$ is called the \emph{arboreal representation
  associated to the map~$f$ and basepoint~$P$},\footnote{The term
  \emph{arboreal representation} for homomorphisms of Galois groups to 
  automorphism groups of trees was coined by Nigel Boston around~2003.}  and we
denote its image by
\[
  G_{f,\infty}(P) := \rho^\arb_{f,P}(G_K) \cong \Gal\bigl(K_{f,\infty}(P)/K\bigr).
\]
One can also realize $G_{f,\infty}(P)$ by taking the inverse limit of the
groups $G_{f,n}(P) := \Gal\bigl(K_{f,n}(P)/K\bigr)$ with respect to the natural
restriction maps.

Characterizing the size of the image of $\rho^\arb_{f,P}$ is one of the
main problems in this area, inspired in part by finite-index results
for Galois representations associated to elliptic curves due to
Serre~\cite{MR0387283,MR1484415} and analogous conjectures for abelian
varieties.

\begin{definition}
The \emph{Odoni index} of the map~$f$ and basepoint~$P$ over the
field~$K$ is the quantity
\[
  \iota_K(f,P) := \bigl[\Aut\bigl(T_{f,\infty}(P)\bigr) : G_{f,\infty}(P)\bigr].
\]
\end{definition}

\begin{question}
\label{arborealquestion}
\begin{parts}
\Part{(a)}
For which choices of fields~$K$, maps $f$, and basepoints $P$ is the
Odoni index~$\iota_K(f,P)$ equal to~$1$, i.e., when is the arboreal
representation~$\rho^\arb_{f,P}$ surjective?
\Part{(b)}
For which choices of fields~$K$, maps $f$, and basepoints $P$ is the
Odoni index~$\iota_K(f,P)$ finite, i.e., when does the image of the
arboreal representation~$\rho^\arb_{f,P}$ have finite index in the full tree automorphism group?
\end{parts}
\end{question}

Virtually nothing is known for $N\ge2$, so we restrict attention to
maps $f:\PP^1\to\PP^1$. In this case, Odoni~\cite{MR805714} showed
in~1985 that if the coefficients of a monic polynomial~$f(X)$ are
independent indeterminates, then $\iota_K(f,0)=1$, and he
also~\cite{MR813379} gave the explicit example $\iota_\QQ(X^2-X+1,0)=1$.

\begin{conjecture}[Odoni, 1985, \cite{MR805714}]
\label{conjecture:odoni}  
Let~$K$ be a Hilbertian field,\footnote{A field~$K$ is
  \emph{Hilbertian} if it is not possible to cover~$\PP^1(K)$ by a
  union of finitely manys sets of the form $F\bigl(C(K)\bigr)$,
  where~$C/K$ is a smooth curve and~$F:C\to\PP^1$ is a finite map of
  degree at least~$2$.  For example, all number fields are Hilbertian,
  but~$\RR$,~$\CC$, and~$\QQ_p$ are not.}  and let $d\ge2$. There
exists a monic polynomial~$f(X)\in K[X]$ of degree~$d$ and a
point~$P\in K$ such that $\iota_K(f,P)=1$.  \textup(Addendum:
\cite[Jones, 2015]{MR3220023} If~$K$ is a number field, then one can
take the coefficients of~$f(X)$ to lie in the ring of integers
of~$K$.\textup)
\end{conjecture}

Looper~\cite{arxiv1609.03398} recently proved Odoni's conjecture
(Conjecture~\ref{conjecture:odoni}) for~$K=\QQ$ and all prime
degrees~$d$, and a number of
researchers~\cite{arxiv1803.01987,arxiv1802.09074,arxiv1803.00434}
independently extended Looper's work.  In particular,
Specter~\cite{arxiv1803.00434} proved Odoni's conjecture for all
number fields, and more generally, for all algebraic extensions~$K/\QQ$
that are unramified outside of a finite set of primes.

We note that for any~$f$ defined over~$K$, there is always a finite
extension field~$L/K$ such that~$\iota_L(f,P)\ge2$.  For example, we
can simply adjoin a point in~$f^{-1}(P)$ to~$K$.  So we turn our
attention to the fundamental question of when~$\iota_K(f,P)$ is
finite.

Over number fields~$K$, it is known that the Odoni
index~$\iota_K(f,P)$ is finite for certain families of quadratic
polynomials~\cite{MR2439638,MR1174401}, and there are isolated
examples of quadratic rational functions~\cite{MR3177912} and
higher-degree polynomials~\cite{arxiv1609.03398}, as well as examples
over global function fields~\cite{arxiv1710.04332}.
See~\cite{MR3220023} for a survey of known results.

On the other hand, there are several conditions that are known to ensure that
$\iota_K(f,P)=\infty$. These include the map being PCF,\footnote{We
  recall that a map is \emph{post-critically finite} (PCF)
  if its critical points are preperiodic.}  maps for which there is a
$K$-automorphism of~$\PP^1$ commuting with~$f$ and fixing $P$, using
a basepoint~$P$ that is a periodic point of~$f$, and choosing a map~$f$
such that the forward orbits of its critical points intersect in
certain ways.  In the case~$N = 1$ and~$d = 2$, we now have enough
knowledge to plausibly conjecture that these are the only exceptions
to finite index.

\begin{conjecture}[{\cite[Conjecture 3.11]{MR3220023}}]
\label{quadratic}
Let $K$ be a field of characteristic~$0$, and suppose that $f \in K(x)$
is a rational function of degree~$2$. Then
\[
  \iota_K(f,P) = \infty
\]
if and only if one of the following holds\textup:
\begin{parts}
\Part{(a)}
The map $f$ is PCF, i.e., the critical points of~$f$ are preperiodic.
\Part{(b)}
There is an $r\ge2$ so that the two critical points $\gamma_1$ and $\gamma_2$ of $f$ satisfy a
relation of the form $f^{r}(\gamma_1) = f^{r}(\gamma_2)$.
\Part{(c)}
The basepoint $P$ of the tree $T_{f,\infty}(P)$ is periodic under $f$.
\Part{(d)}
There is a non-trivial M\"obius transformation
$\f(x)=\frac{\a{x}+\b}{\g{x}+\d}$ defined over $K$ that
satisfies\footnote{See Section~\ref{section:dms} for a detailed
  discussion of $\Aut(f)$, which is the group of automorphisms
  of~$\PP^N$ that commute with~$f$.}
\[
  \f\circ f(x) = f\circ\f(x)\quad\text{and}\quad \f(P)=P.
\]
\end{parts} 
\end{conjecture}

%%%%%%%%%%%%%%%%%%%%%%%%%%%%%%%%%%%%%%%%%%%%%%%%%%%%%%%%%%%%%%%%%%%%%%
\section{Dynatomic Representations\WhoDone{Joe}}
\label{section:dr}
\WhoWrite{Joe}
%%%%%%%%%%%%%%%%%%%%%%%%%%%%%%%%%%%%%%%%%%%%%%%%%%%%%%%%%%%%%%%%%%%%%%

\begin{figure}[ht]
\def\SixCycle{
\put(12.5,0){\circle*{4}}    
\put(37.5,0){\circle*{4}}    
\put(50,21.65){\circle*{4}}  
\put(37.5,43.3){\circle*{4}} 
\put(12.5,43.3){\circle*{4}} 
\put(0,21.65){\circle*{4}}   
\put(16,0){\vector(1,0){19}}  %% x + 3.5, L - 6
\put(39.5,2.5){\vector(100,173){9.5}} %%
\put(48,25){\vector(-100,173){9.5}}  %% 2,3.4
\put(34,43.3){\vector(-1,0){19}} %% 
\put(10,40.5){\vector(-100,-173){9.5}} %%
\put(2,18.25){\vector(100,-173){9.5}}
}
%\begin{center}
\begin{picture}(320,50)(5,30)
  \put(0,30){\SixCycle}
  \put(90,30){\SixCycle}
  \put(180,30){\SixCycle}
  \put(250,50){\makebox(0,0)[l]{$\boldsymbol{\cdots}$}}
  \put(280,30){\SixCycle}
\end{picture}
%\end{center}
\caption{A set of points of period $6$}
\label{figure:6cycles}
\end{figure}
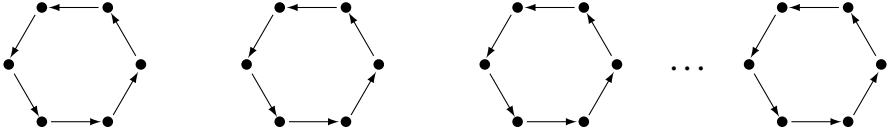

Let~$K$ be a field of characteristic~$0$ with algebraic
closure~$\Kbar$, and let $f:\PP^N\to\PP^N$ be a morphism of
degree~$d\ge2$ defined over~$K$.  We recall that for each~$n\ge1$, we
write
\[
\Per_n^*\bigl(f,\PP^N(\Kbar)\bigr)
:=
\bigl\{ P\in\PP^N(\Kbar):\text{$f^n(P)=P$ and $f^k(P)\ne P$ for $1\le k<n$}\bigr\}
\]
for the set of points of exact period~$n$. We note that~$f$
maps this set to itself, and that the resulting graph can be represented as a union of
$n$-gons.  Figure~\ref{figure:6cycles}
illustrates~$\Per_6^*(f)$, which is a collection of periodic cycles of
period~$6$.

%%%%%%%%%%%%%%%%%%%%%%%%%%%%%%%%%%%%%%%%%%%%%%%%%%%%%%%%%%%%%%%%%%%%%%
%%%%%%%%%%%%%%%%%%%%%%%%%%%%%%%%%%%%%%%%%%%%%%%%%%%%%%%%%%%%%%%%%%%%%%

An automorphism of a collection of $n$-cycles consists of a
permutation of the cycles and a rotation of each cycle.
For $\nu\ge1$, we set the notation
\[
  \Pcal_{n,\nu}^* := \text{a disjoint union of~$\nu$ periodic~$n$-cycles}.
\]
Then the automorphism group of~$\Pcal_{n,\nu}^*$ is most easily
described as a \emph{wreath product}, which by definition is a
semi-direct product,
\[
  \Aut(\Pcal_{n,\nu}^*) \cong
  (\ZZ/n\ZZ) \wr \Scal_\nu := (\ZZ/n\ZZ)^\nu \rtimes \Scal_\nu.
\]
The group law on this semi-direct product is given by
\[
  \bigl( (c_1,\ldots,c_\nu) , \pi \bigr) \cdot
  \bigl( (c'_1,\ldots,c'_\nu) , \pi' \bigr)
  := \bigl( (c_1+c'_{\pi(1)},\ldots,c_\nu+ c'_{\pi(\nu)}) , \pi\pi' \bigr),
\]
where~$c_i$ describes a rotation of the $i$th cycle and~$\pi$
permutes the cycles. In particular, the automorphism group
of~$\Pcal_{n,\nu}^*$ has order
$\#\Aut(\Pcal_{n,\nu}^*)=n^\nu\cdot\nu!$.

The equation $f^n(P)=P$ is equivalent to the intersection of~$N$
hypersurfaces of degree~$d^n+1$ in~$\PP^N$, so by Bezout's theorem it
consists of~$(d+1)^N$ points, counted with multiplicities.  If all of
the points have multiplicity~$1$, then the number of points in
$\Per_n^*(f)$ is given by an inclusion/exclusion formula. One expects
this to be the case for sufficiently large~$n$, as in the following
conjecture.\footnote{Let $P\in\Per^*_n(f)$. Then~$P$ is fixed
  by~$f^n$, so~$f^n$ induces a linear map~$Df^n(P)$ on the tangent
  space of~$\PP^N$ at~$P$, and the periodic point multiplicity of~$P$
  is greater than~$1$ if and only if $Df^n(P)$ has an eigenvalue equal
  to~$1$. For example, the point~$\frac12$ is a fixed point of
  multiplicity~$2$ for the map $f(x)=x^2+\frac14$, since
  $f(x)-x=(x-\frac12)^2$, and we observe that $f'(\frac12)=1$.}
  
\begin{conjecture}
\label{conjecture:numberPernstar}
Let~$K$ be a field of characteristic~$0$,
and let $f:\PP^N\to\PP^N$ be a morphism of degree~$d\ge2$ defined
over~$K$.  Then for all sufficiently large~$n$ we have\footnote{The
  function~$\mu$ is the M\"obius function, which we recall is defined
  as follows: If~$n$ is a product of~$k$ distinct primes,
  then~$\mu(n)=(-1)^k$, otherwise~$\mu(n)=0$.}
\[
  \#\Per_n^*\bigl(f,\PP^N(\Kbar)\bigr)
%% = \nu(n,N,d)
  := \sum_{m\mid n}  \mu\left(\frac{n}{m}\right) (d^m+1)^N.
\]
\end{conjecture}

We note that Conjecture~\ref{conjecture:numberPernstar} is true
for~$N=1$.  This follows from the fact that there are only finitely
many rationally indifferent periodic points, i.e., periodic points
whose multiplier is a root of unity. This in turn is a consequence of
a theorem of Fatou that a rational function~$f(x)\in\CC(x)$ of
degree~$d\ge2$ has at most~$6d-6$ non-repelling periodic
cycles~\cite[Theorem~13.1]{MR2193309}, a bound that was improved to
the best possible~$2d-2$ by Shishikura~\cite{MR892140}.

Since we have assumed that~$f$ is defined over~$K$, we see that the
action of the Galois group~$G_K:=\Gal(\Kbar/K)$ commutes with~$f$,
i.e.,
\[
  \s\bigl(f(P)\bigr) = f\bigl(\s(P)\bigr)
  \quad\text{for all $P\in\PP^N(\Kbar)$ and all $\s\in G_K$.}
\]
It follows that~$G_K$ acts on the set
$\Per_n^*\bigl(f,\PP^N(\Kbar)\bigr)$, and that this action preserves
the cycle structure induced by~$f$.  Hence the action of~$G_K$
on~$\Per_n^*(f)$ yields a homomorphism
\[
  \rho_{n,f}^\dyn : G_K \longrightarrow \Aut(\Pcal_{n,\nu}^*).
\]
The map~$\rho_{n,f}^\dyn$ is called the \emph{$n$-level dynatomic
  representation associated to~$f$}. It is a fundamental problem to
determine the size of its image. For example, does the image have
index bounded independent of~$n$?  And is~$\rho_{n,f}^\dyn$ surjective
for all sufficiently large~$n$?

An issue regarding this last question arises if the map~$f$ admits
automorphisms that are defined over~$K$. We define\footnote{See
  Section~\ref{section:dms} for a detailed discussion of $\Aut(f)$.}
\[
  \Aut_K(f) := \bigl\{ \f\in\PGL_{N+1}(K) : \f\circ f = f\circ\f \bigr\}.
\]
Elements of~$\Aut_K(f)$ act on $\Per_n^*(f)$ and commute with both~$f$
and the action of~$G_K$. This puts a further constraint on the
image of~$\rho_{n,f}^\dyn$, which leads us to define
\[
  \Aut(\Pcal_{n,\nu}^*)^{\Aut_K(f)}
  := \bigl\{ \a\in \Aut(\Pcal_{n,\nu}^*) : \text{$\f\circ\a=\a\circ\f$ for all $\a\in\Aut_K(f)$} \bigr\}.
\]

\begin{question}
\label{question:dynrepn}
Let~$K$ be a field of characteristic~$0$, let $f:\PP^N\to\PP^N$ be a
morphism of degree~$d\ge2$ defined over~$K$, and for~$n\ge1$, let
$\nu=\nu(n,f)=\#\Per_n^*(f)$. Find general necessary and
sufficient conditions on~$f$ for the following statements to
hold\textup:
\begin{parts}
\Part{(a)}
The index of the dynatomic representation $\rho_{n,f}^\dyn$ is bounded, i.e., 
\[
  \sup_{n\ge1}\; \bigl ( \Aut(\Pcal_{n,\nu}^*) : \rho_{n,f}^\dyn(G_K) \bigr) < \infty.
\]
\Part{(b)}
For all sufficiently large~$n$, the dynatomic representation
$\rho_{n,f}^\dyn$ is as surjective as possible, i.e.,
\[
  \rho_{n,f}^\dyn(G_K) = \Aut(\Pcal_{n,\nu}^*)^{\Aut_K(f)}.
\]
\end{parts}
\end{question}

Very little is known about Question~\ref{question:dynrepn}, although
it is not hard to find~$f$ for which the assertions are probably
false.

\begin{example}
Let~$p=2^\ell-1$ be a Mersenne prime. Then $\Per_\ell^*(x^2)$ is the
set of primitive $p$th roots of unity~$\bfmu_p^*$, which the
map~$x\to x^2$ separates into~$p/\ell$ distinct~$\ell$-cycles. The
action of $\Gal(\Qbar/\QQ)$ on~$\bfmu_p^*$ is completely determined by
where it sends a single element of~$\bfmu_p^*$, so the Galois action cannot
come close to being surjective. More precisely, we have
\[
  \#\rho_{\ell,x^2}^\dyn(G_\QQ) = p-1,\quad
  \Aut_\QQ(f)=\{x,x^{-1}\},\quad
  \#\Aut(\Pcal_{\ell,p/\ell}^*)^{\Aut_\QQ(f)} = \ell^{p/\ell}\cdot (p/\ell)!.
\]
\end{example}
%% Note that the map x -> 1/x commutes with the Galois maps x -> x^i 

%%%%%%%%%%%%%%%%%%%%%%%%%%%%%%%%%%%%%%%%%%%%%%%%%%%%%%%%%%%%%%%%%%%%%%
\section{Intersections of Orbits and Subvarieties\WhoDone{Tom}}
\label{section:ios}
\WhoWrite{Tom (started by Joe)}
%%%%%%%%%%%%%%%%%%%%%%%%%%%%%%%%%%%%%%%%%%%%%%%%%%%%%%%%%%%%%%%%%%%%%%

Let~$f:X\to X$ be a morphism, let~$Y\subseteq X$ be a subvariety
of~$X$, and let~$P\in X$ be a wandering point for~$f$. It is natural
to ask under what circumstances the orbit~$\Orbit_f(P)$  intersects
the subvariety~$Y$ at infintely many points. For example, one way that this
can happen is if $\Orbit_f(P)\cap Y$ is non-empty and~$Y$ contains a
subvariety~$Z\subset Y$ of dimension at least~$1$ having the property
that~$f^n(Z)\subset Z$ for some integer~$n\ge1$.

In order to state a precise conjecture, we define a (\emph{one-sided})
\emph{arithmetic progression} in~$\NN$ to be a collection of integers
of the form
\[
  \{ a k + b : k = 0,1,2,\ldots \}
\]
for some fixed integers~$a\ge0$ and~$b\ge0$. Note that we
permit~$a=0$, so a single-element set~$\{b\}$ is an arithmetic
progression.

\begin{conjecture}
[Dynamical Mordell--Lang Conjecture:
    Bell, 2006, \cite{MR2225492};
    Denis, 1994, \cite{MR1259107};
    Ghioca--Tucker, 2009, \cite{arxiv0805.1560}]
\label{conjecture:dynmordelllang}
Let~$X/\CC$ be a quasi-projective variety, let $f:X\to X$ be a morphism,
let~$P\in X$ be a point, and let~$Y\subseteq X$ be a subvariety of~$X$.
Then
\[
  \{ n\ge 0 : f^n(P)\in Y \}
\]
is a finite union of arithmetic progressions.
\end{conjecture}

Before stating a consequence of the conjecture, we generalize the
notion of (pre)periodic point to (pre)periodic subvariety.

\begin{definition}
\label{definition:preperiodicvariety}
Let~$f:X\to X$ be an endomorphism of a quasi-projective variety.  A
subvariety~$Y\subset X$ is \emph{periodic} if~$f^n(Y)\subseteq Y$ for some~$n\ge1$,
and~$Y$ is \emph{preperiodic} if~$f^n(Y)\subseteq f^m(Y)$ for some~$n>m\ge0$.
\end{definition}

Conjecture \ref{conjecture:dynmordelllang} implies that if~$P$ is
wandering and~\text{$\Orbit_f(P) \cap Y$} is infinite, then the
Zariski closure of~\text{$\Orbit_f(P) \cap Y$}  contains a periodic
subvariety~$Z$ of dimension at least~$1$.  In particular, if~$Y$
contains no $f$-periodic subvarieties of dimension at least~$1$, then
conjecturally the intersection~\text{$\Orbit_f(P)\cap Y$} is finite.

There is an extensive literature proving various special cases of
Conjecture~\ref{conjecture:dynmordelllang}; see for
example~\cite{MR3468757,arXiv:0712.2344, arxiv1512.03085,
   arXiv:0704.1333, arxiv0805.1560, MR3235993, arxiv1610.00367,
  MR2861079, arxiv0705.1954, arxiv1604.08216, MR3360335,
  arxiv1109.0207, arxiv1403.3975, MR3263169, arxiv1503.00773,
  arxiv1602.04253}.  Here are two notable examples.

\begin{theorem}
Conjecture~$\ref{conjecture:dynmordelllang}$ is true in the following
cases\textup:
\begin{parts}
  \Part{(a)}
  \textup{(Xie, 2015, \cite{arxiv1503.00773})}:
   The map ~$f: \AA^2 \to \AA^2$ is a morphism defined over~$\Qbar$.
   \Part{(b)}
   \textup{(Bell--Ghioca--Tucker, 2010, \cite{arxiv:0808.3266})}
   The map~$f: X \to X$ is unramified.  Note that this includes the case
   that~$f$ is an automorphism, e.g., a polynomial automorphism of~$\AA^N$.
\end{parts}
\end{theorem}
  
There are also variants of Conjecture~\ref{conjecture:dynmordelllang}
for monoids generated by more than one morphism, as in the following
result.

\begin{theorem}
[Ghioca--Tucker--Zieve, 2016, \cite{MR3468757}] Let
$f_1(x),f_2(x)\in\CC[x]$ be polynomials of degree greater than one.
Let~$\alpha \in \CC$ be a wandering point of~$f_1$, and let~$\beta \in
\CC$ be a wandering point of~$f_2$.  Suppose that there are infinitely
many pairs~$(n_1,n_2) \in \NN^2$ such that~$f_1^{n_1}(\alpha) =
f_2^{n_2}(\beta)$.  Then there are positive integers~$\ell_1$
and~$\ell_2$ such that~$f_1^{\ell_1}(x) = f_2^{\ell _2}(x)$.
\end{theorem}

On the other hand, if one takes finitely many commuting morphisms
\[
 f_1, \dots, f_g: X \longrightarrow X,
\]
then it is impossible to say much of anything about the structure of the set
\begin{equation}
  \label{eqn:n1ngNg}
  \bigl\{ (n_1, \dots, n_g)\in\NN^g : f_1^{n_1} \circ \dots \circ f_g^{n_g} (P) \in Y \bigr\}.
\end{equation}
Indeed, Scanlon and Yasufuku \cite{MR3250036} (see also
\cite{MR2861079}) have shown that the set~\eqref{eqn:n1ngNg} can be
arbitrarily complicated, in a way that can be made precise using the
language of polynomial-exponential equations.

Conjecture~\ref{conjecture:dynmordelllang} says that if a point~$P$ is
not contained in a proper~$f$-preperiodic subvariety, then the
intersection \text{$\Orbit_f(P)\cap Y$} with any
subvariety~$Y\subsetneq X$ is finite.  This raises the question of
whether there exist \emph{any} algebraic points whose orbit does not
lie in a proper subvariety, i.e., whose orbit is Zariski dense.  An
obvious obstruction to the existence of such points would be a
fibering as in the following definition.

\begin{definition}
Let~$X$ be a quasi-projective variety, and let~$f:X\to X$ be an
endomorphism. We say that~$f$ is a \emph{fibered map} if there exists
a positive dimensional variety~$Z$ and a dominant rational map
$\psi:X\to Z$ such that~$\psi\circ f = \psi$.
\end{definition}

The following general conjecture, which proposes that
fibering is the only possible obstruction, builds on an earlier
conjecture due to Zhang~\cite{zhang:distalgdyn}.

\begin{conjecture}
[Amerik--Bogomolov--Rovinsky, 2011, \cite{MR2862064};
  Medvedev--Scanlon, 2009, \cite{arxiv0901.2352}]  
\label{newconjecture}
Let~$X$ be a quasi-projective variety defined over an algebraically
closed field~$K$ of characteristic~$0$. Let~$f:X\to X$ be an
endomorphism defined over~$K$, and assume that~$f$ is not a fibered
map.
\begin{parts}
\Part{(a)} \textup{(Weak form)}
There exists a point~$P\in X(K)$ whose orbit~$\Orbit_f(P)$
is Zariski dense in~$X$.
\Part{(b)} \textup{(Strong form)}
The set
\[
  \bigl\{ P\in X(K) : \text{$\Orbit_f(P)$ is Zariski dense in $X$} \bigr\}
\]
is Zariski dense in~$X$.
\end{parts}
\end{conjecture}

Conjecture~\ref{newconjecture}(a) has been proven in a handful of
cases~\cite{MR2862064, MR2784670, arxiv1310.5775, arxiv0901.2352}, but
the general case remains very much open, even in the case that~$f$ is
an automorphism.

\begin{remark}  
The classical Mordell--Lang conjecture was proven by
Faltings~\cite{MR718935,faltings:mordelllangconj}, and the rank~$1$ case
fits into the set-up of this section.  Let~$A/\CC$ be
an abelian variety, let~$Q\in A(\CC)$ be a point, and let $T_Q:A\to A$
be the translation-by-$Q$ map defined by $T_Q(P)=P+Q$.
Let~$Y\subseteq A$ be a subvariety of~$A$. Then Faltings' theorem
implies that Conjecture~\ref{conjecture:dynmordelllang} is true for
the map~$T_Q$, i.e., the set of~$n\ge0$ such that~$T_Q^n(P)=nP+Q\in Y$
is a finite union of arithmetic progressions. The general case of the
classical Mordell--Lang conjecture has a similar formulation using the
monoid of maps generated by a finite list~$T_{Q_1},\ldots,T_{Q_r}$ of
translations.\footnote{The usual statement of the classical
  Mordell--Lang conjecture says that for any finitely generated
  subgroup $\G\subset A(\CC)$, the Zariski closure of~\text{$Y\cap\G$}
  in~$A$ is a finite union of translates of abelian subvarieties
  of~$A$ by points of finite order.}
\end{remark}

%%%%%%%%%%%%%%%%%%%%%%%%%%%%%%%%%%%%%%%%%%%%%%%%%%%%%%%%%%%%%%%%%%%%%%
\section{(Pre)Periodic Points on Subvarieties\WhoDone{Tom}}
\label{section:ppsv}
\WhoWrite{Tom (started by Joe)}
%%%%%%%%%%%%%%%%%%%%%%%%%%%%%%%%%%%%%%%%%%%%%%%%%%%%%%%%%%%%%%%%%%%%%%

As we have seen in Sections~\ref{section:ub} and~\ref{section:dr}, and
will further explore in Sections~\ref{section:dmc}
and~\ref{section:goodred}, the periodic and preperiodic points of a
map form a collection of specially marked points for the map, just as
the torsion points on an abelian variety are special.  The classical
Manin--Mumford conjecture, which was proven
Raynaud~\cite{MR688265,MR717600}, describes the distribution of
torsion points of an abelian variety~$A$ lying on a subvariety
of~$A$.\footnote{More precisely, the classical Manin--Mumford
  conjecture that was proven by Raynaud says that if~$A/\CC$ is an
  abelian variety and~$Y\subset A$ is a subvariety, then the Zariski
  closure of~\text{$Y\cap A_\tors$} is a finite union of translates of
  abelian subvarieties of~$A$ by points of finite order.}  Replacing
torsion points with preperiodic points,
cf.\ Example~\ref{example:Gprepertorsion}, leads to various dynamical
analogues.

Zhang~\cite{zhang:distalgdyn,MR1311351} formulated a version for
\emph{polarizable endomorphisms},\footnote{Formally, $f:X\to X$
  is \emph{polarizable} if there is an ample line bundle~$\Lcal$ on~$X$ and
  an integer $d\ge2$ such that $f^*\Lcal\cong\Lcal^{\otimes d}$.}
which are endomorphisms $f: X \to X$ of a projective variety that
extend to non-invertible endomorphisms of the projective space in
which~$X$ sits; see~\cite{MR1995861}. Unfortunately, this natural
conjecture turns out to be false, but we give the statement anyway
because it illustrates some of the subtleties and pitfalls that can
arise when generalizing conjectures and theorems from arithmetic
geometry to the dynamical setting.

\begin{conjecture}
[{\bfseries False} Dynamical Manin--Mumford Conjecture]
\label{DMM}
Let $X/\CC$ be a projective variety, let~$f:X\to X$ be a polarizable
morphism, and let~$Y$ be a subvariety of~$X$. If
\[
\# \bigl( Y \cap \PrePer(f) \bigr) = \infty,
\]
then~$Y$ contains a periodic subvariety of~$X$ of dimension at
least~$1$.  \textup(See
Definition~$\ref{definition:preperiodicvariety}$ for the definition of
periodic subvariety.\textup)
\end{conjecture}

Conjecture~\ref{DMM}, although plausible, is false.  One can construct
a counterexample by starting with an elliptic curve~$E$ having complex
multiplication and taking~$Y$ to be the diagonal in~$E \times E$; see
\cite{MR2854724}.  However, the following modification may be true.

\begin{conjecture}
[Ghioca--Tucker--Zhang, 2011, {\cite[Conjecture 2.4]{MR2854724}}]
\label{reformulation}
Let~$X/\CC$ be a projective variety, let~$f: X \to X$ be a polarizable
endomorphism, and let~$Y$ be a
subvariety of~$X$ which has no component contained entirely within the singular
locus of~$X$. Then the following are equivalent\textup:
\begin{parts}
  \Part{(a)}
  $Y$ is preperiodic under~$f$.
  \Part{(b)}
  The intersection $Y\cap \PrePer_{f}(X)$ contains a Zariski dense
  subset of smooth points such that for every~$x$ in the subset, the
  tangent subspace of~$Y$ at~$x$ is preperiodic under the induced
  action of~$f$ on the
  Grassmanian~$\operatorname{Gr}_{\dim(Y)}\left(T_{X,x}\right)$, where
  we write~$T_{X,x}$ for the tangent space of~$X$ at the point~$x$.
\end{parts}
\end{conjecture}

Conjecture~\ref{reformulation} has been proven in a few
special cases, including the following.

\begin{theorem}
  \label{theorem:MMtrue}
  Conjecture~$\ref{reformulation}$ is true in the following cases\textup:
  \begin{parts}
    \Part{(a)}
    \textup{(Ghioca--Tucker--Zhang, 2011, {\cite{MR2854724}})}\\
    $f$ is a group endomorphism of an abelian variety.
    \Part{(b)}
    \textup{(Ghioca--Nguyen-Ye, 2015, \cite{arxiv1511.06081};
      Ghioca--Nguyen--Ye, 2017, \cite{arxiv:1705.04873})}\\
    $f$ is an endomorphism of $(\PP^1)^n$.
  \end{parts}
\end{theorem}

We remark that the original Conjecture~\ref{DMM} is false for both of
the cases covered by Theorem~\ref{theorem:MMtrue}.

%%%%%%%%%%%%%%%%%%%%%%%%%%%%%%%%%%%%%%%%%%%%%%%%%%%%%%%%%%%%%%%%%%%%%%
\section{Dynamical (Dynatomic) Modular Curves \WhoDone{Michelle+Joe}}
\label{section:dmc}
\WhoWrite{Michelle + Joe}
%%%%%%%%%%%%%%%%%%%%%%%%%%%%%%%%%%%%%%%%%%%%%%%%%%%%%%%%%%%%%%%%%%%%%%
A natural way to study periodic points of rational maps on~$\PP^1$ is
to classify pairs~$(f,P)$, where~$f$ runs over all rational maps of
some fixed degree and~$P$ runs over all points of some fixed period
for~$f$.  A degree~$d$ rational map $f:\PP^1\to\PP^1$ is described by
a pair of homogeneous polynomials $f(X,Y)=[F(X,Y),G(X,Y)]$. We view
the coefficients of~$F$ and~$G$ as indeterminates, say
\[
  F(X,Y) = a_0X^d + a_1X^{d-1}Y+\cdots+a_d
  \quad\text{and}\quad
  G(X,Y) = b_0X^d + b_1X^{d-1}Y+\cdots+b_d,
\]
and we write the $n$th iterate of~$f$ as
\[
  f^n(X,Y) = [F_n(X,Y),G_n(X,Y)],
\]
so~$F_n$ and~$G_n$ are in the polynomial ring
$\ZZ[a_0,\ldots,a_d,b_0,\ldots,b_d,X,Y]$. The points whose period
divide~$n$ are then the roots of the polynomial
\[
  YF_n(X,Y) - XG_n(X,Y) = 0,
\]
but if we want to get the points of exact period~$n$, we need to
remove points of lower period. This is done by inclusion-exclusion,
modeled after the construction of the classical cyclotomic
polynomials.

\begin{definition}
With notation as above, the \emph{$n$th dynatomic polynomial} is
given by\footnote{Alternatively, we can
set~$\F_1=YF-XG$ and  define~$\F_n$ recursively via the
relation $\prod_{d\mid n}\F_d=YF_n-XG_n$.}
\[
\F_n(X,Y) := \prod_{d\mid n} \bigl( YF_n(X,Y) - XG_n(X,Y) \bigr)^{\mu(n/d)},
\]
where $\mu$ is the M\"obius function; see Section~\ref{section:dr}.
If we want to specify the dependence on~$f$, we will
write~$\F_n(f;X,Y)$. More generally, for any $m\ge0$ and $n\ge1$, we
characterize preperiodic points with tail length~$m$ and cycle length~$n$
using the \emph{generalized $(m,n)$ dynatomic polynomial}
\[
  \F_{m,n}(X,Y) := \frac{ \F_n\bigl( F_m(X,Y), G_m(X,Y) \bigr) }
  { \F_n\bigl( F_{m-1}(X,Y), G_{m-1}(X,Y) \bigr) }.
\]
(For~$m=0$, we set $\F_{0,n}:=\F_n$.)
\end{definition}

\begin{proposition}
\label{proposition:dynpolyispoly}
\textup{(\cite[Theorem~4.5]{MR2316407})}
The dynatomic polynomial is a polynomial. More precisely,
\[
  \F_n(X,Y) \in \ZZ[a_0,\ldots,a_n,b_0,\ldots,b_n,X,Y].
\]
\end{proposition}

For a given map~$f\in\End_d^1(\Kbar)$, the roots of its associated
dynatomic polynomial $\F_n(f;X,Y)$ include all of the points of exact
period~$n$ for~$f$, but it may happen that there are roots giving
points of period strictly smaller than~$n$. When this happens,~$\F_n$
has a multiple root. Following Milnor, the roots of~$\F_n(f;X,Y)$ are
said to be points of \emph{formal period~$n$ for~$f$}.

\begin{remark}
See~\cite{arxiv:0801.3643,arxiv:1011.5155} for a generalization of
Proposition~\ref{proposition:dynpolyispoly} to dynatomic $0$-cycles
on~$\PP^N$ and on other varieties.
\end{remark}

We now fix a degree $d\ge2$ and restrict attention to the special, but
enlightening, case of polynomials maps of the form\footnote{This is
  the dehomogenized form of the map $f_c(X,Y)=[X^d+cY^d,Y^d]$, where
  we set $X=x$ and $Y=1$.}
\[
  f_c(x) = x^d + c.
\]
The $n$th dynatomic polynomial for~$f_c(x)$ is a polynomial in the
variables~$x$ and~$c$, so we write it as $\F_n(x,c)\in\ZZ[x,c]$.  For
example, taking $d=2$ and $n=3$, we find that $\F_3(x,c)$ is given by
\begin{multline*}
 \F_3(x,c) = \smash[b]{\frac{f^3_c(x)-x}{f_c(x)-x}}
 = x^6 + x^5 - (3 c - 1) x^4 - (2 c - 1) x^3 \\*
 + (3 c^2 - 3 c + 1) x^2 + (c^2 - 2 c + 1) x - c^3 + 2 c^2 - c + 1.
\end{multline*}

Before stating the next result, we need one definition.

\begin{definition}
\label{definition:gonality}
Let $C$ be an irreducible algebraic curve.  The \emph{gonality of~$C$}
is the smallest degree of a dominant rational map $C\to\PP^1$. It is
an analogue for algebraic curves of the degree of a number field
over~$\QQ$. More generally, for any irreducible variety~$X$ of
dimension~$N$, we define the \emph{gonality of~$X$} to be the smallest
degree of a dominant rational map $X\to\PP^N$.
\end{definition}

\begin{theorem}
\label{theorem:X1nirredgon}
Let $d\ge2$, let~$\F_n(x,c)$ be the $n$th dynatomic polynomial
associated to the map $f_c(x)=x^d+c$, and let $X_1^\dyn(f_c;n)$ be the
smooth projective curve \textup(Riemann surface\textup) obtained by
starting with the affine curve defined by the equation $\F_n(x,z)=0$,
taking its closure in~$\PP^2$, and then resolving the
singularities. Further, let $X_0^\dyn(f_c;n)$ be the quotient of
$X_1^\dyn(f_c;n)$ by the equivalence relation
$(x,z)\sim\bigl(f_c(x),z\bigr)$.
\begin{parts}
\Part{(a)}  
\textup{(Bousch, 1992, \cite{bousch:thesis}; Lau--Schleicher, 1994,
  \cite{lauschleicher:irreducibility,MR3646529}; Morton, 1996,
  \cite{morton:dynamicalmoduli})} The polynomial~$\F_n(X,Y)$ does not
factor in~$\CC[X,Y]$, and hence the curves $X_1^\dyn(f_c;n)$
and~$X_0^\dyn(f_c;n)$ are geometrically irreducible.
\Part{(b)}  
\textup{(Morton, 1996, \cite{morton:dynamicalmoduli})}
There are explict, albeit quite complicated, formulas for the genera of
the curves~$X_1^\dyn(f_c;n)$ and~$X_0^\dyn(f_c;n)$.
\Part{(c)}  
\textup{(Doyle--Poonen, 2017,\cite{arxiv1711.04233})}
The gonalities of $X_1^\dyn(f_c;n)$ and $X_0^\dyn(f_c;n)$ satisfy
\begin{align*}
  \operatorname{Gonality}\bigl(X_1^\dyn(f_c;n)\bigr) &\xrightarrow[\;{n\to\infty}\;]{} \infty,\\*
  \operatorname{Gonality}\bigl(X_0^\dyn(f_c;n)\bigr) &\ge \left(\frac12-\frac1{2d}-o(1)\right)n
  \quad\text{as $n\to\infty$.}
\end{align*}
\end{parts}
\end{theorem}

\begin{remark}
The notation $X_1^\dyn(f_c;n)$ and $X_0^\dyn(f_c;n)$ is modeled after the
classical notation for the elliptic modular curves $X_1(n)$ and
$X_0(n)$.  The points of~$X_1(n)$ classify pairs $(E,P)$, with~$E$ an
elliptic curve and $P\in E$ a point of exact order~$n$. The
subscript~$1$ on~$X_1$ indicates that~$X_1(n)$ classifies a single
point of order~$n$, just as~$X_1^\dyn(f_c;n)$ classifies a single point of
period~$n$. More generally, there are
curves~$X_1^\dyn(f_c;n_1,\ldots,n_r)$ that classify maps with marked
points of several indicated periods, and even more generally, curves
$X_1^\dyn(\G)$ that classify maps with marked points and/or cycles
with orbits having a specified graph
structure. See~\cite{doyle17dynmod,moduliportrait2017} for
further details. 
\end{remark}

We next consider the set of parameter values~$c$ for which the
polynomial~$f_c(x)$ is post-critically finite (PCF). The critical
points of~$f_c(x)$ are~$\infty$, which is fixed, and~$0$, so~$f_c(x)$
is PCF if and only if there are integers~$m\ge0$ and~$n\ge1$ such that
$\F_{m,n}(f_c;0)=0$. The quantity~$\F_{m,n}(f_c;0)$ is a polynomial
in~$c$, and it is not hard to show that it has simple roots, but it
occasionally has an extraneous factor.

\begin{definition}
Let $d\ge2$ and $f_{d,c}(x)=x^d+c$. The \emph{$(m,n)$ Gleason
  polynomial} of~$f_{d,c}$ is
\[
G_{d,m,n}(c) := \begin{cases}
  \F_{m,n}(f_{d,c};0) / \F_{0,n}(f_{d,c};0)^{d-1}
  &\text{if $m\ge1$ and $n\bigm| (m-1)$,} \\
  \F_{m,n}(f_{d,c};0)
  &\text{otherwise.}\\
\end{cases}
\]
\end{definition}

The roots of $G_{d,m,n}(c)$ give the values of~$c$ for which the
non-trivial critical point of~$x^d+c$ has tail length~$m$ and formal
period length~$n$.

\begin{question}
\label{question:gleasonpoly}
Let $G_{d,m,n}(c)\in\ZZ[c]$ be the Gleason polynomial of~$x^d+c$.
\begin{parts}
\Part{(a)}
Is it true that $G_{2,m,n}$ is irreducible in~$\QQ[c]$ for all $m\ge0$
and $n\ge1$?
\Part{(b)}
For which~$(d,m,n)$ is $G_{d,m,n}$ irreducible in~$\QQ[c]$?
\Part{(c)}
More generally, what is the complete factorization of $G_{d,m,n}$
in~$\QQ[c]$?
\end{parts}
\end{question}

A few partial results are known regarding
Question~\ref{question:gleasonpoly}.  Goksel~\cite{arxiv1806.01208}
shows that for~$m\ge1$, the Gleason polynomial~$G_{d,m,n}$ is
irreducible in~$\QQ[c]$ if $n=1$ and~$d$ is prime, or if~$n=d=2$.
Buff~\cite{bufflpcfup} notes that~$G_{d,0,3}$ is reducible in~$\QQ[c]$
if and only if $d\equiv1\pmodintext6$.  See~\cite{buffepsteinkoch2018}
for generalizations of these results by Buff, Epstein, and Koch

It is natural to consider dynamical modular curves associated to other
interesting 1-parameter families of maps, for example the maps
\[
  f_b(x) = bx+x^{-1}
\]
discussed in Conjecture~\ref{conjecture:perbdaut2} that satisfy
$f_b(-x)=-f_b(x)$. The points on the dynatomic curve~$\F_n(x,b)=0$
associated to~$f_b$ classify points of formal period~$n$ for~$f_b$,
but in contrast to Theorem~\ref{theorem:X1nirredgon}, these curves may
be reducible.

\begin{theorem}
\label{theorem:manes2comp}
\textup{(Manes, 2009, \cite{MR2524186})} Let~$\F_n(f_b;x,b)$ be the
$n$th dynatomic polynomial associated to~$f_b(x)=bx+x^{-1}$, and let
$X_1^\dyn(f_b;n)$ be the smooth projective curve obtained by
completing and desingularizing the affine curve defined by the
equation $\F_n(f_b;x,z)=0$. Then for all $n\ge2$, the curve
$X_1^\dyn(f_b;2n)$ has at least two irreducible components, and in
fact, the polynomial $\F_n(f_b;x,z)=0$ has at least two irreducible
factors in~$\QQ[x,z]$.
\end{theorem}

\begin{conjecture}
\label{conjecture:manes2comp}
Let~$f_b(x)=bx+x^{-1}$ be as in Theorem~$\ref{theorem:manes2comp}$.
\begin{parts}
\Part{(a)}
If $n\ge1$ is odd, then the polynomial $\F_n(f_b;x,z)=0$ is
irreducible in~$\CC[x,z]$.
\Part{(b)}
If $n\ge4$ is even, then the polynomial $\F_n(f_b;x,z)=0$ has
exactly two irreducible factors in~$\CC[x,z]$.
\end{parts}
\end{conjecture}

See Conjecture~\ref{conjecture:compandgen}  for a super-sized, albeit less precise,
generalization of Conjecture~\ref{conjecture:manes2comp}.

%%%%%%%%%%%%%%%%%%%%%%%%%%%%%%%%%%%%%%%%%%%%%%%%%%%%%%%%%%%%%%%%%%%%%%
\section{Dynamical Moduli Spaces\WhoDone{Michelle+Joe}}
\label{section:dms}
\WhoWrite{Michelle + Joe}
%%%%%%%%%%%%%%%%%%%%%%%%%%%%%%%%%%%%%%%%%%%%%%%%%%%%%%%%%%%%%%%%%%%%%%

Rather than studying iteration of a single rational function~$f$, we
saw in Section~\ref{section:dmc} that it is often advantageous to
consider families of maps, or even the family of all maps. A rational
map $f:\PP^N\to\PP^N$ of degree~$d$ is specified by an~$(N+1)$-tuple
\[
  f = [f_0,\ldots,f_N]
\]
of degree~$d$ homogeneous polynomials
$f_0,\ldots,f_N\in\Kbar[X_0,\ldots,X_n]$, Each~$f_i$ has the form
\[
  f_i = \sum_{j_0+j_1+\cdots+j_N=d} a_{ij_0j_1\cdots j_N} X_0^{j_0}X_1^{j_1}\cdots X_N^{j_N}.
\]
We identify~$f$ with the list of
coefficients~$a_{ij_0\cdots{j_N}}$, but we note that multiplying the
coefficients by a non-zero constant gives the same map~$f$. Thus the
map~$f$ is identified with the point in projective space given by its
coefficients listed in some specified order,
\[
  f \longleftrightarrow [a_{ij_0j_1\cdots j_N}] \in \PP^\nu(\Kbar),
  \qquad\text{$\nu=\nu(N,d):= \dbinom{N+d}{d}(N+1)-1$.}
\]

Of course, not every point in~$\PP^\nu(\Kbar)$ gives a rational map of
degree~$d$, since the corresponding polynomials~$f_0,\ldots,f_N$ may
have a common factor, but ``most'' do, and indeed, most give a
morphism, as explained by the following result.

\begin{theorem}
The subset of~$\PP^\nu$ corresponding to degree~$d$ morphisms
of~$\PP^N$, which we denote by~$\End_d^N$, is a Zariski open subset
of~$\PP^\nu$.\footnote{More precisely, there is a single geometrically
  irreducible polynomial~$R(f)$ in the coefficients of~$f$, called the
  \emph{Macaulay resultant} of~$f$, having the property that $f\in\End_d^N$
  if and only if $R(f)\ne0$. For $N=1$, if we write $f=[F(X,Y),G(X,Y)]$,
  then $R(f)$ is just the classical resultant of~$F$
  and~$G$.}{$^,$}\footnote{It is also true that the subset of~$\PP^\nu$
  corresponding to degree~$d$ dominant rational maps of~$\PP^N$ is a
  Zariski open subset
  of~$\PP^\nu$; see~\cite[Proposition~7]{arxiv0908.3835}.}
\end{theorem}

If we start with a dominant rational map $f:\PP^N\dashrightarrow\PP^N$
and make a simultaneous change-of-variables to the domain and range,
we get a new map whose dynamical properties are the same as the
original one. A linear change of variables of~$\PP^N$ is given by an
invertible $(N+1)$-by-$(N+1)$ matrix, and two matrices determine the
same variable change if they are scalar multiples of one another. Thus
the automorphism group of~$\PP^N$ is the projective linear group
\[
  \Aut(\PP^N) = \PGL_{N+1}.
\]
We define an action of $\f\in\PGL_{N+1}(\Kbar)$ on a rational map
$f:\PP^N\to\PP^N$ by setting
\begin{equation}
  \label{eqn:conjaction}
  f^\f := \f^{-1}\circ f \circ \f.
\end{equation}
This action commutes with iteration in the sense that
\[
  (f^\f)^n = (f^n)^\f.
\]
In studying families of dynamical systems on~$\PP^N$, we really want
to study equivalence classes of maps modulo this action
of~$\PGL_{N+1}$. Identifying the set of maps~$f$ with (a subset of) the projective
space~$\PP^\nu$ as above, the conjugation action~\eqref{eqn:conjaction} defines
a complicated, but linear, action of~$\PGL_{N+1}$ on~$\PP^\nu$; more precisely,
it defines a group homomorphism
\[
  \PGL_{N+1} \longhookrightarrow \PGL_{\nu+1} = \Aut(\PP^\nu).
\]

\begin{definition}
The \emph{moduli space of degree $d$ dynamical systems on $\PP^N$} is
the quotient space
\[
  \Moduli_d^N := \End_d^N / \PGL_{N+1}.
\]
\end{definition}

This raises a delicate question: to what extent does the quotient
$\Moduli_d^N$ exist and have nice properties? This question is
partially answered by the following amalgam of results.

\begin{theorem}
Let $N\ge1$ and $d\ge2$.
\begin{parts}
\Part{(a)}
\textup{(Milnor \cite{milnor:quadraticmaps})}
$\Moduli_d^1(\CC)$ exists as a complex orbifold.
\Part{(b)}
\textup{(Silverman \cite{silverman:modulirationalmaps},
  Petsche--Szpiro--Tepper \cite{MR2567424}, Levy \cite{MR2741188})}
$\Moduli_d^N$ exists as a ``nice'' quotient variety. In technical
terms,~$\Moduli_d^N$ exists as a geometric quotient scheme over
$\Spec(\ZZ)$ in the sense of geometric invariant theory
\textup(GIT\textup).
\Part{(c)}
\textup{(Milnor \cite{milnor:quadraticmaps}, Silverman \cite{silverman:modulirationalmaps})}
$\Moduli_2^1$ is isomorphic to the affine plane $\AA^2$, and its natural GIT compactification is~$\PP^2$.
\Part{(d)}
\textup{(Levy \cite{MR2741188})} $\Moduli_d^1$ is a rational variety,
i.e., there exists a dominant rational map
$F:\PP^{2d-2}\dashrightarrow\Moduli_d^1$ such that~$F$ has a dominant
rational inverse.
\end{parts}
\end{theorem}

\begin{definition}
The \emph{automorphism group} of a dominant rational map $f:\PP^N\dashrightarrow\PP^N$
defined over~$\Kbar$ is the group
\[
  \Aut(f) := \bigl\{ \f\in\PGL_{N+1}(\Kbar) : f^\f = f \bigr\}.
\]
It is easy to check that $\Aut(f^\psi)=\psi^{-1}\Aut(f)\psi$, so the
$\PGL_{N+1}$-equivalence class of~$f$ determines the group~$\Aut(f)$ up to
$\PGL_{N+1}$-conjugacy.
\end{definition}

\begin{remark}
If~$f$ is an endomorphism, then it is known that~$\Aut(f)$ is finite,
and indeed its order may be bounded in terms of~$d$
and~$N$~\cite{MR2741188}, but for dominant rational maps, the
group~$\Aut(f)$ may be infinite.
\end{remark}

\begin{question}
  We collect some of the natural questions that one might ask
  about the moduli spaces~$\Moduli_d^N$ and about automorphism groups.
  \begin{parts}
  \Part{(a)}
    For $N\ge2$, are all \textup(some, none\textup) of the moduli
    spaces~$\Moduli_d^N$ rational varieties?
  \Part{(b)}
    For each $(N,d)$, describe the possible groups $\Aut(f)$ for
    $f\in\Moduli_d^N$, or at least give a good upper bound for
    \[
      U(N,d):=\sup_{f\in\Moduli_d^N}\#\Aut(f)
    \]
    in terms of~$N$ and~$d$.  It is known that~$U(N,d)$ is
    finite~\cite{MR2741188}, that
    $U(1,d)\le\max\{60,2d+2\}$~\cite[Example~2.54]{MR2884382}, and that
    $U(2,d)\le60d^6$~\cite[Theorem~6.2]{arxiv1509.06670}.
  \Part{(c)}
    Describe the subvarieties of $\Moduli_d^N$ corresponding to
    maps~$f$ whose automorphism group~$\Aut(f)$ is non-trivial.
    See~\cite{MR3709645} for a detailed analysis for $N=1$
    and~\cite{arxiv1607.05772} for some partial
    results for~$N=d=2$.  For example,
    \[
    \{f\in\Moduli_2^1\cong\AA^2:\Aut(f)\ne1\}=\text{cuspidal cubic curve in~$\AA^2$,}
    \]
    with~$\Aut(f)\cong\ZZ/2\ZZ$ except for the map~$f$ at the cusp,
    for which~$\Aut(f)\cong\Scal_3$.
  \Part{(d)}
    Most of the spaces~$\Moduli_d^N$ are singular.  Describe the
    singularity locus.  The case~$N=1$ is analyzed
    in~\cite{MR3709645}.
  \end{parts}
\end{question}

An arithmetic question closely related to dynamical moduli spaces is
that of fields of definition versus field of moduli, a problem that
also arises in the theory of abelian varieties.

\begin{definition}
Let $f\in\End_d^N(\Qbar)$, and let
\[
  \pi : \End_d^N(\Qbar)\longrightarrow\Moduli_d^N(\Qbar)
\]
be the natural projection map. The \emph{field of moduli of~$f$} is
the field $\QQ\bigl(\pi(f)\bigr)$, i.e., the field generated by the
coordinates of the point~$\pi(f)$. A \emph{field of definition for~$f$} is
any field~$L$ for which there is some~$\f\in\PGL_{N+1}(\Qbar)$ such
that~$f^\f\in\End_d^N(L)$.
\end{definition}

It is easy to see that every field of definition of~$f$ contains its
field of moduli, but there exist maps whose field of moduli is not a
field of definition. For example, one can show that the field of
moduli of the map $f(x)=i(x^3-1)/(x^3+1)$ is~$\QQ$, although~$\QQ$ is not a
field of definition for~$f$.

\begin{definition}
For~$f\in\End_d^N(\Qbar)$, let~$K_f$ denote the field of moduli
of~$f$, and define
\[
\operatorname{FOM/FOD}(N,d) := \sup_{f\in\End_d^N(\Qbar)}
\inf_{\substack{\text{$L$ is a field of}\\\text{definition of $f$}\\}}  [L:K_f].
\]
\end{definition}

In other words,~$\operatorname{FOM/FOD}(N,d)$ is the smallest
integer~$D$ such that every map in~$\End_d^N(\Qbar)$ has a field of
definition whose degree over its field of moduli is bounded by~$D$.
Of course, it is not clear a priori that~$\operatorname{FOM/FOD}(N,d)$
is finite.

\begin{theorem}
\begin{parts}
\Part{(a)}
For $N=1$ we have the exact values
\[
\operatorname{FOM/FOD}(1,d) = \begin{cases}
  1 &\text{if $d\ge2$ is even. \textup{(Silverman, 1995, \cite{silverman:fieldofdef})},} \\
  2 &\text{if $d\ge3$ is odd. \textup{(Hidalgo, 2014, \cite{MR3230378})}.} \\
\end{cases}
\]
\Part{(b)}
For all $N\ge1$ and all $d\ge2$, the quantity $\operatorname{FOM/FOD}(N,d)$ is finite
\textup{(Doyle--Silverman, 2018, \cite{fomfodPN2018})}.
\end{parts}
\end{theorem}

The fact that the bound for $\operatorname{FOM/FOD}(1,d)$ does not depend on~$d$ suggests
the following question.

\begin{question}
For a given $N\ge1$, is
\[
  \sup_{d\ge2} \operatorname{FOM/FOD}(N,d) < \infty?
\]
If so, what is its value as a function of~$N$?
\end{question}

Just as we did for dynamical modular curves in Section~\ref{section:dmc}, we
can define moduli spaces \emph{with level structure} by specifying
maps together with points or cycles of various shapes. For simplicity,
we restrict to the case of a single periodic point.

\begin{definition}
Let $N\ge1$, $n\ge1$, and $d\ge2$. We recall that~$\Per_n^*(f)$
denotes the set of points of exact period~$n$ for~$f$, and we define
\[
  \End_d^N[n] := \bigl\{ (f,P) : \text{$f\in\End_d^N $ and $P\in\Per_n^*(f)$} \bigr\}
  \subset \End_d^N\times\PP^N
\]
to be the space that classifies maps with a marked point of exact
period~$n$.  We define an action of~$\f\in\PGL_{N+1}$
on~$(f,P)\in\End_d^N[n]$ by the rule
\[
  (f,P)^\f := \bigl(f^\f,\f^{-1}(P)\bigr),
\]
and we let
\[
  \Moduli_d^N[n] := \End_d^n[n]/\PGL_{N+1}
\]
be the associated quotient space.  There is a natural
map
\[
  \text{$i_n:\Moduli_d^N[n]\to\Moduli_d^N$\quad induced by\quad $(f,P)\mapsto f$.}
\]
\end{definition}

\begin{example}
Let $N=1$ and $d=2$, so $\Moduli_2^1\cong\AA^2$. The set of polynomial
maps $\{f_c(x):=x^2+c\}$ gives a curve~$L$ in~$\Moduli_2^1$, and
Theorem~\ref{theorem:X1nirredgon} says that for all~$n\ge1$, the
inverse image curve~$i_n^{-1}(L)\cong X_1^\dyn(f_c;n)$ is irreducible,
and that its genus and gonality go to infinity as~$n\to\infty$.
Similarly, the set of rational maps $\{f_b(x):=bx+x^{-1}\}$ gives a
curve~$C$ in~$\Moduli_2^1$, and Theorem~\ref{theorem:manes2comp} says
that for all even~$n\ge2$, the inverse image curve~$i_n^{-1}(C)\cong
X_1^\dyn(f_c;n)$ has at least two irreducible components.  
\end{example}

Dynamical modular curves associated to natural 1-parameter families of
maps may be reducible. The following conjecture, the first part of
which was raised as a question by Poonen, says that there exists a
uniform bound for the amount of reducibility as we add more level
structure.

\begin{conjecture}
\label{conjecture:compandgen}
Let~$N\ge1$, let $d\ge2$, and let $\G\subset\Moduli_d^N$ be an
irreducible subvariety, e.g.,~$\G$ could be an irreducible curve. For
each~$n\ge1$, write~$i_n^{-1}(\G)$ as a union of irreducible
varieties,
\[
  i_n^{-1}(\G)=\G_{n,1}\cup\G_{n,2}\cup\cdots\cup\G_{n,\nu},
  \quad\text{where $\nu=\nu(N,d,n;\G)$.}
\]
\vspace{-10pt}
\begin{parts}
\Part{(a)}
\textup{(Poonen, 2012, \cite{arxiv1206.7104})}
%% In Poonen's article, he only says that Question 1.4 has a dynamical
%% analogue, so he doesn't exactly state this.
The number of irreducible components of $i_n^{-1}(\G)$ is bounded
independently of~$\G$ and~$n$, i.e., there is a constant~$C_1(N,d)$ such
that
\[
  \sup_{\substack{\textup{irreducible $\G\subset\Moduli_d^N$}\\n\ge1\\}} \nu(N,d,n;\G) \le C_1(N,d).
\]
\Part{(b)}
There is a constant $C_2(N,d)$ such that for all $n\ge C_2(N,d)$, all
irreducible $\G\subset\Moduli_d^N$, and all $1\le j\le\nu(N,d,n\;\G)$,
the variety~$\G_{n,j}$ is of general type.\footnote{The precise
  definition of \emph{general type} is somewhat technical, although
  for a smooth curve it simply means that the genus is at
  least~$2$. In general, suffice it to say that general type implies
  that the geometry of~$\G_{n,j}$ is very complicated.}
\Part{(c)}
The gonality \textup(see Definition~$\ref{definition:gonality}$) of the
components of~$i_n^{-1}(\G)$ satisfy
\[
  \lim_{n\to\infty} \min_{1\le i\le \nu(N,d,n;\G)}
  \operatorname{Gonality}(\G_{n,i}) = \infty .
\]
\Part{(d)}
Let~~$\G$ be an irreducible curve, so~$i_n^{-1}(\G)$ is a union of
irreducible curves. Then
\[
  \lim_{n\to\infty} \min_{1\le i\le \nu(N,d,n;\G)}
  \operatorname{Genus}(\G_{n,i}) = \infty .
\]
\end{parts}
\end{conjecture}

\begin{question}
In Conjecture~\textup{\ref{conjecture:compandgen}(c,d)}, give
\textup(accurate\textup) lower bounds for the genera and gonalities of
the components~$\G_{n,i}$ as~$n\to\infty$ when~$\G$ is an irreducible
curve.
\end{question}

%%%%%%%%%%%%%%%%%%%%%%%%%%%%%%%%%%%%%%%%%%%%%%%%%%%%%%%%%%%%%%%%%%%%%%
\section{Unlikely Intersections in Dynamics\WhoDone{Joe+Laura}}
\label{section:spdms}
\WhoWrite{Joe, based on Laura's notes}
%%%%%%%%%%%%%%%%%%%%%%%%%%%%%%%%%%%%%%%%%%%%%%%%%%%%%%%%%%%%%%%%%%%%%%
In Sections~\ref{section:ios} and~\ref{section:ppsv} we discussed
subvarieties containing a large number of preperiodic points
or a large number of points in an orbit. An orbit~$\Orbit_f(P)$
consists of a countable number of points, as does the set~$\PrePer(f)$
of preperiodic points, and the theme in those sections was that it is
unlikely for these countable sets to have large intersection with a
proper subvariety unless there is a geometric reason to expect the
intersection to be large.

The study of \emph{unlikely intersections} in arithmetic geometry has
seen a surge of interest, and there are many analogous problems that
can be formulated for dynamical systems. The guiding principle is that
families of objects should intersect infrequently unless there is
a geometric reason forcing the converse.

The following theorem, which is a dynamical analogue of a theorem of
Masser and Zannier for elliptic curves~\cite{MR2766181}, illustrates
the type of result that one seeks.  It it worth noting that the
finiteness in Theorem~\ref{theorem:BDMunlikely} takes place in moduli
space, i.e., in the~$t$-space that parameterizes the quadratic
polynomials~$x^2+t$.

\begin{theorem}
  [Baker--DeMarco, 2011, \cite{arxiv0911.0918}]
  \label{theorem:BDMunlikely}
For $t\in\CC$, let $f_t(x)=x^2+t$ be the usual quadratic polynomial.  
Let $a,b\in\CC$ complex numbers with $a^2\ne b^2$. Then
\[
  \bigl\{ t\in\CC : \text{$a\in\PrePer(f_t)$ and $b\in\PrePer(f_t)$}
  \bigr\}~\text{is finite.}
\]
\end{theorem}

More generally, we might allow the points and the maps to vary with
the parameter. So  let $T$ be a non-singular algebraic curve, and 
let
\begin{equation}
  \label{eqn:aTPbTPfTEnd}
  a:T\longrightarrow \PP^1,\quad
  b:T\longrightarrow \PP^1,\quad
  f:T\longrightarrow \End_d^1,
\end{equation}
be non-constant algebraic maps.  We view~$a_t$ and~$b_t$ as defining
families of points parameterized by~$t\in T$, and similarly we
view~$f$ as a family of rational maps $f_t:\PP^1\to\PP^1$
parameterized by~$t\in T$.

The independence of~$a$ and~$b$ in Theorem~\ref{theorem:BDMunlikely},
which we stated as~$a^2\ne b^2$, is equivalent to saying
that~$f_t(a)\ne f_t(b)$.  More generally, for varying families of maps
and points, we want to define what it means for~$a$ and~$b$ to be
\emph{$f$-dynamically related}. For example, this will include the case
that $f^n\circ a=f^m\circ b$ for some~$n\ge0$ and~$m\ge0$, but this is
not the only situation. For example, if there is some
$\f:T\to\End_d^1$ that commutes with~$f$, i.e.,
$f_t\circ\f_t=\f_t\circ f_t$ for all $t\in T$, then~$a$
and~$\f\circ{a}$ are dynamically related for every~$a$. Similarly,
if~$a_t\in\PrePer(f_t)$ for every~$t\in T$, then~$a$ is
dynamically related to every~$b$. The varied nature of these examples
gives some justification for the technical nature of the following
formal definition.

\begin{definition}
Let $T$ be a non-singular algebraic curve, and let~$a$,~$b$, and~$f$
be maps as in~\eqref{eqn:aTPbTPfTEnd}.  Let $\overline{k(T)}$ be an
algebraic closure of the function field of~$T$, so we may view~$a$
and~$b$ as elements of $\PP^1\bigl(\overline{k(T)}\bigr)$, and we may
view~$f$ as a map
$f:\PP^1\bigl(\overline{k(T)}\bigr)\to\PP^1\bigl(\overline{k(T)}\bigr)$. Then
we say that~\emph{$a$ and~$b$ are $f$-dynamically related} if there is
an algebraic curve
\[
  C \subset (\PP^1\times\PP^1)\bigl(\overline{k(T)}\bigr)
\]
satisfying
\[
  (a,b)\in C \quad\text{and}\quad (f\times f)(C) \subseteq C.
\]
\end{definition}

\begin{conjecture}
[Ghioca--Hsia--Tucker, 2013; \cite{MR3095224}; DeMarco, 2016, \cite{MR3531361}]
\label{conjecture:unlikelyorbits}
Let $T$ be a non-singular algebraic curve, and let~$a$,~$b$, and~$f$
be maps as in~\eqref{eqn:aTPbTPfTEnd} for some~$d\ge2$ and with~$k$ an
algebraically closed field.
Then the following are equivalent\textup:
\begin{parts}
  \Part{(a)}
  The set
  $\bigl\{ t\in T(k) : \text{$a_t$ and $b_t$ are both in $\PrePer(f_t)$} \bigr\}$ is infinite.
  \Part{(b)}
  The maps $a$ and $b$ are $f$-dynamically related.
\end{parts}
\end{conjecture}

Using the theory of canonical heights, which are defined and discussed
in Section~\ref{section:ch}, we can formulate a strengthened version
of Conjecture~\ref{conjecture:unlikelyorbits}.  In order to relate the
conjectures, we note that the canonical height is a
function~$\hhat_f:\PP^1(\Qbar)\to[0,\infty)$ with the property that
$\hhat_f(\a)=0$ if and only $\a\in\PrePer(f)$. Bogomolov's original
small height conjectures, which were for points on subvarieties of
tori and abelian varieties, were proven by Ullmo~\cite{MR1609514}
and Zhang~\cite{MR1609518}.

\begin{conjecture}
[Dynamical Bogomolov Conjecture]    
\label{conjecture:bogounlikelyorbits}
Let $T$,~$a$,~$b$, and~$f$ be maps as in~\eqref{eqn:aTPbTPfTEnd} for
some~$d\ge2$, with all maps defined over~$\Qbar$.
Then the following are equivalent\textup:
\begin{parts}
  \Part{(a)}
  For every $\e>0$, the set
  $\bigl\{ t\in T(\Qbar) : \text{$\hhat_{f_t}(a_t)\le\e$ and $\hhat_{f_t}(b_t)\le\e$} \bigr\}$ is infinite.
  \Part{(b)}
  The maps $a$ and $b$ are $f$-dynamically related.
\end{parts}
\end{conjecture}

We recall that~$P\in\PP^1$ is a \emph{critical point} for a rational
map $f:\PP^1\to\PP^1$ if~$f$ fails to be locally one-to-one in a
neighborhood of~$P$, and that~$f$ is said to be \emph{post-critically
  finite} (PCF) if all of its critical points are preperiodic.  PCF
maps appeared briefly in Section~\ref{section:ar} and will reappear in
Section~\ref{section:ch}. As suggested by Table~\ref{table:dictionary}
in the Introduction, we view the collection of~PCF maps in the moduli
space~$\Moduli_d^1$ as being in some ways analogous to the collection
of~CM abelian varieties in the moduli space~$\Acal_g$ of abelian
varieties.  In the setting of abelian varieties, there are famous
conjectures of Andr\'e and Oort which assert that a
subvariety~$X\subset\Acal_g$ contains a Zariski dense set of~CM points
if and only if there is an underlying geometric reason. In evocative
terminology, the classical Andr{\'e}--Oort conjecture is often phrased
as saying that special points of~$\Acal_g$ are Zariski dense in~$X$ if
and only if~$X$ itself is special.  This formulation would have an
obvious dynamical analogue if only we knew exactly what it means
for~$X\subset\Moduli_d^1$ to be a dynamically special subvariety!

In order to formulate a more precise conjecture, we return to the
construction of the spaces~$\End_d^1$ and~$\Moduli_d^1$ in
Section~\ref{section:dms}. A map~$f\in\End_d^1$ has exactly~$2d-2$
critical points, counted with multiplicities. We define
\[
\End_d^\crit:= \left\{ (f,P_1,\ldots,P_{2d-2}) :
\begin{tabular}{@{}l@{}}
  $f\in\End_d^1$, and $P_1,\ldots,P_{2d-2}$\\
  are critical points of $f$\\
\end{tabular} \right\}
\]
to be the space of maps with marked critical points. The
group~$\PGL_2$ acts as usual, cf.\ the construction
of~$\Moduli_d^N[n]$ in Section~\ref{section:dms}, and one can
construct the GIT quotient
space~$\Moduli_d^\crit:=\End_d^\crit/\PGL_2$.\footnote{We are being
  slightly imprecise here. It is easy to construct the GIT quotient
  mdouli space of maps~$f$ having~$2d-2$ \emph{distinct} marked
  critical points, but one must be more careful when dealing with
  critical points that have higher
  multiplicity. See~\cite{moduliportrait2017}, for example, for a
  discussion.}

\begin{conjecture}[Dynamical Andr{\'e}--Oort Conjecture]
\label{conjecture:andreoort}
\textup{(Baker--DeMarco, ${}\approx2011$, \cite{arxiv1211.0255,MR2884382})}
Let $X\subseteq\Moduli_d^\crit$ be an algebraic subvariety such that
the PCF maps in~$X$ are Zariski dense in~$X$. Then~$X$ is cut out by
``\emph{critical orbit relations}.''
\end{conjecture}

Formulas of the form $f^n(P_i)=f^m(P_j)$ define critical
point relations, but other relations may arise from symmetries of~$f$,
and even subtler relations can come from ``hidden relations'' due to
sub-dynamical systems. See~\cite[Remark~6.58]{MR2884382} for an
example due to Ingram.  Although there is not, as yet, a complete
description of what is meant by the phrase ``critical orbit
relations'' in Conjecture~\ref{conjecture:andreoort}, we note that
Conjecture~\ref{conjecture:unlikelyorbits} implies the following
precise conjecture for 1-dimensional families.

\begin{conjecture}[Dynamical Andr{\'e}--Oort Conjecture for 1-Dimensional Families]
\label{conjecture:andreoort1dim}
\textup{(Baker--DeMarco, 2013, \cite{arxiv1211.0255})}
Let $T$ be a non-singular algebraic curve, and let
\[
  T \longrightarrow\Moduli_d^\crit,\quad
  t \longmapsto (f_t,P_{1,t},\ldots,P_{2d-2,t})
\]
be a family of rational maps with marked critical points.  Then the
following are equivalent\textup:
\begin{parts}
  \Part{(a)}
  There are infinitely many~$t\in T$ such that~$f_t$ is PCF.
  \Part{(b)}
  For every $i$ and $j$, the maps $P_i:T\to\PP^1$ and $P_j:T\to\PP^1$
  are $f$-dynamically related.
\end{parts}
\end{conjecture}

Building on a number of earlier partial results
\cite{arxiv1211.0255,arxiv0911.0918,MR3801434,MR3095224,MR3801489},
Favre and Gauthier have recently announced a proof of
Conjecture~\ref{conjecture:andreoort1dim} for all families of
polynomials, but the conjecture remains open for families of rational
maps.

%%%%%%%%%%%%%%%%%%%%%%%%%%%%%%%%%%%%%%%%%%%%%%%%%%%%%%%%%%%%%%%%%%%%%%
\section{Good Reduction of Maps and Orbits\WhoDone{Joe}}
\label{section:goodred}
\WhoWrite{Joe}
%%%%%%%%%%%%%%%%%%%%%%%%%%%%%%%%%%%%%%%%%%%%%%%%%%%%%%%%%%%%%%%%%%%%%%

Let~$p\in\ZZ$ be a prime. In this section we use tildes to denote the
reduction modulo~$p$ map,
\[
  \ZZ\longrightarrow\FF_p,\quad c\longmapsto\tilde c,
\]
and if we need to specify the prime~$p$, we write~$\tilde c\bmod p$.
Similarly, we define a reduction modulo~$p$ map
\[
  \ZZ[X_0,\ldots,X_N]\longrightarrow\FF_p[X_0,\ldots,X_N],\quad
  F(X_0,\ldots,X_N)\longmapsto\tilde F(X_0,\ldots,X_N),
\]
obtained by reducing the coefficients of a polynomial.

\begin{definition}
Let $P\in\PP^N(\QQ)$. The \emph{reduction of~$P$ modulo~$p$},
denoted~$\tilde P$, is defined as follows. First choose homogeneous
coordinates~$P=[c_0,\ldots,c_N]$ for~$P$ satisfying
\[
  c_0,\ldots,c_N\in\ZZ\quad\text{and}\quad \gcd(c_0,\ldots,c_N)=1,
\]
and then set
\[
  \tilde P = [\tilde c_0,\ldots,\tilde c_N] \in \PP^N(\FF_p).
\]
This gives a well-defined \emph{reduction modulo~$p$ map}
$\PP^N(\QQ)\to\PP^N(\FF_p)$.
\end{definition}

\begin{definition}
Let~$f\in\End_d^N(\QQ)$, i.e.,~$f:\PP^N\to\PP^N$ is a morphism defined
over~$\QQ$. The \emph{reduction of~$f$ modulo~$p$},
denoted~$\tilde{f}$, is defined as follows. First write~$f$ in the
form
\[
  f=[f_0,\ldots,f_N]\quad\text{with}\quad f_0,\ldots,f_N\in\ZZ[X_0,\ldots,X_N],
\]
and if necessary, divide~$f_0,\ldots,f_N$ by an appropriate
integer to ensure that the gcd of the collection of all of the coefficients
of~$f_0,\ldots,f_N$ is equal to~$1$. Such an~$f$ is said to be written
in \emph{normalized form}. We then set
\[
  \tilde f=[\tilde f_0,\ldots,\tilde f_N]\quad\text{with}\quad \tilde
    f_0,\ldots,\tilde f_N\in\FF_p[X_0,\ldots,X_N],
\]
where~$\tilde f_i$ is obtained by reducing the coefficients of~$f_i$
modulo~$p$.
\end{definition}

\begin{remark}
The fact that a map $f=[f_0,\ldots,f_N]\in\End_d^N(\QQ)$ has
degree~$d$ means that~$f_0,\ldots,f_N$ are homogeneous polynomials of
degree~$d$ having no common factors in~$\QQ[X_0,\ldots,X_N]$. However,
when we reduce modulo~$p$, it may happen that the polynomials pick up
a common factor.  For example, $F=X^2+3Y^2$ and $G=X^2+XY$ have no
common factors in $\QQ[X,Y]$, but reducing modulo~$3$, we find that
$\tilde{F}\bmod3=X^2$ and $\tilde{G}\bmod3=X^2+XY$ acquire a common
factor of~$X$ in $\FF_3[X,Y]$.
Thus~$\tilde f(X,Y)=[X^2,X^2+XY]=[X,X+Y]$ has degree~$1$, which shows
that reduction modulo~$p$ does not necssarily respect the degree of a
map.
\end{remark}

\begin{definition}
Let $f:\PP^N\to\PP^N$ be a morphism defined over~$\QQ$, and let~$p$ be
a prime. The map~$f$ has (\emph{naive}\footnote{We discuss later why
  this definition is ``naive.''}) \emph{good reduction modulo~$p$} if
\[
  \deg(\tilde f\bmod p) = \deg(f).
\]
\end{definition}

The following elementary result illustrates the utility of good
reduction. In words, it says that if~$f$ has good reduction, then
reduction commutes with both iteration and evaluation.

\begin{proposition}
\label{proposition:goodredeval}
Let $f:\PP^N\to\PP^N$ be a rational map defined over~$\QQ$,
let~$n\ge1$, let~$P\in\PP^N(\QQ)$, and let~$p$ be a prime of good
reduction for~$f$. Then~$f^n$ has good reduction, and we have
\[
  \widetilde{f^n}=\tilde f^n
  \quad\text{and}\quad
  \widetilde{f^n(P)}=\tilde f^n(\tilde P).
\]
\end{proposition}

Let~$f$ have good reduction, and let~$P\in\PP^N(\QQ)$ be an
$f$-periodic point of period~$n$.  It follows easily from
Proposition~\ref{proposition:goodredeval} that~$\tilde P$ is $\tilde
f$-periodic and that its $\tilde f$-period divides~$n$. It is possible
to say much more.  The following theorem is an amalgamation of
results. We have stated it for maps of~$\PP^N$ defined over~$\QQ$, but
it holds much more generally for self-maps of algebraic varieties over
number fields; see especially~\cite{MR2741181}.

\begin{theorem}
\label{theorem:permodpnmrpe}
\textup{(Morton-Silverman, 1994/5,
  \cite{mortonsilverman:rationalperiodicpoints,mortonsilverman:dynamicalunits};
  Pezda, 1994, \cite{pezda:polynomialcycles2}; Zieve, 1996,
  \cite{zieve:thesis} Hutz, 2009, \cite{MR2741181})} Let
$f:\PP^N\to\PP^N$ be a morphism defined over~$\QQ$, let~$p$ be a prime
of good reduction for~$f$, let~$P\in\PP^N(\QQ)$ be an $f$-periodic
point of exact period~$n$, and let~$m$ be the exact period of the
reduced point~$\tilde P\in\PP^N(\FF_p)$ for the reduced
map~$\tilde{f}$. Then
\[
  \text{$n = m$\quad or\quad $n=mr$\quad or\quad $n=mrp^e$,}
\]
where~$r$ divides~$p^N-1$, and where there is an explicit upper bound for~$e$.
\textup(For example, if $N=1$ and $p\ge3$, then $e\le1$.\textup)
%% and if $N=1$ and $p=2$, then $e\le 3$.
\end{theorem}

Using the fact that~$m$ in Theorem~\ref{theorem:permodpnmrpe} is at
most~$\#\PP^N(\FF_p)$, we can use Theorem~\ref{theorem:permodpnmrpe}
to prove a ``semi-uniform'' version of the Dynamical Uniform
Boundedness Conjecture (Conjecture~\ref{conjecture:MSUBC}) for
periodic points. More precisely, the theorem gives a bound for
$\#\Per\bigl(f,\PP^N(K)\bigr)$ that depends only
on~$N$,~$\deg(f)$,~$[K:\QQ]$, and the number of primes at which~$f$
has bad reduction.  So in lieu of full uniform boundness, one might
ask the following question.

\begin{question}
Find a ``nice'' function~$C(N,d,D,s)$ so that every degree~$d$
morphism $f:\PP^N\to\PP^N$ defined over a field~$K$ of degree at
most~$D$ and having at most~$s$ primes of bad reduction satisfies
\[
  \#\PrePer\bigl(f,\PP^N(K)\bigr) \le C(N,d,D,s).
\]
\end{question}

If one restricts to polynomial maps on~$\PP^1$, then
Benedetto~\cite{MR2339471} proves that one can take
$C_{\textup{poly}}(1,d,D,s)$ equal to a mutiple of
$d^3(D+s)\log(D+s)$, while Canci and Paladino~\cite{MR3556260} give a
weaker bound for rational maps on~$\PP^1$.

We called our earlier definition of good reduction the ``naive
definition.'' The reason it is naive is due to the fact that
conjugating a map~$f:\PP^N\to\PP^N$ by a linear
transformation~$\f\in\PGL_{N+1}$ gives a map~$f^\f$ having the same
dynamical behavior as~$f$. But when we  are reducing
modulo~$p$, conjugation by an element of~$\PGL_{N+1}(\QQ)$ may turn
good reduction into bad reduction, or vice versa.

\begin{example}
The map
\[
  f(X,Y) = [pX^2+XY,Y^2]
\]
has bad reduction at~$p$, since $\tilde f=[XY,Y^2]=[X,Y]$.
But if we conjugate~$f$ by the linear map $\f(X,Y)=[X,pY]$, we find that
\begin{multline*}
  f^\f(X,Y)=\f^{-1}\circ f \circ \f(X,Y)
  = \f^{-1}\circ f(X,pY)
  = \f^{-1}(pX^2+pXY,p^2 Y^2) \\
  = [pX^2+pXY,pY^2]
  =[X^2+XY,Y^2].
\end{multline*}
Thus~$f^\f$ has (naive) good reduction, so it would have been
advisable to use~$f^\f$ instead of~$f$. This leads to a better
definition.
\end{example}

\begin{definition}
Let $f:\PP^N\to\PP^N$ be a rational map defined over~$\QQ$, and
let~$p$ be a prime. We say that~$f$ has \emph{good reduction
  modulo~$p$} if there is a change of variables $\f\in\PGL_{N+1}(\QQ)$
such that the conjugate map~$f^\f$ has good reduction modulo~$p$.
\end{definition} 

\begin{remark}
There  are algorithms due to Benedetto~\cite{MR3263950},
Bruin--Molnar~\cite{arxiv1204.4967}, Rumely~\cite{MR3361223}, and
Szpiro--Tepper--Williams~\cite{MR3119234}, using a variety of methods,
that find a change of variables for a given map~$f$ on~$\PP^1$ making
the bad reduction as good as possible, and in particular determining
if~$f$ has good reduction. They play the role in arithmetic dynamics
of a similar algorithm of Tate~\cite{MR0393039} that finds minimally
bad equations for elliptic curves.
\end{remark}

A famous theorem of Shafarevich~\cite[Theorem~IX.6.1]{MR2514094} says
that if~$S$ is a finite set of primes in the ring of integers of a
number field~$K$, then there are only finitely many elliptic curves
defined over~$K$ having good reduction at all primes not
in~$S$. Shafarevich further conjectured, and Faltings~\cite{MR718935}
proved, the much deeper result that the same is true for abelian
varieties of arbitrary fixed dimension.  For those who are familiar
with these theorems, it may come as a surprise that there are lots of
maps $f:\PP^N\to\PP^N$ that have good reduction at every prime. For
example, if $f:\PP^1\to\PP^1$ has the dehomogenized form
\begin{equation}
  \label{eqn:fmonicpolyP1}
  f(x)=x^d+a_1x^{d-1}+\cdots+a_d\quad\text{with $a_1,\ldots,a_d\in\ZZ$,}
\end{equation}
then~$f$ has good reduction at every prime, since it is clear that
$\deg(\tilde f\bmod p)=d$ for all~$p$. Thus a naive generalization of
Shafarevich's conjecture to dynamical systems is false.

However, polynomial maps~\eqref{eqn:fmonicpolyP1} form only a small
part of the totality of rational maps on~$\PP^1$. Thus one might ask,
within the collection of all dynamical systems, how widespread are the
maps having good reduction outside some fixed finite set of places. As
in Section~\ref{section:dms}, we parameterize degree~$d$ dynamical
systems on~$\PP^N$ by the points of the variety~$\Moduli_d^N$.  Then
for example, we have~$\dim\Moduli_d^1=2d-2$, while the set of good
reduction polynomial maps~\eqref{eqn:fmonicpolyP1} in $\Moduli_d^1$
lies on a subvariety of dimension~$d-1$.

\begin{definition}
Let $N\ge1$ and $d\ge2$. The \emph{Shafarevich dimension} for
degree~$d$ dynamical systems on~$\PP^N$, denoted  $\shafdim_d^N$, is the quantity
\[
  \sup_{\substack{\text{$K$ a number field}\\  \text{$S$ a finite set of primes}   \\}} \hspace{-.75em}
  \dim\overline{ \bigl\{ f\in\Moduli_d^N(K) : \text{$f$ has good reduction at all primes not in $S$} \bigr\} },
\]
where the overline denotes the Zariski closure.
\end{definition}

Good reduction polynomials~\eqref{eqn:fmonicpolyP1} show that
$\shafdim_d^1\ge d-1$. The next result gives a small improvement,
while suggesting a fundamental question.

\begin{proposition}
\textup{(Silverman, 2017, \cite[Proposition~12]{arxiv1703.00823})} The
Shafarevich dimension for degree~$d$ dynamical systems
on~$\PP^1$ satisfies
\[
  \shafdim_2^1=2,\quad \shafdim_3^1=4,\quad\text{and}\quad
  \shafdim_d^1\ge d+1\quad\text{for all $d\ge4$.}
\]
\end{proposition}

\begin{question}
\label{question:shafdimdN}
What is the exact value of $\shafdim_d^N$\textup{?} In particular, is it true that
\[
  \shafdim_d^1=d+1\quad\text{for all $d\ge3$\textup{?}}
\]
\end{question}

\begin{remark}
Petsche~\cite{MR2999309} shows that for certain codimension~$3$
families in $\Moduli_d^1$, there are only finitely many maps having
good reduction outside of~$S$.  We also mention that Petsche and
Stout~\cite{MR3293730} more-or-less define the Shafarevich dimension
$\shafdim_d^1$ on~$\PP^1$ and ask if $\shafdim_d^1=2d-2$, i.e., they
ask if good reduction maps are Zariski dense in~$\Moduli_d^1$. They
prove that this is true for $d=2$.
\end{remark}

\begin{remark}
Two maps $f_1,f_2\in\End_d^N(K)$ are said to be \emph{$K$-twists} of
one another if there is a $\Kbar$-change of
variables~$\f\in\PGL_{N+1}(\Kbar)$ such that~$f_2=f_1^\f$. We
consider~$f_1$ and~$f_2$ to be $K$-equivalent if it is
possible to find such a~$\f$ in~$\PGL_{N+1}(K)$. For example, the
maps~$f_1(x)=x^3$ and~$f_2(x)=2x^3$ are inequivalent $\QQ$-twists of
one another, since $f_2(x)=f_1^\f(x)$ for $\f(x)=\sqrt2 x$.
Stout~\cite{MR3273498} has proven that for a given~$f\in\End_d^N(K)$,
there are only finitely many equivalence classes of maps in
$\End_d^N(K)$ that are $K$-twists of~$f$ and have good reduction at
all primes not in~$S$. 
\end{remark}

For any finite set of primes~$S$, we have seen that there are
infinitely many self-maps~$f$ of~$\PP^N$ having good reduction at all
primes not in~$S$. As noted earlier, this means that the natural
dynamical analogue of the Shafarevich--Faltings finiteness theorem for
abelian varieties is false.  Various authors have attempted to rescue
the finiteness by adding a condition that various parts of certain
$f$-orbits have good reduction.  In order to state a sharp result for
maps of~$\PP^1$, we need a definition, which we state for general
number fields.

\begin{definition}
Let $K$ be a number field, let~$S$ be a finite set of primes of the
ring of integers of~$K$, and let~$R_S$ be the ring of $S$-integers
of~$K$. For~$n\ge1$ and $d\ge1$, define~$\GR_d^1[n](K,S)$ to be the
set of triples~$(f,Y,X)$, where $f\in\End_d^1(K)$ and $Y\subseteq
X\subset\PP^1(K)$ are finite sets, satisfying the following conditions:
\begin{parts}
  \Part{\textbullet}
  $X=Y\cup f(Y)$.
  \Part{\textbullet}
  $X$ is $\Gal(\Kbar/K)$-invariant.
  \Part{\textbullet}
  $\sum_{P\in Y} e_f(P)=n$, where~$e_f(P)$ is the ramification index of~$f$ at~$P$.
  \Part{\textbullet}
   $f$  has good reduction at all primes not in $S$.
  \Part{\textbullet}
  The points in $X$ remain distinct when reduced modulo~$\gp$ for every
  prime $\gp\notin S$.
\end{parts}
\end{definition}

There is a natural action of~$\PGL_2(R_S)$ on~$\GR_d^1[n](K,S)$,
where~$\PGL_2(R_S)$ is (essentially) the group of $2$-by-$2$
matrices~$\f$ with~$\det(\f)\in R_S^*$. This action is given by the
formula
\[
  \f\cdot(f,Y,X) = \bigl(f^\f,\f^{-1}(Y),\f^{-1}(X)\bigr).
\]
The following gives a strong Dynamical Shafarevich Theorem for~$\PP^1$
with orbit points; see also earlier work of
Szpiro--Tucker~\cite[(2008)]{arxiv0603436},
Szpiro--West~\cite[(2017)]{arxiv1705.05489}, and
Petsche--Stout~\cite[(2015)]{MR3293730}.

\begin{theorem}
\label{theorem:shafconjP1}
\textup{(Silverman, 2017, \cite{arxiv1703.00823})}
Let~$K/\QQ$ be a number field,  let~$S$ be a finite set of primes
of~$K$, and let~$d\ge2$. Then for all $n\ge2d+1$, the set 
\[
  \GR_d^1[n](K,S)/\PGL_{2}(R_S)~\text{is finite.}
\]
\end{theorem}

\begin{question}
Is there a natural generalization of Theorem~$\ref{theorem:shafconjP1}$
to maps of~$\PP^N\/$\textup{?} \textup(We mention that the naive generalization is false;
see~\cite[Section~8]{arxiv1703.00823}.\textup)
\end{question}

It is not hard to construct examples showing that
Theorem~\ref{theorem:shafconjP1} is false for~$n=2d$;
see~\cite[Proposition~11]{arxiv1703.00823}. However, if one further
specifies the exact orbit structure of the map $f:Y\to X$, then some
configurations with~$n=2d$ permit infinitely many good reduction
triples, while other configurations allow only finitely many. Without
giving precise definitions, which may be found
in~\cite{moduliportrait2017,arxiv1703.00823}, we informally say that
an orbit structure is specified by a \emph{portrait}~$\Pcal$,
for example
\begin{equation}
  \label{eqn:portrait112}
  \xymatrix{
    \Pcal: & {\bullet} \ar@(dr,ur)[]_{} \\
  }
  \qquad
  \xymatrix{
    {\bullet} \ar[r]_{} & {\bullet}  \\
  }
  \qquad
  \xymatrix{
    {\bullet} \ar@(dl,dr)[r]_{}   & {\bullet}   \ar@(ur,ul)[l]_{}  \\
  }
\end{equation}
and we define the set~$\GR_d^1[\Pcal](K,S)$ to be the subset of an
appropriate~$\GR_d^1[n](K,S)$ for which $f:Y\to X$ models the
portrait~$\Pcal$.  Thus for the portrait~\eqref{eqn:portrait112}, an
element~$(f,Y,X)$ of~$\GR_d^1[\Pcal](K,S)$ consists of a good
reduction map~$f:\PP^1\to\PP^1$, a set~$X$ containing~$5$ points that
remain distinct modulo all primes not in~$S$, and a subset~$Y\subset
X$ containing of~$4$ points.  Further, the points in~$Y$ consist of a
fixed point of~$f$, a $2$-cycle for~$f$, and a fourth point that~$f$
maps to the point in~$X$ that is not in~$Y$.

\begin{question}
\label{question:critGRd1P}
Give criteria, in terms of an integer~$d\ge2$ and the geometry of a
portrait~$\Pcal$, that distinguish between the following two
situations\textup:
\begin{parts}
\Part{(1)}
$\GR_d^1[\Pcal](K,S)/\PGL_2(R_S)$ is finite for all number fields~$K$ and
all finite sets of primes~$S$.
\Part{(2)}
$\GR_d^1[\Pcal](K,S)/\PGL_2(R_S)$ is infinite for some number fields~$K$ and
some finite sets of primes~$S$.
\end{parts}
\end{question}

\begin{remark}
For a given portrait~$\Pcal$, one can define a Shafarevich dimension
$\shafdim_d^1[\Pcal]$ associated to~$\GR_d^1[\Pcal]$. Then
Questions~\ref{question:shafdimdN} and~\ref{question:critGRd1P} may be
generalized by asking for a formula, or algorithm, to
compute~$\shafdim_d^1[\Pcal]$ in terms of~$d$ and the geometry
of~$\Pcal$. See~\cite{arxiv1703.00823} for details, as well as for a
complete list of the values of~$\shafdim_2^1[\Pcal]$ for the~$35$
portraits that contain at most~$4$ points and are allowable for
degree~$2$ maps.
\end{remark}

We close this section with another sort of dynamical good reduction
problem. In Section~\ref{section:dmc} we discussed the dynatomic
polynomial equation,
\[
  Y_1^\dyn(n) : \F_n(x,z)=0,
\]
whose solutions~$(\a,c)$ characterized points~$\a$ of formal
period~$n$ for the quadratic polynomial~$f_c(x)=x^2+c$. The
curve~$Y_1^\dyn(n)$ is an irreducible non-singular plane curve, and as
such, we can ask for which primes~$p$ does the curve~$Y_1^\dyn(n)$
have bad reduction. Dynamically, this should correspond to some sort
of intrinsic collapse among the periodic $n$-cycles of~$f_c$. By way
of analogy, we note that the classical elliptic modular curve~$X_1(n)$
has bad reduction at exactly the primes dividing~$n$, which is quite
reasonable, since the $n$-torsion points on an elliptic curve collapse
when reduced modulo~$p$ if~$p$ divides~$n$. We might similarly expect
that~$Y_1^\dyn(n)$ has bad reduction at the primes dividing~$n$, but
it also often has bad reduction at other sporadic-looking primes.

Morton~\cite{morton:dynamicalmoduli} defines a discriminant~$D_n$
whose prime divisors include all primes of bad reduction
for~$Y_1^\dyn(n)$, but recent work~\cite{arxiv1703.04172} shows that
in many instances, this list is too large. For example, the
curve~$Y_1^\dyn(5)$ has bad reduction only at~$5$ and~$3701$, but the
set of prime divisors of~$D_5$ is~$\{3,5,11,31,3701,4712,86131\}$.
Similarly, although~$D_7$ has many prime divisors, the dynatomic
modular curve $Y_1^\dyn(7)$ has bad reduction at exactly two primes,
\[
  7 \quad\text{and}\quad 84562621221359775358188841672549561.
\]

\begin{theorem}
\label{theorem:DKOPRSW}
\textup{(Doyle--Krieger--Obus--Pries--Rubinstein-Salzedo--West~\cite{arxiv1703.04172})}
Let $D_n$ be Morton's discriminant associated to the $n$th iterate
of~$x^2+c$; see~\cite{arxiv1703.04172,morton:dynamicalmoduli}
for the precise definition of~$D_n$.
\begin{parts}
\Part{(a)}
Let $n\ge4$, and let~$p$ be an odd prime divisor
of~$n$. Then~$Y_1^\dyn(n)$ has bad reduction at~$p$.
\Part{(b)}
Let $p$ be an odd prime divisor of~$D_n$ such that~$p^2$ does not
divide~$D_n$. Then $Y_1^\dyn(n)$ has bad reduction at~$p$.
\end{parts}
\end{theorem}

Theorem~\ref{theorem:DKOPRSW} is a small part of the material
in~\cite{arxiv1703.04172}, which considers more generally the
polynomials~$x^m+c$, analyzes good and bad reduction of
both~$Y_1^\dyn(n)$ and~$Y_0^\dyn(n)$, and proves in many cases
that~$\widetilde{Y_0^\dyn(n)}\bmod p$ is irreducible.  However, we
still lack an intrinsic explanation for large sporadic primes of bad
reduction.

\begin{question}
\label{question:Y0Y1badred}
Let $Y_0^\dyn(n)$ and $Y_1^\dyn(n)$ be the dynamical modular curves
classifying points, respectively cycles, of formal period~$n$ attached
to the family of maps~$x^2+c$, or~$x^m+c$, or some other natural
$1$-parameter family of rational maps on~$\PP^1$ such as $bx+x^{-1}$.
\begin{parts}
\Part{(a)}
Describe the primes of bad reduction for $Y_0^\dyn(n)$ and $Y_1^\dyn(n)$.
\Part{(b)}
Give a natural geometric/dynamical explanation for the ``sporadic''
primes of bad reduction, which we imprecisely define as primes not
dividing~$n$ and not dividing a finite set of bad primes associated to
the family~$f$.
\Part{(c)}
Describe the primes~$p$ for which the reductions of $Y_0^\dyn(n)$
and/or $Y_1^\dyn(n)$ modulo~$p$ are geometrically irreducible, e.g.,
for the latter, this is asking for the primes~$p$ such that the
polynomial~$\F_n(x,z)$ does not factor in the ring~$\bar\FF_p[x,z]$.
\end{parts}  
\end{question}

Note that part of Question~\ref{question:Y0Y1badred}(b) is to
precisely define what it means for a bad reduction prime to be
sporadic.  For example, the results in~\cite{arxiv1703.04172} suggest
that for~$x^m+c$, the sporadic primes for $n$-periodic points are
those that do not divide~$2nm$.

%%%%%%%%%%%%%%%%%%%%%%%%%%%%%%%%%%%%%%%%%%%%%%%%%%%%%%%%%%%%%%%%%%%%%%
\section{Orbits of Rational Maps\WhoDone{Joe}}
\label{section:orm}
\WhoWrite{Joe}
%%%%%%%%%%%%%%%%%%%%%%%%%%%%%%%%%%%%%%%%%%%%%%%%%%%%%%%%%%%%%%%%%%%%%%

Let $f:\PP^N\dashrightarrow\PP^N$ be a dominant rational map. We
recall that $\PP^N(K)_f$ denotes the set of points~$P\in \PP^N(K)$ for
which the entire forward orbit~$\Orbit_f(P)$ is well-defined.  Working
over~$\CC$, it is clear that the wandering points in~$\PP^N(\CC)_f$
are Zariski dense, since the the complement
$\bigcup_{n\ge0}f^{-n}\bigl(I(f)\bigr)$ is a countable union of proper
subvarieties of~$\PP^N$.  If we replace~$\CC$ with~$\Qbar$, this
statement becomes much harder to prove, since~$\PP^N(\Qbar)$ itself is
countable.  Indeed, it is not even clear, a priori,
that~$\PP^N(\Qbar)_f$ is non-empty.

\begin{theorem}
\label{theorem:Qbarwander}
\textup{(Amerik, 2011, \cite{MR2784670})}
Let $f:\PP^N\dashrightarrow\PP^N$ be a dominant rational map defined over~$\Qbar$.
Then the wandering points in~$\PP^N(\Qbar)_f$  are Zariski dense in~$\PP^N$.
\end{theorem}

Amerik's proof uses deep results of Hruskovski involving model theory
and logic. It would be interesting to find a more elementary,
dynamically inspired, proof.

One might ask for more, namely that there are individual points~$P$
whose orbit~$\Orbit_f(P)$ is Zariski dense. This cannot be true in
general, since there is a natural obstruction if~$f$ 
fibers over a base. This leads to the following deep conjecture.

\begin{conjecture}[Orbit Density Conjecture for $\AA^N$]
  \textup{(Medvedev--Scanlon, 2009,
    {\cite[Conjecture~5.10]{arxiv0901.2352}};
    Amerik--Bogomolov--Rovinsky, 2011,
    {\cite[Conjecture~1.2]{MR2862064}})}
\label{conjecture:orbdensityPN}
Let $f:\AA^N\to\AA^N$ be a morphism defined over~$\Qbar$, i.e.,
$f=(f_1,\ldots,f_N)$, where the~$f_i\in\Qbar[x_1,\ldots,x_N]$ are
polynomials with~$\Qbar$-coefficients. Assume that there is no
rational function $g\in\Qbar(x_1,\ldots,x_N)$ satisfying $g\circ
f=g$. Then there exists a point~$P\in\PP^N(\Qbar)$ whose
orbit~$\Orbit_f(P)$ is Zariski dense in~$\PP^N$.
\end{conjecture}

\begin{theorem}
\textup{(Xie, 2017, \cite{MR3705271})}
Conjecture~$\ref{conjecture:orbdensityPN}$ is true for~$N=2$.
\end{theorem}

Theorem~\ref{theorem:Qbarwander} and
Conjecture~\ref{conjecture:orbdensityPN} say that there are a lot of
points with coordinates in~$\Qbar$ having certain properties. The
field~$\Qbar$ is of great arithmetic interest, but frequently one
would like to find points defined over a specific number field.

\begin{question}
To what extent are Theorem~$\ref{theorem:Qbarwander}$ and
Conjecture~$\ref{conjecture:orbdensityPN}$ valid if one replaces~$\Qbar$
with a number field~$K$?
\end{question}

\subsection{Generalization to Algebraic Varieties}

Amerik's theorem (Theorem~\ref{theorem:Qbarwander}) is true more
generally for arbitrary varieties, as is the orbit density conjecture,
for which we now give a general formulation.

\begin{conjecture}[Orbit Density Conjecture]
  \textup{(Medvedev and Scanlon, 2009,
    {\cite[Conjecture~5.10]{arxiv0901.2352}}; Amerik, Bogomolov and
    Rovinsky, 2011, {\cite[Conjecture~1.2]{MR2862064}})}
\label{conjecture:orbdensity}
Let $\Kbar$ be an algebraically closed field of characteristic~$0$,
let~$X/K$ be a quasi-projective variety, and let~$f:X\to X$ be a
dominant morphism. Assume that there does not exist a non-constant
rational function~$g\in \Kbar(X)$ satisfying $g\circ f=g$.  Then there
exists a point~$P\in X(\Kbar)$ whose orbit~$\Orbit_f(P)$ is Zariski
dense in~$X$.
\end{conjecture}

%%%%%%%%%%%%%%%%%%%%%%%%%%%%%%%%%%%%%%%%%%%%%%%%%%%%%%%%%%%%%%%%%%%%%%
\section{Dynamical Degrees of Rational Maps\WhoDone{Joe}}
\label{section:ddrm}
\WhoWrite{Joe}
%%%%%%%%%%%%%%%%%%%%%%%%%%%%%%%%%%%%%%%%%%%%%%%%%%%%%%%%%%%%%%%%%%%%%%

We recall that a rational map $f:\PP^N\dashrightarrow\PP^N$ is
described by a list of homogeneous polynomials $f_0,f_1,\ldots,f_N\in
K[X_0,\ldots,X_N]$ of the same degree, with the proviso that any
common factors of~$f_0,\ldots,f_N$ are removed. The \emph{degree
  of~$f$} is then the common degree of~$f_0,\ldots,f_N$, and~$f$ is
\emph{dominant} if~$f_0,\ldots,f_N$ do not themselves satisfy a
non-trivial polynomial relation.

In general, if~$\deg(f)=d$, then one might expect that the degree of
the $n$th iterate~$f^n$ would satisfy $\deg(f^n)=d^n$.  This is
true if~$f$ is a morphism, and in any case, we always have an
inequality $\deg(f^n)\le d^n$, but cancellation of common factors
may cause the degree of~$f^n$ to be smaller than expected.

\begin{example}
\label{example:Fibdegmap}
The degree~$2$  rational map $f:\PP^2\dashrightarrow\PP^2$ given by
the formula $f(X,Y,Z) = [YZ,XY,Z^2]$ has second iterate
\[
  f^2(X,Y,Z) = f(YZ,XY,Z^2) = [XYZ^2,XY^2Z,Z^4] = [XYZ,XY^2,Z^3].
\]
Thus $\deg(f^2)=3<4$, since we have canceled a factor of~$Z$.  An easy
induction shows that $\deg(f^n)$ equals the $(n+2)$nd Fibonacci
number.
%% Similarly, we find that $f^3(X,Y,Z) = [XY^2Z^2,X^2Y^3,Z^5]$, so $\deg(f^3)=5$.
\end{example}

It is an interesting problem to determine the degree
sequence~$\deg(f^n)$ of a rational map. The average growth is measured
by the following quantity, which was defined independently by
Arnold~\cite{MR1139553}, Bellon--Viallet~\cite{MR1704282}, and
Russakovskii--Shiffman~\cite{MR1488341}.

\begin{definition}
Let $f:\PP^N\dashrightarrow\PP^N$ be a dominant rational map. The
\emph{dynamical degree of~$f$} is
\[
  \d_f := \lim_{n\to\infty} \Bigl( \deg(f^n) \Bigr)^{1/n}.
\]
The quantity $\log(\d_f)$ is sometimes called the \emph{algebraic
  entropy of~$f$.}
\end{definition}

For example, using the classical closed formula for Fibonacci numbers,
we see that the map in Example~\ref{example:Fibdegmap} has dynamical
degree~$\d_f=\frac12(1+\sqrt5)$, the so-called Golden Ratio. It is a
nice exercise, using the inequality
$\deg(f^{n+m})\le\deg(f^n)\deg(f^m)$, to prove that the limit
defining~$\d_f$ converges.  Note that we have $\deg(f^n)\approx\d_f^n$
when~$n$ is large.

\begin{remark}
We recall that a rational map
$f=[f_0,f_1,\ldots,f_N]:\PP^N\dashrightarrow\PP^N$ is a
\emph{morphism} if there are no solutions to $f_0=f_1=\cdots=f_N=0$
in~$\PP^N(\Kbar)$. In this case, it is not hard to see that
$\deg(f^n)=\deg(f)^n$ for all~$n$, and hence~$\d_f=\deg(f)$.  More
generally, a rational map~$f$ is said to be \emph{algebraically
  stable} if $\deg(f^n)=\deg(f)^n$ for all~$n$, or equivalently
if~$\d_f=\deg(f)$.
\end{remark}

Based on examples such as Example~\ref{example:Fibdegmap},
Bellon--Viallet conjectured that the sequence of degrees $\deg(f^n)$
for $n=1,2,3,\ldots$ should satisfy a linear recursion with constant
coefficients. This turns out to be false.

\begin{theorem}
\label{theorem:hasselblattpropp}
\textup{(Hasselblatt--Propp, 2007, \cite{MR2358970}; Bedford--Kim,
2008, \cite{MR2449533})} For every $N\ge3$ there exist rational
maps $f:\PP^N\dashrightarrow\PP^N$ for which the sequence of degrees
$\deg(f^n)$ does not satisfy a linear recursion with constant
coefficients.
\end{theorem}

\begin{example}
Hasselblatt and Propp provide the explicit example
\[
  f : \PP^3\dashrightarrow\PP^3,\quad f(X,Y,Z,W) = [YW,ZW,X^2,XW],
\]
whose degree sequence $2,3,4,6,9,12,17,25,\ldots$ does not satisfy a
linear recurrence with constant coefficients. However, the
construction of Hasselblatt and Propp does not provide a negative
answer to the following question.
\end{example}

\begin{question}
Let $f:\PP^N\dashrightarrow\PP^N$ be a dominant rational map.  Do
there exist a finite collection of constant coefficient linear
recursions
\[
  L_n^{(1)},L_n^{(2)},\ldots,L_n^{(r)}
\]
and an index function $\iota:\NN\to\{1,2,\ldots,r\}$ so that for all $n\ge1$ we
have
\[
%% \deg(f^n) = L_n^{(\iota(n))} ?
\deg(f^n) = L_n^{(\iota_n)} ?
\]
\end{question}

\begin{remark}
Complementing Theorem~\ref{theorem:hasselblattpropp}, there are
non-trivial situations for which it is known that the degree sequence
does satisfy a linear recursion. This is true, for example, for
birational maps of surfaces~\cite[Diller--Favre (2001)]{MR1867314} and
for polynomial endomorphsims of~$\AA^2$~\cite[Favre--Jonsson
  (2011)]{arxiv0711.2770}.
\end{remark}

And although the original Bellon--Viallet conjecture that the degree
sequence is a linear recursion turns out to be false, the following
intriguing weaker conjecture is still open.

\begin{conjecture}[Integrality Conjecture]  
\label{conjecture:dfalgint}
\textup{(Bellon--Viallet, \cite{MR1704282}, 1999)} Let
$f:\PP^N\dashrightarrow\PP^N$ be a dominant rational map. Then the dynamical
degree~$\d_f$ of~$f$ is an algebraic integer, i.e.,~$\d_f$ is the root
of a monic polynomial having integer coefficients.
\end{conjecture}

For endomorphisms of~$\PP^2$, Conjecture~\ref{conjecture:dfalgint} is
known when~$f$ is a birational
map~\cite[Diller--Favre~(2001)]{MR1867314} and when~$f$ is a
polynomial map~\cite[Favre--Jons\-son~(2007)]{MR2339287}.

Suppose that the polynomials defining~$f:\PP^N\dashrightarrow\PP^N$
have integer coefficients. Then for each prime~$p$, we can reduce the
coefficients modulo~$p$ to create a rational
map~$f_p:\PP^N\dashrightarrow\PP^N$ defined over the finite
field~$\FF_p$. For all but finitely many primes~$p$, the map~$f_p$
will still be dominant, so we can compute its dynamical
degree~$\d_{f_p}$.  We clearly have $\d_{f_p}\le\d_f$, but the
inequality may be strict, since reducing modulo~$p$ may remove some
terms that allow us to cancel an extra factor. For example,
$f(X,Y)=[3X^2+XY,Y^2]$ has degree~$2$, but~$f_3(X,Y)=[XY,Y^2]=[X,Y]$
has degree~$1$.

\begin{theorem}
\label{theorem:xiedptod}  
\textup{(Xie, 2015, \cite{arxiv1106.1825})}
Let $f:\PP^2\dashrightarrow\PP^2$ be a dominant rational map defined
by polynomials with integer coefficients. Then
\begin{equation}
  \label{eqn:dptod}
  \lim_{p\to\infty} \d_{f_p}= \d_f.
\end{equation}
\end{theorem}

Lacking evidence, we pose the following as a question,
rather than a conjecture.

\begin{question}
Does~\eqref{eqn:dptod} hold for dominant rational maps
$f:\PP^N\dashrightarrow\PP^N$ when $N\ge3$?
\end{question}

\begin{example}
Xie~\cite[Section~5]{arxiv1106.1825} provides the following example to
show that the limit formula~\eqref{eqn:dptod} cannot be strengthened
to the statement that~$\d_{f_p}=\d_f$ for all sufficiently
large~$p$. Let
\[
  f  = [XY, XY - 2Z^2, YZ + 3Z^2].
\]
Then Xie proves that~$\d_f=2$, but $\d_{f_p}<2$ for all primes~$p$.
More precisely, for~$p>2$, he shows that~$\d_{f_p}$ is the largest
real root of the polynomial~$\l^{\e(p)}-2\l^{{\e(p)}-1}+1$,
where~$\e(p)$ is the order of~$2$ in the multiplicative
group~$\FF_p^*$.
\end{example}

It is interesting to investigate how dynamical degrees vary in 
families of maps.

\begin{example}
Consider the $3$-parameter family of rational maps
\[
  f_{a,b,c}\bigl([X,Y,Z]\bigr) = [XY,XY+aZ^2,bYZ+cZ^2].
\]
For most choices of~$a,b,c$, we have~$\d_{f_{a,b,c}}=2$, but there are
exceptions. It is proven in~\cite{arxiv1609.02119} that
for~$a,b,c\in\CC$, we have~$\d_{f_{a,b,c}}<2$ if and only if there is
a root of unity~$\xi\in\CC$ such that $c^2=(\xi+\xi^{-1})^2ab$.
\end{example}

\begin{conjecture}
\label{conjecture:dftledfe}
Let $f_T:\PP^N\to\PP^N$ be a family of dominant rational maps, depending
rationally on a list of parameters $T=(T_1,\ldots,T_k)$. Then for all $\e>0$, the set
\begin{equation}
  \label{eqn:tCkdftledfe}
  \bigl\{ t\in\CC^k  : \d(f_t)\le \d(f)-\e \bigr\}
\end{equation}
is not Zariski dense in~$\CC^k$. In particular, for one-parameter
families, i.e., ~$k=1$, the set~\eqref{eqn:tCkdftledfe} is finite.
\end{conjecture}

One approach to proving Conjecture~$\ref{conjecture:dftledfe}$ would
be to prove the following uniform lower bound for~$\d_f$, from which
it is easy to deduce Conjecture~$\ref{conjecture:dftledfe}$;
see~\cite{arxiv1609.02119,arxiv1106.1825}.

\begin{conjecture}
\label{conjecture:dfgedfk}
Let $N\ge1$. There exists a constant~$\g_N>0$ such that for all
dominant rational maps $f:\PP^N\dashrightarrow\PP^N$ we have
\[
  \d(f) \ge \g_N \cdot \min_{0\le n<N} \frac{\deg(f^{n+1})}{\deg(f^n)}.
\]
N.B. The constant~$\g_N$ is not allowed to depend on~$f$.
\end{conjecture}

\begin{theorem}
\textup{(Xie, 2015, \cite{arxiv1106.1825})}
For~$N=2$, Conjecture~$\ref{conjecture:dfgedfk}$, and thus also
Conjecture~$\ref{conjecture:dftledfe}$, are true.
\end{theorem}

The definition of~$\d_f$ says that the degree sequence~$\deg(f^n)$
grows roughly like~$\d_f^n$, but this may not capture the true growth
rate. For example, the map $f(X,Y,Z)=[X^dY,Y^dZ,Z^{d+1}]$ has $\d_f=d$
and $\deg(f^n)\sim d^{-1}n\d_f^n$. This leads to a refined degree
question, posed in varying forms in~\cite{MR2358970,arxiv1111.5664}. 

\begin{question}
\label{question:degfnnelldn}
Let $f:\PP^N\dashrightarrow\PP^N$ be a dominant rational map.  Is it
true that there exists an integer~$\ell_f$ satisfying $0\le\ell_f\le
N$ such that\footnote{We recall that the notation $F(n)\asymp G(n)$
  means that there are positive constants~$c_1,c_2,c_3$, depending on~$F$
  and~$G$, so that $c_1F(n) \le G(n) \le c_2F(n)$ for all~$n\ge c_3$.}
\[
  \deg(f^n) \asymp n^{\ell_f} \d_f^n,
\]
where the implied constants may depend on~$f$, but not on~$n$?
\end{question}

Question~\ref{question:degfnnelldn} has an affirmative answer for
monomial maps~\cite{arxiv1011.2854,arxiv1007.0253}, and for maps of
surfaces if~$f$ satisfies $(\d_f^{(1)})^2>\d_f^{(2)}$
\cite[Boucksom--Favre--Jonsson (2008)]{arxiv0608267}.  If~$f$ is a
birational map, then it is known that the sequence of
degrees~$\deg(f^n)$ is bounded if and only if there is a subsequence
of~$n$ along which it is bounded, and that if~$\deg(f^n)$ is
unbounded, then there is a uniform lower bound for its growth
rate~\cite[Cantat--Xie (2018)]{arxiv1802.08470}.  (This result on
birational maps holds more generally for any variety.)

\subsection{Generalization to Algebraic Varieties}
Let~$X$ be a smooth projective variety of dimension~$N$, and let
$f:X\dashrightarrow X$ be a dominant rational map. Then~$f$ induces a
pullback map~$f^*$ on the group~$\Pic(X)$ of line bundles on~$X$, but
in general it is not true that $(f^n)^*$ is the $n$th iterate
of~$f^*$. Letting~$\Lcal$ be any ample sheaf on~$X$, the
\emph{dynamical degree of~$f$} is defined using intersection theory by
the limit
\[
  \d_f := \lim_{n\to\infty} \Bigl( (f^n)^*\Lcal \cdot \Lcal^{N-1}\Bigr)^{1/n}.
\]
More generally, this~$\d_f$ is the first dynamical degree, and for
$1\le i\le N$, we define the \emph{$i$th dynamical degree of~$f$} by
\begin{equation}
  \label{eqn:highorderdyndeg}
  \d_f^{(i)} := \lim_{n\to\infty} \Bigl( \bigl( (f^n)^*\Lcal\bigr)^i \cdot \Lcal^{N-i}\Bigr)^{1/n}.
\end{equation}
It was noted by Guedj~\cite{MR2179389} that the the sequence of
dynamical degrees is concave, i.e., there is a~$k$ such that
\[
  \d_f^{(1)} \le \d_f^{(2)} \le \cdots \le \d_f^{(k)}
  \ge \d_f^{(k+1)} \ge \cdots \ge \d_f^{(N)}.
\]
For a description of the full sequence of dynamical degrees when~$f$
is a monomial map, see \cite[Favre--Wulcan~(2012)]{arxiv1011.2854} and
\cite[Lin~(2012)]{arxiv1010.6285}.

\begin{question}
Let~$X$ be a smooth projective variety of dimension~$N$, and let
$f:X\dashrightarrow X$ be a dominant rational map. Is it true
that the higher order dynamical degrees~$\d_f^{(i)}$ are algebraic
integers for all~$1\le i\le N$\textup{?}
\end{question}

%%%%%%%%%%%%%%%%%%%%%%%%%%%%%%%%%%%%%%%%%%%%%%%%%%%%%%%%%%%%%%%%%%%%%%
\section{Arithmetic Degrees of Orbits\WhoDone{Joe}}
\label{section:ado}
\WhoWrite{Joe}
%%%%%%%%%%%%%%%%%%%%%%%%%%%%%%%%%%%%%%%%%%%%%%%%%%%%%%%%%%%%%%%%%%%%%%
For a fixed dimension~$N$, the degree of a rational
map~$f:\PP^N\dashrightarrow\PP^N$ is a measure of the complexity
of~$f$, since for example, it takes on the order of~$d^N$ coefficients to
specify~$f$.  In an analogous manner, the complexity of a list of
integers is measured roughly by the number of bits it takes to specify
the integers. The mathematical term for this number theoretic
complexity is \emph{height}, which we now describe more precisely.
For further details and proofs, see for
example~\cite{hindrysilverman:diophantinegeometry,lang:diophantinegeometry,MR2316407,MR2514094}.

\begin{definition}
Let $P=[a_0,a_1,\ldots,a_N]\in\PP^N(\QQ)$, where we choose homogeneous
coordinates for~$P$ satisfying~$a_i\in\ZZ$ and
$\gcd(a_0,\ldots,a_N)=1$. Then the \emph{height of~$P$} is
\[
  h(P) := \log \max \bigl\{ |a_0|, |a_1|, \ldots, |a_N| \bigr\}.
\]
There is a similar, but more complicated, definition for
points~$P\in\PP^N(\Qbar)$ whose coordinates are algebraic numbers. But
the intuition is that it takes on the order of~$h(P)$ bits to describe
the point~$P$.
\end{definition}

It is clear that for any~$B$, there are only finitely
many~$P\in\PP^N(\QQ)$ satisfying~$h(P)\le B$, and the same is true
of~$\PP^N(K)$ for any number field~$K$, and even for the union over
all~$K$ of a fixed degree over~$\QQ$. Further, for a rational map
$f:\PP^N\dashrightarrow\PP^N$ of degree~$d$, an elementary triangle
inequality argument shows that there is a constant~$A_f$ such that
\begin{equation}
  \label{eqn:hfPledhPO1}
  h\bigl(f(P)\bigr) \le d\cdot h(P) + A_f
  \quad\text{for all $P\in\PP^N(\Qbar)$ with $f(P)$ defined.}
\end{equation}
If~$f$ is a morphism, i.e., defined everywhere, a more elaborate
argument using the Nullstellensatz can be used to prove that there is
a complementary lower bound with another constant~$B_f$,
\begin{equation}
  \label{eqn:hfPgedhPO1}
  h\bigl(f(P)\bigr) \ge d\cdot h(P) - B_f
  \quad\text{for all $P\in\PP^N(\Qbar)$.}
\end{equation}

Applying~\eqref{eqn:hfPledhPO1} repeatedly, we see that the arithmetic
complexity, i.e, the height of~$f^n(P)$, is never much more than a
multiple of~$d^n$. By analogy with the motivation used to define
dynamical degree, we measure the average growth rate of the arithmetic
complexity of~$f^n(P)$ by the following quantity, whose logarithm is a sort of
\emph{arithmetic entropy} of the orbit~$\Orbit_f(P)$.

\begin{definition}
Let $f:\PP^N(\Qbar)\dashrightarrow\PP^N(\Qbar)$ be a dominant rational
map. We recall that~$\PP^N(\Qbar)_f$ denotes the set of
points~$P\in\PP^N(\Qbar)$ for which the entire forward
orbit~$\Orbit_f(P)$ is well-defined. Then the \emph{arithmetic degree
  of the $f$-orbit of a point~$P\in\PP^N(\Qbar)_f$} was defined
in~\cite{arxiv1111.5664} to be the limit
\begin{equation}
  \label{eqn:afPdef}
  \a_f(P) := \lim_{n\to\infty} h\bigl(f^n(P)\bigr)^{1/n}.
\end{equation}
Since it is not yet known in general that this limit exists, we
also define the \emph{upper arithmetic degree} using the limsup,
\[
  \bar\a_f(P) := \limsup_{n\to\infty} h\bigl(f^n(P)\bigr)^{1/n}.
\]
\end{definition}

We note that Theorem~\ref{theorem:Qbarwander} ensures that
$\PP^N(\Qbar)_f$ contains lots of wandering points, and for maps~$f$
that don't fiber, Conjecture~\ref{conjecture:orbdensity} says that
there will even be lots of points in $\PP^N(\Qbar)_f$ having Zariski
dense orbits.

We know that~$f^n$ is given by polynomials of degree roughly~$\d_f^n$.
Combined with~\eqref{eqn:hfPledhPO1}, this suggests
that~$h\bigl(f^n(P)\bigr)$ should be no larger than some fixed
multiple of~$\d_f^n$. This is the intuition behind the following
inequality, which is proven in~\cite{arxiv1111.5664}:
\begin{equation}
  \label{eqn:afPledf}
  \bar\a_f(P) \le \d_f.
\end{equation}
(But see Section~\ref{subsection:algvar:arithdeg} for the more
difficult case when~$\PP^N$ is replaced by a general variety.)  In
words, the inequality~\eqref{eqn:afPledf} says that the arithmetic
complexity of an orbit is at most equal to the geometric complexity of
the map. Part~(d) of the next conjecture provides a natural geometric
condition that implies equality.

\begin{conjecture}[Arithmetic Degree Conjectures (on $\PP^N$)]
\label{conjecture:afPdfconjs}  
\textup{(Kawaguchi--Silverman, 2006, \cite{kawsilarithcompl2006}; see also \cite{MR3456169,MR3483624})}
Let $f:\PP^N\dashrightarrow\PP^N$ be a dominant rational map defined over~$\Qbar$.
\begin{parts}
\Part{(a)}
Let $P\in\PP^N(\Qbar)_f$.
The limit~\eqref{eqn:afPdef} defining~$\a_f(P)$ converges.
\Part{(b)}
The set $\bigl\{ \a_f(P) : P\in\PP^N(\Qbar)_f \bigr\}$ is finite.
\Part{(c)}
Let $P\in\PP^N(\Qbar)_f$. The arithmetic degree~$\a_f(P)$ is an algebraic integer.
\Part{(d)}  
Let $P\in\PP^N(\Qbar)_f$ be a point whose orbit~$\Orbit_f(P)$ is Zariski dense in~$\PP^N$.
Then $\a_f(P)=\d_f$.
\end{parts}
\end{conjecture}

Conjecture~\ref{conjecture:afPdfconjs} is known for monomial maps,
which are maps $f=[f_0,\ldots,f_N]$ such that their coordinate
functions~$f_i$ are monomials~\cite{arxiv1111.5664}, and parts of the
conjecture are known in various other cases; see for
example~\cite{MR3091603,MR3189467,MR3456169,MR3483624}.

We may ask for more refined information regarding the possible growth
rates of~$h\bigl(f^n(P)\bigr)$. The following is an arithmetic
analogue of Question~\ref{question:degfnnelldn}.

\begin{question}
\label{question:hfnPasymp}
Let $f:\PP^N\dashrightarrow\PP^N$ be a dominant rational map defined
over~$\Qbar$, and let $P\in\PP^N(\Qbar)_f$ be a point whose
orbit~$\Orbit_f(P)$ is Zariski dense in~$\PP^N$.  Is it true that
there exist integers~$k_f$ and ~$\ell_f$ with $0\le\ell_f\le N$ such
that
\begin{equation}
  \label{eqn:hfnPdfnnellflognkf}
  h\bigl(f^n(P)\bigr) \asymp  \d_f^n\cdot n^{\ell_f}\cdot (\log n)^{k_f},
\end{equation}
where the implied constants depend on~$f$ and~$P$, but not on~$n$?
Further, if~$\d_f>1$, is it true that~$k_f=0$?
\end{question}

If $f:\PP^N\to\PP^N$ is a morphism,
then~\eqref{eqn:hfnPdfnnellflognkf} holds in the stronger form
$h\bigl(f^n(P)\bigr)=\d_f^n\hhat_f(P)+O(1)$; see
Section~\ref{section:ch}.  The situation is more complicated for
general varieties.  Sano~\cite{arxiv1801.02831} proves that if~$f:X\to
X$ is a morphism of a smooth projective variety defined over~$\Qbar$,
and if~$P\in X(\Qbar)$ satisfies~$\a_f(P)>1$, then a version
of~\eqref{eqn:hfnPdfnnellflognkf} holds with~$k_f=0$.

\begin{example}
To see that the $\log$ term in Question~\ref{question:hfnPasymp} may
be necessary, consider the map $f:\PP^4\dashrightarrow\PP^4$ given
by
\[
  f(X,Y,Z,W)=\bigl[ X(Y+Z),W(Y+Z), WZ, W^2]\quad\text{and}\quad P = [1,0,1,1].
\]
One can compute $\deg(f^n)=n+1$, so $\d_f=1$ and $\ell_f=1$, while
$f^n(P)=[n!,n,1,1]$, so $h\bigl(f^n(P)\bigr)\sim
n\log(n)$. See~\cite[Example~17]{arxiv1111.5664} for details.
\end{example}

We next seek to quantify the amount by which a rational map may
decrease the height of a point.  Let $f:\PP^N\dashrightarrow\PP^N$ be
a dominant rational map defined over~$\Qbar$. The \emph{height
  expansion ratio of~$f$} is
\[
  \mu(f) := \sup_{\emptyset\ne U\subset\PP^N} \liminf_{\substack{P\in
      U(\Qbar)\\ h(P)\to\infty\\}} \frac{h\bigl(f(P)\bigr)}{h(P)},
\]
where the sup is over all non-empty Zariski open subsets
of~$\PP^N$. Thus~$\mu(f)$ measures how much~$f$ can reduce the heights
of points on a large subset of~$\PP^N$. For example,
from~\eqref{eqn:hfPledhPO1} we see that $\mu(f)\le\deg(f)$,
while~\eqref{eqn:hfPgedhPO1} implies that if~$f$ is a morphism, then
$\mu(f)=\deg(f)$. It turns out that~$\mu(f)$ cannot be arbitrarily small.

\begin{theorem}
\textup{(Silverman, 2011, \cite{arxiv0908.3835})}
We have
\[
  \overline\mu_d(\PP^N) :=
  \inf_{\substack{f:\PP^N\dashrightarrow\PP^N\\\textup{dominant}\\\deg(f)=d\\}} \mu(f) > 0.
\]
\end{theorem}

It is known that $\overline\mu_d(\PP^N)\le d^{-(N-1)}$ for all $N\ge2$
and all $d\ge2$; see \cite[Proposition~10]{arxiv0908.3835}
and~\cite{MR2975153}.  This leads to the obvious question.

\begin{question}
What is the exact value of $\overline\mu_d(\PP^N)$? 
\end{question}

\subsection{Generalization to Algebraic Varieties}
\label{subsection:algvar:arithdeg}
In order to generalize Conjecture~\ref{conjecture:afPdfconjs} to a
projective variety~$X$, we need an appropriate height function.

\begin{definition}
Let $X$ be a smooth projective variety, let $i:X\hookrightarrow\PP^N$
be an embedding, and let $f:X\dashrightarrow X$ be a dominant rational
map, with everything defined over~$\Qbar$.  Then the \emph{arithmetic
  degree of the $f$-orbit of~$P\in X(\Qbar)_f$} is
\begin{equation}
  \label{eqn:afPdefgen}
  \a_f(P) := \lim_{n\to\infty} h\bigl(i\circ f^n(P)\bigr)^{1/n}.
\end{equation}
We write~$\bar\a_f(P)$ for~\eqref{eqn:afPdefgen} with the limit
replaced by limsup.
\end{definition}

It is not hard to verify that~$\a_f(P)$, if it exists, does not depend
on the embedding~$X$ in projective space. The main conjectures on
arithmetic degrees carry over verbatim to the general case.

\begin{conjecture}[Arithmetic Degree Conjectures]
\label{conjecture:afPdfconjsvariety}  
\textup{(Kawaguchi--Silverman, 2006,
  \cite{kawsilarithcompl2006,MR3456169,MR3483624})} Let
$f:X\dashrightarrow X$ be a dominant rational map of a smooth
projective variety~$X$ defined over~$\Qbar$.
\begin{parts}
\Part{(a)}
Let $P\in X(\Qbar)_f$.
The limit~\eqref{eqn:afPdefgen} defining~$\a_f(P)$ converges.
\Part{(b)}
The set $\bigl\{ \a_f(P) : P\in X(\Qbar)_f \bigr\}$ is finite.
\Part{(c)}
Let $P\in X(\Qbar)_f$. The arithmetic degree~$\a_f(P)$ is an algebraic integer.
\Part{(d)}  
Let $P\in X(\Qbar)_f$ whose orbit~$\Orbit_f(P)$ is Zariski dense
in~$X$.  Then $\a_f(P)=\d_f$.
\end{parts}
\end{conjecture}

The inequality
\begin{equation}
  \label{eqn:afPledfPXQbar}
  \bar\a_f(P) \le \d_f\quad\text{for $P\in X(\Qbar)_f$}
\end{equation}
relating arithmetic and dynamical degrees is true in the more general
setting of Conjecture~\ref{conjecture:afPdfconjsvariety}, but the
proof is no longer elementary as it was for~$X=\PP^N$, even if we
assume that~$f$ is a morphism. See~\cite{MR3456169} for a proof
of~\eqref{eqn:afPledfPXQbar} when~$f$ is a morphism,
and~\cite{arxiv1606.00598} for the general case of dominant rational
maps.

Progress on Conjecture~\ref{conjecture:afPdfconjsvariety} is
fragmentary. Parts~(a,b,c) are true when~$f$ is a
morphism~\cite{MR3456169}.  The more difficult Part~(d) has been
proven when~$f$ is a monomial map (\emph{op.\ cit.}), when~$X$ is an
abelian variety, i.e., when~$X$ has the structure of an algebraic
group~\cite{MR3456169,MR3614521}, and when~$f$ is a morphism and
either~$X$ is surface~\cite{arxiv1701.04369} or~$X$ is higher dimensional
with additional structure~\cite{arxiv1802.07388}.

%%%%%%%%%%%%%%%%%%%%%%%%%%%%%%%%%%%%%%%%%%%%%%%%%%%%%%%%%%%%%%%%%%%%%%
\section{Canonical Heights\WhoDone{Patrick+Joe}}
\label{section:ch}
\WhoWrite{Patrick + Joe}
%%%%%%%%%%%%%%%%%%%%%%%%%%%%%%%%%%%%%%%%%%%%%%%%%%%%%%%%%%%%%%%%%%%%%%

We recall from Section~\ref{section:ado} that the height~$h(P)$ of a
point~$P\in\PP^N(\Qbar)$  measures  the arithmetic complexity
of~$P$, and that a degree~$d$ morphism $f:\PP^N\to\PP^N$ has the effect,
roughly, of multiplying the height by~$d$. More precisely,
\[
  h\bigl(f(P)\bigr) = dh(P) + O(1)\quad\text{for all $P\in\PP^N(\Qbar)$,}
\]
where the big-$O$ bound depends on~$f$, but is independent of~$P$.
This suggests that~$h\bigl(f^n(P)\bigr)$ grows like~$d^n$, and in
Section~\ref{section:ado} we took the limit of~$h(f^n(P))^{1/n}$ to
measure this growth. When~$f$ is a morphism, we can do much better.

\begin{definition}
Let $f:\PP^N\to\PP^N$ be a morphism of degree~$d$ defined
over~$\Qbar$, and let $P\in\PP^N(\Qbar)$. The (\emph{dynamical})
\emph{canonical height} of~$P$ for the map~$f$ is the quantity
\begin{equation}
  \label{eqn:hhatdef}
  \hhat_f(P) := \lim_{n\to\infty} \frac{1}{d^n}h\bigl(f^n(P)\bigr).
\end{equation}
\end{definition}

The definintion of~$\hhat_f$ is modeled after Tate's construction of
the canonical (N\'eron--Tate) height on abelian varieties.
See~\cite[\S3.4]{MR2316407} for the telescoping sum proof that the
limit~\eqref{eqn:hhatdef} converges.

The dynamical canonical height has several agreeable properties:
\begin{align}
  \label{eqn:hffPdhfP}
  \hhat_f\bigl(f(P)\bigr)&=d\hhat_f(P).  \\
  \label{hfPhPO1}
  \hhat_f(P) &= h(P) + O(1).  \\
  \label{eqn:hfP0Ppp}
  \hhat_f(P)=0\quad&\Longleftrightarrow\quad P\in\PrePer(f).
\end{align}
We note that~$\hhat_f$ is uniquely determined by~\eqref{eqn:hffPdhfP} and ~\eqref{hfPhPO1}

The classical Lehmer conjecture gives a lower bound for the height of
algebraic numbers that are not roots of unity. It is the $f(x)=x^2$
case of the following conjecture, which provides a quantitative
converse to~\eqref{eqn:hfP0Ppp}.

\begin{conjecture}[Dynamical Lehmer Conjecture]
\label{conjecture:dynlehmer}
Let $f:\PP^N\to\PP^N$ be a morphism of degree~$d\ge2$ defined over a
number field~$K$. There is a constant~$C(f)>0$ such that for
all~$P\in\PP^N(\Qbar)$ that are wandering for~$f$, we have\footnote{We
  recall from Section~\ref{section:background} that~$K(P)$ is the
  field generated by the  coordinates of~$P$.}
\[
  \hhat_f(P) \ge \frac{C(f)}{\bigl[K(P):K\bigr]}.
\]
\end{conjecture}

In order to state the next conjecture, we need a way to measure the
arithmetic complexity of a rational map $f:\PP^N\to\PP^N$ defined
over~$\Qbar$. It is not enough to simply use the coefficients of the
polynomials defining~$f$, since as we saw in
Section~\ref{section:dms}, changing coordinates on~$\PP^N$ has the
effect of drastically changing the coefficients of~$f$ without
affecting the underlying dynamical system. The solution is to view~$f$
as a point~$\langle f\rangle$ in~$\Moduli_d^N$, the space of ``maps up to change-of-coordinates.''
The space~$\Moduli_d^N$ is itself an algebraic variety, so we fix an embedding
\[
  \iota : \Moduli_d^N \longhookrightarrow \PP^L,
\]
and then we define the \emph{moduli height of~$f$} to be the height of
the associated point in projective space,
\[
  h_\Moduli(f) := h\bigl(\iota\bigl( \langle f\rangle \bigr)\bigr).
\]
One can show that any two such heights are commensurable, in the sense
that if $h_\Moduli'$ corresponds to some other embedding~$\iota'$ of~$\Moduli_d^N$,
then there are positive constants~$c_1,c_2,c_3,c_4$ such that
\[
  c_1 h_\Moduli(f) - c_2 \le h_\Moduli'(f) \le c_3 h_\Moduli(f) + c_4
  \quad\text{for all $f\in\Moduli_d^N(\Qbar)$.}
\]

Informally, the next conjecture says that the arithmetic complexity of
a wandering orbit is always larger than the arithmetic complexity of
the map. It is a generalization of a conjecture that Serge Lang
formulated for elliptic curves~\cite{MR518817} and provides a
different quantitative converse to~\eqref{eqn:hfP0Ppp}.

\begin{conjecture}[Dynamical Lang Height Conjecture]
\label{conjecture:dynlanght}  
Let $N\ge1$, $D\ge1$, and $d\ge2$. There are positive constants
$c_1(D,N,d)$ and $c_2(D,N,d)$ such that for all number fields~$K/\QQ$
with~$[K:\QQ]\le D$, all degree~$d$ morphisms $f:\PP^N\to\PP^N$
defined over~$K$, and all points $P\in\PP^N(K)$ whose
orbit~$\Orbit_f(P)$ is Zariski dense in~$\PP^N$, we have
\[
  \hhat_f(P) \ge c_1(D,N,d) h_\Moduli(f) + c_2(D,N,d).
\]
\end{conjecture}

See~\cite{MR2501344,MR3277939} for numerical evidence supporting
Conjecture~\ref{conjecture:dynlanght} in the case $N=D=1$ and
$d\in\{2,3\}$, including conjectural values for the constants. We
also note that the Zariski density assumption in
Conjecture~\ref{conjecture:dynlanght} is necessary. For example, the
map~$f=[X^2,Y^2,\g Z^2]$ has height roughly~$h(\g)$, while the
$f$-orbit of the point~$P=[\a,1,0]$ has canonical height
$\hhat_f(P)=h(\a)$, so we can make~$h_\Moduli(f)$ arbitrarily large
while~$\hhat_f(P)$ remains bounded.  Ingram~\cite{MR2504779} and
Looper~\cite{arxiv1709.08121} have proven weak versions of
Conjecture~\ref{conjecture:dynlanght} for polynomial maps~$f$
of~$\PP^1$, where the constants depend also on the number of primes of
bad reduction of~$f$.

The dynamics of a rational map $f:\PP^1\to\PP^1$ is greatly influenced
by the orbits of its critical points. We
write~$\Crit(f)$ for the set of critical points of~$f$ in~$\PP^1$,
counted with appropriate multiplicities.

\begin{definition}
Let $f:\PP^1\to\PP^1$ be a rational map of degree at least~$2$ defined
over~$\Qbar$.  The \emph{critical height} of~$f$ is
\[
  \hhat^\crit(f) := \sum_{P\in\Crit(f)} \hhat_f(P).
\]
\end{definition}

We observe from~\eqref{eqn:hfP0Ppp} that $\hhat^\crit(f)=0$ if and
only if every critical point of~$f$ is preperiodic, i.e., if and only
if~$f$ is post-critically finte~(PCF).  In order to state the next
result, which provides an interpretation of~$\hhat^\crit$ and was
conjectured  by Silverman circa~2010, we need a somewhat technical
definition. However, the only part of the definition that is required
in order to understand Theorem~\ref{theorem:hcrithmod} is that most
maps are not Latt\`es maps.

\begin{definition}
For each integer~$m$ and each elliptic curve~$E/\CC$, we construct a
rational map~$f_{E,m}$ by choosing a Weierstrass equation
$y^2=x^3+Ax+B$ for~$E$ and requiring~$f_{E,m}$ to satisfy the equation
\[
  f_{E,m}\bigl(x(P)\bigr) = x(mP)\quad\text{for all $P\in E(\CC)$.}
\]
The map~$f_{E,m}$ has degree~$m^2$ and is called a (flexible)
\emph{Latt\`es map}. The collection of all Latt\`es maps
in~$\Moduli_{m^2}^1$ forms a 1-dimensional family, so most maps are
not Latt\`es maps, since $\dim\Moduli_d^1=2d-2$. Latt\`es
maps have many special properties; in particular, every Latt\`es map
is PCF. For further information about Latt\`es maps, see for
example~\cite{milnor:lattesmaps} or~\cite[Chapter~6]{MR2316407}.
\end{definition}

\begin{theorem}
  \label{theorem:hcrithmod}
  \textup{(Ingram, 2016, \cite{arxiv1610.07904})}
Let $d\ge2$. There are positive constants $c_1,c_2,c_3,c_4$, depending
only on~$d$, such that for all non-Latt\`es degree~$d$ rational maps
$f:\PP^1\to\PP^1$ defined over~$\Qbar$, we have
\[
  c_1 h_\Moduli(f) - c_2 \le \hhat^\crit(f) \le c_3 h_\Moduli(f) + c_4.
\]
\end{theorem}

We remark that the easier upper bound part of
Theorem~\ref{theorem:hcrithmod} had been proven earlier by
Silverman~\cite[Theorem~6.31]{MR2884382}.  We also note that an
important corollary of Theorem~\ref{theorem:hcrithmod} is that outside
of the Latt\`es locus, the set of~PCF maps in~$\Moduli_d^1$ is a set
of bounded height. This weaker result had been proven earlier by
Benedetto, Ingram, Jones, and Levy~\cite{MR3265554}.  See
also~\cite{arxiv1803.06859} for a recent quantitative version of
Theorem~\ref{theorem:hcrithmod}.

\begin{question}
Is there a natural generalization of Theorem~$\ref{theorem:hcrithmod}$
for morphisms $f:\PP^N\to\PP^N$\textup{?} The moduli height~$h_\Moduli$ was
defined earlier for all~$\Moduli_d^N$, but it is not at all clear what
should take the place of the critical height when~$N\ge2$.
\end{question}

Conjecture~\ref{conjecture:dynlanght} asks for the smallest positive
value of the canonical height~$\hhat_f(P)$ as~$f$ varies over all
degree~$d$ morphisms $f:\PP^N\to \PP^N$ and~$P$ varies over all
$f$-wandering points.  Questions about the variation of the canonical
height become easier if we require that the maps and points vary in an
algebraic family. Thus let $B$ be a quasi-projective variety, let
$f_t:\PP^N\to \PP^N$ be a family of morphisms parametrized by
$t\in{B}$, and let $P:B\to \PP^N$ be a morphism, so we may view~$P_t$
as being a family of points parametrized by~$B$. It is then natural to
ask about the behavior of the function
\[
  B(\Kbar)\longrightarrow\RR,\quad
  t\longmapsto \hhat_{f_t}(P_t),
\]
or to relate it to other natural functions on $B$.

Alternatively, we may view~$f$ as an endomorphism of $\PP^N$ defined
over the function field~$K(B)$ of~$B$, and then $P$ is a point in
$\PP^N\bigl(K(B)\bigr)$.  In the case that $B$ is a curve, each
element of $\b\in K(B)$ defines a map $\b:B\to\PP^1$. We define the
height~$h(\b)$ to be the degree of the map~$\b$. We can similarly
define the height of a point~$P\in\PP^N\bigl(K(B)\bigr)$, and by
taking limits, we obtain a canonical height~$\hhat_f(P)$ for~$P$
relative to the map~$f$.  We also define a normalized height
on~$B(\Qbar)$ by fixing a non-constant rational map $\b:B\to\PP^1$ and
defining
\[
  h_B:B(\Qbar)\longrightarrow\RR,\quad
  h_B(t) := \frac{1}{\deg\b} h\bigl(\b(t)\bigr).
\]
More generally, for every divisor~$D\in\Div(B)\otimes\RR$, one can
define an associated Weil height function $h_{B,D}:B(\Qbar)\to\RR$,
well-defined up to~$O(1)$;
see~\cite{hindrysilverman:diophantinegeometry,lang:diophantinegeometry}
for details.

Generalizing an earlier result of Silverman for families of abelian
varieties, Call and Silverman proved the following result relating
three different height functions associated to a 1-dimensional family
of maps and points.

\begin{theorem}[Call--Silverman, 1993, {\cite[Theorem~4.1]{callsilv:htonvariety}}]
\label{theorem:callsilv}
Let $B$ be a smooth curve, and let $f:\PP^N\to\PP^N$ and $P:B\to\PP^N$ be a family of
maps and a family of points parametrized by~$B$. Then
\[
  \lim_{\substack{t\in B(\Qbar)\\ h_{B}(t)\to \infty\\}}\frac{\hhat_{f_t}(P_t)}{h_{B,}(t)}=\hhat_{f}(P).
\]
\end{theorem}

Theorem~\ref{theorem:callsilv} is an asymptotic formula for a ratio of
heights. To what extent can it be strengthened?  If $f$ is a family
of Latt\`{e}s maps, Tate~\cite{MR692114} proved that one can associate to~$f$ and~$P$ a
divisor $E_{f,P}\in \Div(B)\otimes\QQ$ such that
\begin{equation}
  \label{eqn:tatehftpt}
  \hhat_{f_t}(P_t)=h_{B, E_{f,P}}(t)+O(1).
\end{equation}
Thus the map $t\mapsto\hhat_{f_t}(P_t)$ is a Weil height on the base
curve~$B$. In particular, if~$B\cong\PP^1$, then any two divisors of
the same degree give the same height, up to~$O(1)$, so in this case
Tate's result for Latt\`es maps implies that
\[
  \hhat_{f_t}(P_t)=\hhat_{f}(P)h(t)+O(1).
\]
For an arbitrary base curve~$B$, there is a weaker estimate in which
the~$O(1)$ is replaced by $O(\sqrt{h_B(t)})$.

In general, let $f_t:\PP^N\to\PP^N$ be a family of morphisms
parametrized by a non-singular variety~$B$, and let $P:B\to\PP^N$ be a
morphism giving a family of points parametrized by~$B$. We define a
divisor associated to~$f$ and~$P$ by choosing a hyperplane~$H$
in~$\PP^N$ and taking the limit
\begin{equation}
  \label{eqn:EfPdef}
   E_{f,P} := \lim_{n\to\infty} \frac{1}{d^n}(f^n\circ P)^*H \in \Div(B)\otimes\RR.
\end{equation}
The limit should exist, with the height~$h_{B,E}$ on~$B$ associated
to~$E$ not depending on the choice of~$H$, up to the usual~$O(1)$.

\begin{conjecture}
\label{conj:variation}
Let $f_t:\PP^N\to\PP^N$ be a family of morphisms
parametrized by a non-singular curve~$B$,  let $P:B\to\PP^N$ be a
morphism, and let $E_{f,B}$ be the divisor defined by~\eqref{eqn:EfPdef}.
Then
\[
  \hhat_{f_t}(P_t)=h_{B, E_{f,B}}(t)+O(1) \quad\text{as $t$ ranges over $B(\Kbar)$.}
\]
\end{conjecture}
 
More generally, one might extend Conjecture~\ref{conj:variation}
to a base~$B$ of higher dimension, noting that there is no guarantee
that~$E_{f,B}$ will be ample, and indeed, in the case of constant
families it is the zero divisor.

Some progress on Conjecture~\ref{conj:variation} has been made,
including for the following type of map.

\begin{definition}
A morphism $f:\PP^N\to\PP^N$ is a \emph{regular polynomial
  endomorphism} (RPE) if there is a hyperplane $H\subset\PP^N$ satisfying
$f^{-1}(H)=H$.  Equivalently, let $f_0:\AA^N\to \AA^N$ be a polynomial
map whose extension~$f:\PP^N\to\PP^N$ is a morphism. Then~$f$ is an
RPE.
\end{definition}

A hyperplane~$H\subset\PP^N$ is isomorphic to~$\PP^{N-1}$, so
restricting an RPE~$f$ to its invariant hyperplane~$H$ gives a
morphism $\tilde f:\PP^{N-1}\to\PP^{N-1}$, and thus defines a point
$\langle\tilde f\rangle\in\Moduli_d^{N-1}$. We say that a family
$f_t:\PP^N\to\PP^N$ of RPEs parametrized by $t\in B$ is \emph{fibral}
if $\langle\tilde f_t\rangle$ is constant as~$t$ varies over the
parameter space~$B$. For example, all polynomials~$f(x)\in\Kbar[x]$
are RPEs of~$\PP^1$, and in this case one can show that every family
of RPEs is fibral.

\begin{theorem}[Ingram~\cite{MR3181564,MR3436155}]
\label{th:variation}
Conjecture~$\ref{conj:variation}$ is true for fibral families of
regular polynomial endomorphisms of~$\PP^N$, and in this case the
divisor $E_{f,P}$ is in~$\Div(B)\otimes\QQ$.
\end{theorem}

Although the restriction to fibral families may seem strong, it is
noted in~\cite{MR3436155} that extending Theorem~\ref{th:variation} to
arbitrary families of RPEs would be enough to extend it to arbitrary
families of endomorphisms.

For a more refined conjecture, we can decompose the canonical
height~$\hhat_f$ associated to a morphism~$f:\PP^N\to\PP^N$ into a sum
of \emph{canonical local heights},
\[
\hhat_f = \sum_{v} \hat \lambda_{f,v},
\]
where the sum is over equivalence classes of absolute values~$v$ on
the field~$K$, and $\hat \lambda_{f,v}$ is a function
$\PP^N(K_v)\to\RR\cup\{\infty\}$ on the points of~$\PP^N$ defined over
the $v$-completion~$K_v$ of~$K$. Roughly speaking, the local canonical
height of a point $P\in\PP^N(K_v)$ is characterized by
\begin{equation}
  \label{eqn:hatlambdafv}
\hat\lambda_{f,v}(P)
= -\log(\text{$v$-adic distance from $P$ to the hyperplane at $\infty$}) + O(1)
\end{equation}
and a transformation formula for the difference
$\hat\lambda_{f,v}\bigl(f(P)\bigr)-(\deg f)\hat\lambda_{f,v}(P)$.
See~\cite{callsilv:htonvariety} for details.

\begin{question}
Given a family of maps $f:\PP^N\to\PP^N$ and a family of
points~$P$, parametrized by~$B$, describe the map
\[
  B(K_v)\longrightarrow\RR\cup\{\infty\},\quad
  t \longmapsto \hat\lambda_{f_t,v}(P_t).
\]
In particular, to what extent is it a Weil function associated to some
divisor in $\Div(B)\otimes\RR$?
\end{question}

The following result gives progress on this problem, contingent on a
technical condition.

\begin{theorem}[Mavraki-Ye~\cite{arxiv1502.04660}]
Let $f_t:\PP^1\to\PP^1$ be a family of morphisms
parametrized by a non-singular curve~$B$,   let $P:B\to\PP^1$ be a
morphism, and let~$D$ be a divisor of degree~$1$ on~$B$.
Suppose further that the pair $(f, P)$ is \emph{quasi-adelic}. Then
\[
  \hhat_{f_t}(P_t)=\hhat_{f}(P)h_{B, D}(t)+O(1).
\]
\end{theorem}

We refer the reader to~\cite{arxiv1502.04660} for the somewhat
technical definition of quasi-adelic, which has to do with the
continuity, or more generally summability, properties of certain
canonical measures on Berkovich~$\PP^1$.

Conjecture~\ref{conj:variation} suggests another natural question. It
turns out that in Tate's formula~\eqref{eqn:tatehftpt} and in
Theorem~\ref{th:variation}, the divisor~$E_{f,P}\in\Div(B)\otimes\RR$
is in fact an element of~$\Div(B)\otimes\QQ$.  The former follows from
the existence of the N{\'e}ron model for elliptic curves, the latter
from the existence of B\"{o}ttcher coordinates for polynomial
endomorphisms. These imply in these cases that the canonical height
$\hhat_{f}(P)$ of~$P$ over the function field of~$B$ is a rational
number.

\begin{question}
\label{question:ffhhatratQ}
Let $K(B)$ be the function field of an algebraic curve, let
$f:\PP^N\to\PP^N$ be a morphism of degree $d\ge2$ defined over the
field~$K(B)$, and let $P\in\PP^N\bigl(K(B)\bigr)$.  Is the canonical
height $\hhat_f(P)$ always in~$\QQ$? If not, give necessary and
sufficient conditions for $\hhat_f(P)\in\QQ$.
\end{question}

Question~\ref{question:ffhhatratQ} is open, even for~$N=1$,
but the following result for local canonical heights is known.

\begin{theorem}[DeMarco-Ghioca~\cite{arxiv1602.05614}]
\label{th:demarcoghioca}
Let $\CC_p$ be a completion of an algebraic closure of the field
of~$p$-adic numbers~$\QQ_p$, and let~$v$ be the usual valuation
on~$\CC_p$, so in particular~$v(\CC_p^*)=\QQ$. Then for every $d\ge2$,
there exists a rational function $f:\PP^1\to\PP^1$ and a point
$P\in\PP^1(\CC_p)$ such that the local canonical height
$\hat{\lambda}_{f,v}(P)$ is irrational.
\end{theorem}

We note that the value of the local canonical height depends on
various choices, but the statement of Theorem~\ref{th:demarcoghioca}
is independent of those choices. Theorem~\ref{th:demarcoghioca} does
not answer Question~\ref{question:ffhhatratQ}, since the global
canonical height is a sum of local canonical heights, but it shows
that the local argument used to prove that
$E_{f,P}\in\Div(B)\otimes\QQ$ for polynomial maps in
Theorem~\ref{th:variation} will not work in general for rational maps.

\subsection{Generalization to Algebraic Varieties}
\label{subsection:algvar:ch}

In this subsection we assume that the reader is familiar with various
standard terms from algebraic geometry.

Let $X/\Qbar$ be a smooth projective variety, let~$f:X\to X$ be a
$\Qbar$-morphism, and let~$D\in\Div_\Qbar(X)\otimes\RR$ be an
$\RR$-linear sum of divisors on~$X$.  One can associate to the divisor
class of~$D$ an equivalence class of height functions
$h_D:X(\Qbar)\to\RR$, well-defined up to bounded functions; see for
example~\cite{hindrysilverman:diophantinegeometry,lang:diophantinegeometry}.

\begin{definition}
Suppose that there is a real number~$\l>1$ with the property that
\begin{equation}
  \label{eqn:fDsimlambdaD}
  f^*D \sim \l D\quad\text{in $\Pic_\Qbar(X)\otimes\RR$.}
\end{equation}
Then the \emph{dynamical canonical height} of~$P$ for the map~$f$ and
divisor (class)~$D$ is the quantity
\begin{equation}
  \label{eqn:hhatdDP}
  \hhat_{f,D}(P) := \lim_{n\to\infty} \frac{1}{\l^n}h\bigl(f^n(P)\bigr).
\end{equation}
\end{definition}

It is known that the limit~\eqref{eqn:hhatdDP}
converges~\cite{callsilv:htonvariety}, and that it satisfies
\[
  \hhat_{f,D}\bigl(f(P)\bigr) = \l \hhat_{f,D}(P)
  \quad\text{and}\quad
  \hhat_{f,D}(P)=h_D(P)+O(1).
\]
Further,
\begin{equation}
  \label{eqn:Damplehhat}
  \text{$D$ is ample} \quad\Longrightarrow\quad \hhat_{f,D}(P)\ge0 \quad\text{and}\quad
  \hhat_{f,D}(P)=0\Longleftrightarrow P\in\PrePer(f).
\end{equation}

Unfortunately, not every map admits an ample divisor satisfying
$f^*D\sim \l D$, but we can find a nice divisor if we relax the
ampleness condition to instead allow~$D$ to be nef, and we relax the
linear equivalence~\eqref{eqn:fDsimlambdaD} to be algebraic
equivalence.

\begin{proposition}
\label{proposition:nefht}
Let $f:X\to X$ be a finite endomorphism of a smooth variety,
let~$\d_f$ be the associated dynamical degree, and suppose
that~$\d_f>1$.
\begin{parts}
\Part{(a)}
\textup{\cite[Remark~29]{MR3483624}} There exists a non-zero nef
divisor class~$\theta_f\in\NS(X)\otimes\RR$ with the property
that\footnote{A divisor is \emph{nef} (\emph{numerically effective})
  if it lies in the closure of the ample cone
  in~$\NS(X)\otimes\RR$. Equivalently, the divisor~$D$ is nef if for
  every ample divisor~$H$ and every~$n\ge1$, the divisor~$nD+H$ is
  ample. Thus non-ample nef divisors are on the edge of being ample.}
\[
  f^*\theta_f \equiv \d_f \theta_f \quad\text{in $\NS(X)\otimes\RR$.}
\]
%% We call~$\theta_f$ a \emph{nef eigenclass for~$f$}.
\Part{(b)}
\textup{\cite[Theorem~5]{MR3483624}} Let~$\theta_f\in\NS(X)\otimes\RR$
be as in~\textup{(a)}. Then the canonical height
limit~\eqref{eqn:hhatdDP} with $D\equiv\theta_f$ and~$\l=\d_f$
converges to give a value~$\hhat_{f,\theta_f}(P)$. For any ample
divisor~$H$ on~$X$, this \emph{nef canonical height} satisfies
\[
  \hhat_{f,\theta_f}\circ f = \d_f\hhat_{f,\theta_f}
  \quad\text{and}\quad
  \hhat_{f,\theta_f} = h_{\theta_f} + O(\sqrt{h_H}).
\]
\end{parts}
\end{proposition}

%% \begin{remark}
%% The convergence of~$\hhat_{f,\theta_f}(P)$ in
%% Proposition~\ref{proposition:nefht} is not a consequence of our
%% earlier statements, because the algebraic equivalence
%% $f^*\theta_f\equiv\d_f\theta_f$ does not imply that the
%% difference~$h_\theta\circ f-\d_f h_\theta$ is bounded. 
%% \end{remark}

\begin{question}
\label{question:nefht}
Let~$X$,~$f$,~$\theta_f$, and~$\hhat_{f,\theta_f}$ be as in
Proposition~\ref{proposition:nefht}.
\begin{parts}
\Part{(a)}
Is it true that $\hhat_{f,\theta_f}(P)\ge0$ for all~$P\in X(\Qbar)$? 
\Part{(b)}
Describe the set of points
\begin{equation}
  \label{eqn:hhatfthetaPle0}
  \bigl\{ P\in X(\Qbar) : \hhat_{f,\theta_f}(P)\le 0 \bigr\},
\end{equation}
beyond the obvious fact that it contains~$\PrePer(f)$.
\end{parts}
\end{question}

\begin{remark}
Here are three cases where the answer to Question~\ref{question:nefht}
is known.%
\begin{parts}
\Part{(1)}
If~$\theta_f$ is ample, then it is easy to see that
$\hhat_{f,\theta_f}\ge0$ and that the set~\eqref{eqn:hhatfthetaPle0}
equals~$\PrePer(f)$.
\Part{(2)}
Let~$X$ be an abelian variety and let~$f$ be an isogeny. Then
$\hhat_{f,\theta_f}\ge0$, and there is a proper abelian
subvariety~$Y\subsetneq X$ such that the
set~\eqref{eqn:hhatfthetaPle0} in~(b) is equal to the set of
translates $Y+X_\tors$; see~\cite{MR3456169}.
\Part{(3)}
Let $X$ be a K3 surface with $\rank\NS(X)=2$, and let~$f:X\to X$ be an
automorphism of infinite order obtained by composing non-commuting
involutions.  Then $\hhat_{f,\theta_f}\ge0$ and the
set~\eqref{eqn:hhatfthetaPle0} equals~$\PrePer(f)$;
see~\cite{silverman:K3heights}.
\end{parts}
\end{remark}

Even less is known for rational maps that are not morphisms.
The following was suggested by Mattias Jonsson.

\begin{question}
\label{question:highdyndeg}
Let $X$ be a smooth projective variety of dimension~$N$, let
$f:X\dashrightarrow X$ be a dominant rational map, and suppose that
the higher order dynamical degrees~$\d_f^{(1)},\ldots,\d_f^{(N)}$
defined earlier~\eqref{eqn:highorderdyndeg} for a strictly log concave
sequence, i.e., they satisfy $(\d_f^{(i)})^2 >
\d_f^{(i-1)}\d_f^{(i+1)}$ for all $1\le i\le N-1$.
Also let~$D$ be an ample divisor. 
\begin{parts}
\Part{(a)}
Is it true that $\deg_D(f^n):=(f^n)^*D\cdot
D^{N-1}\asymp(\d_f^{(1)})^n$?
\Part{(b)}
Suppose that~$X$,~$f$, and~$D$ are defined over~$\Qbar$, and
let~$P\in{X}(\Qbar)_f$. Does the limit
\[
  \hhat_{f,D}(P) := \lim_{n\to\infty} (\d_f^{(1)})^{-n} h_D\bigl(f^n(P)\bigr)
\]
converge? If so, is it independent of~$A$ up to a multiplicative
constant.
\Part{(c)}
Assuming that the limit in~\text{(b)} exists,
describe the set of~$P\in X(\Qbar)$ satisfying $\hhat_{f,D}(P)=0$, beyond the
fact that it contains~$\PrePer(f)$. 
\end{parts}
\end{question}

For $\dim(X)=2$, the answer to Question~\ref{question:highdyndeg}(a) is
yes~\cite{arxiv0608267}, and for birational~$f$ there are partial
results for~(b) and~(c) in~\cite{arxiv1505.03559}.

%%%%%%%%%%%%%%%%%%%%%%%%%%%%%%%%%%%%%%%%%%%%%%%%%%%%%%%%%%%%%%%%%%%%%%
\section{$p$-adic and Nonarchimedean Dynamics\WhoDone{Rob}}
\label{section:wdpd}
\WhoWrite{Rob}
%%%%%%%%%%%%%%%%%%%%%%%%%%%%%%%%%%%%%%%%%%%%%%%%%%%%%%%%%%%%%%%%%%%%%%

As noted in Section~\ref{section:introduction}, the study of dynamics
over global fields leads naturally to questions about dynamics over
local fields. In addition, the rich theory of complex dynamics on the
one hand provides a model for a theory of $p$-adic dynamics, and on
the other hand raises the question of how much of the archimedean
theory generalizes to nonarchimedean fields. Moreover, besides
$p$-adic fields, there are questions in complex dynamics whose
solution uses dynamics over Laurent series fields such as $\CC((t))$,
which are nonarchimedean, but not $p$-adic. Thus, we consider the
following general setting.

Let~$\CC_v$ be an algebraically closed field that is complete with
respect to an absolute value~\text{$|\cdot|_v$} satisfying the
nonarchimedean triangle inequality:
\[
  |a+b|_v \le \max\bigl\{ |a|_v,|b|_v\bigr\}
  \quad\text{for all $a,b\in\CC_v$.}
\]
For example, for any prime number~$p$, the field~$\CC_p$
is the completion of an algebraic closure of the~$p$-adic
field~$\QQ_p$. We set the additional notation
\[
  \Ocal_v := \bigl\{a\in\CC_v : |a|_v\le 1 \bigr\}\quad\text{and}\quad
  \Mcal_v := \bigl\{a\in\CC_v : |a|_v< 1 \bigr\}
\]
for, respectively, the ring of integers and the maximal ideal
for~$\CC_v$.  

The projective line~$\PP^1(\CC_v)$ over~$\CC_v$ is totally
disconnected and not locally compact, but it can be embedded in the
\emph{Berkovich projective line}~$\PBerk$, which is a compact,
Hausdorff, and path-connected topological space that
contains~$\PP^1(\CC_v)$ as a dense subset.  We briefly describe some
of the salient features of~$\PBerk$ and its endomorphisms, and refer
the reader to \cite[Chapters~1--2]{bakerrumelyberkbook},
\cite[Chapters~6--7]{benedettobook}, or \cite[Section~5.10]{MR2316407}
for further details.

The points of~$\PBerk$ that lie in~$\PP^1(\CC_v)$ are said to be of
Type~I.  The points of Type~II and~III in~$\PBerk$ correspond to
closed disks in~$\CC_v$. Type~II points correspond to \emph{rational
  disks}, which are disks whose radius is in~$|\CC_v^{\times}|_v$,
while Type~III points correspond to \emph{irrational disks}, which are
disks whose radius is not in~$|\CC_v^{\times}|_v$. Finally, there are
Type~IV points that will play only a minor role in this
survey.\footnote{Type~IV points correspond to decreasing
  chains~$D_1\supset{D}_2\supset\cdots$ of disks in~$\CC_v$ with empty
  intersection. They are needed to ensure that~$\PBerk$ is compact.}
The Type~II point~$\zeta(0,1)$ corresponding to the closed unit
disk~$\bar{D}(0,1)\subset\CC_v$ is called the \emph{Gauss point}
of~$\PBerk$.

A rational function~$f(z)\in\CC_v(z)$ defines a
map~$f:\PP^1(\CC_v)\to\PP^1(\CC_v)$, and this map extends continuously
to a function
\[
  f:\PBerk\to\PBerk .
\]
Each point~$\zeta\in\PBerk$ maps to its image~$f(\zeta)$ with a
\emph{local degree}, or \emph{multiplicity}, which is an integer
between~$1$ and~$\deg(f)$. If the local degree of~$f$ at~$\zeta$ is
divisible by the residue characteristic of~$\CC_v$, we say
that~$\zeta$ is \emph{wild}.

For a Type~I point~$\zeta$, the image~$f(\zeta)$ is defined using the
map~$f:\PP^1(\CC_v)\to\PP^1(\CC_v)$, and the local degree at~$\zeta$
is the usual multiplicity (or ramification index) of~$f$ at~$\zeta$.
To describe how~$f$ acts on points of Type~II and~III, we note that
if~$D_1\subset\CC_v$ is a disk, then the image~$f(D_1)$ is either a
disk~$D_2\subset\CC_v$, the complement~$\PP^1(\CC_v)\smallsetminus
D_2$ of a disk, or all of~$\PP^1(\CC_v)$.  In the first two cases, the
action of~$f$ on~$\PBerk$ maps the point~$\zeta_1\in\PBerk$
corresponding to~$D_1$ to the point~$\zeta_2\in\PBerk$ corresponding
to~$D_2$.  In particular, if~$D_1$ and~$D_2$ are closed disks
with~$D_2=f(D_1)$, then the local degree of~$f$ at~$\zeta_1$
is~$\#f^{-1}(\zeta)$ for any point~$\zeta\in D_2$.

Let~$f\in\CC_v(z)$, and let~$a\in\PP^1(\CC_v)$ be an $f$-periodic
point of of exact period~$n$. For ease of exposition, we assume
that~$\infty\notin\Orbit_f(a)$.\footnote{The multiplier~$\lambda_f(a)$
  is invariant under change of coordinates, so
  if~$\infty\in\Orbit_f(a)$, we can compute~$\lambda_f(a)$ by first
  changing coordinates.}  Then the \emph{multiplier of~$f$ at~$a$} is the quantity
\[
  \lambda_f(a) = \big( f^n \big)'(a) = \prod_{i=0}^{n-1} f'\big( f^i(a) \big) .
\]
We say that~$a$ is
\begin{itemize}
\item \emph{attracting} if~$\bigl|\lambda_f(a)\bigr|_v < 1$,
\item \emph{indifferent}, or \emph{neutral}, if~$\bigl|\lambda_f(a)\bigr|_v = 1$,
\item \emph{repelling} if~$\bigl|\lambda_f(a)\bigr|_v > 1$.
\end{itemize}

Periodic points of Types II, III, and IV
have multipliers that are positive integers; more precisely, the multiplier of such
a periodic point~$\zeta$ is defined to be the local degree of~$f^n$ at~$\zeta$,
where~$n$ is the period of~$\zeta$ under~$f$. Periodic points of Type~III
and~IV necessarily have multiplier~1, because any rational function fixing a
point~$\zeta$ of Type~III or~IV must act as a homeomorphism on a neighborhood of~$\zeta$.

The multiplier of a Type~II point is more interesting and is computed
as follows. Suppose that~$\zeta=\zeta(a,r)$ is the point corresponding
to the rational closed disk $\bar{D}(a,r)$, where~$a\in\CC_v$
and~$r\in |\CC_v^{\times}|_v$.  Replacing~$f$ by~$\phi^{-1} \circ f
\circ \phi$ for an appropriate $\phi(z) = cz+a$, we may
assume~$\zeta=\zeta(0,1)$ is the Gauss point.  Write~$f^n=g/h$ as a
quotient of relatively prime polynomials~$g,h\in\Ocal_v[z]$ such that
at least one coefficient~$u$ of~$g$ or~$h$ satisfies~$|u|_v=1$.  Then
we can reduce the coefficients of~$g$ and~$h$ modulo~$\Mcal_v$ to
obtain a rational function~$\bar{g}/\bar{h}\in (\Ocal_v/\Mcal_v)(z)$
with coefficients in the residue field.  The assumption that~$f^n$
fixes the Gauss point implies that~$\bar{g}/\bar{h}$ is nonconstant,
and then the multiplier of~$f$ at~$\zeta$ is the degree of the
rational function~$\bar{g}/\bar{h}$.  We say that a periodic
point~$\zeta$ of Type~II, III, or~IV is \emph{indifferent} if its
multiplier is~$1$, and \emph{repelling} if its multiplier is at
least~$2$.  (Points of Type~II, III, and~IV are never attracting.)

The (\emph{Berkovich}) \emph{Fatou set} of~$f\in\CC_v(z)$ is the set~$\Fcal_f$
of those points~$\zeta\in\PBerk$ which have a neighborhood~$U$ such that
\[
  \PBerk \smallsetminus \bigg( \bigcup_{n\geq 0} f^n(U) \bigg)
  \quad \text{is infinite}.
\]
The (\emph{Berkovich}) \emph{Julia set} of~$f$ is the complement
$\Jcal_f = \PBerk\smallsetminus \Fcal_f$.
For alternative equivalent definitions,
see for example~\cite[Section~10.5]{bakerrumelyberkbook}.
The Fatou and Julia sets are both invariant under~$f$, i.e., 
\[
  f(\Fcal_f) = f^{-1}(\Fcal_f) = \Fcal_f
  \quad\text{and}\quad
  f(\Jcal_f) = f^{-1}(\Jcal_f) = \Jcal_f .
\]
For any~$n\geq 1$, the Fatou set of~$f^n$ coincides with the Fatou set of~$f$,
and similarly for the Julia set.

All repelling periodic points lie in the Julia set, while the Fatou
set contains all attracting periodic points, all indifferent periodic
points of Types~I, III, and~IV, and most of the Type~II indifferent
points. However, if the residue field of~$\CC_v$ is not algebraic over
a finite field, e.g., if~$\CC_v$ has residue characteristic zero, then
some Type~II indifferent periodic points may lie in the Julia set.  By
way of contrast, in classical complex dynamics, 
most indifferent periodic points lie in the Julia set.  The Berkovich Fatou set is
always nonempty, because~$f$ necessarily has a nonrepelling fixed
point of Type~I; see \cite[Corollary~1.3]{benedetto:hyperbolic}.  The
Julia set is also nonempty, because~$f$ necessarily has a repelling
fixed point of Type~I or~II; see \cite[Th\'{e}or\`{e}me~B]{MR2018827}.
However, the Julia set may contain no Type~I points, i.e.,
it is possible to have~$\Jcal_f\cap \PP^1(\CC_v)=\emptyset$.
This happens, for example, if~$f$ has good reduction; cf.\ Section~\ref{section:goodred}.

It is true that the repelling periodic points are dense in the Julia
set~\cite[Section~12.3]{benedettobook}, but if we restrict our attention
to~$\PP^1(\CC_v)$, i.e., to the Type~I points, then things are 
more complicated. We pose a conjecture.

\begin{conjecture}
\label{conj:repdense}
Let~$f\in\CC_v(z)$ be a rational function of degree at least~$2$, and
let~$\Jcal_f\subset\PBerk$ be the Julia set
of~$f$. Then~$\Jcal_f\cap\PP^1(\CC_v)$ is equal to the closure
in~$\PP^1(\CC_v)$ of the set of Type~I repelling periodic points
of~$f$.
\end{conjecture}

Hsia \cite{hsia:periodicpoint} proved that~$\Jcal_f\cap\PP^1(\CC_v)$
is contained in the closure of \emph{all} Type~I periodic points, but
in contrast to the situation over~$\CC$, a rational map over~$\CC_v$
may have infinitely many non-repeling periodic points. On the other
hand, B\'{e}zivin \cite{MR1864626} proved that if~$f$ has at least one
Type~I repelling periodic point, then Conjecture~\ref{conj:repdense}
is true.

Let~$U$ be a connected component of~$\Fcal_f$. Then the image~$f(U)$
is also a connected component of~$\Fcal_f$. Further, there is an integer
$m \in \{1,\ldots, d\}$ such that every point in~$f(U)$ has exactly~$m$ preimages
in~$U$, counted with multiplicity; we say that~$f$ is~$m$-to-$1$ on~$U$.
A periodic component~$U$ of period~$n\geq 1$ is said to be
\begin{itemize}
\item \emph{attracting} if~$f^n$ is~$m$-to-$1$ on~$U$ with~$m\geq 2$;
\item \emph{indifferent} if~$f^n$ is~$1$-to-$1$ on~$U$.
\end{itemize}
Rivera-Letelier proved the following classification theorem
over~$\CC_p$, and his proof can be extended
to arbitrary~$\CC_v$ \cite[Chapter~9]{benedettobook}.  To state it, we
define a \emph{domain of Cantor type}\footnote{The name comes from the
  fact that the boundary of such a domain is homeomorphic to the
  Cantor set.}  to be an increasing union~$\bigcup_{j\geq 1} U_j$,
where each~$U_j$ is the complement in~$\PBerk$ of~$N_j$ disjoint
closed Berkovich disks, with~$N_j < N_{j+1}$
and~$\overline{U}_j\subset U_{j+1}$.

\begin{theorem}
[Rivera-Letelier, 2000, \cite{MR2040006}]
Let~$f\in\CC_v(z)$ be a rational function of degree at least~$2$,
and let~$U\subseteq\Fcal_f$ be a periodic Fatou component of period~\text{$n\geq1$}.
\begin{parts}
\Part{(a)}
If~$U$ is an attracting component, then there is
a Type~I attracting periodic point~$x\in\Fcal_f$
such that all points of~$U$ are attracted to~$x$ under iteration of~$f^n$.
Moreover,~$U$ is either a rational open disk or a domain of Cantor type.
\Part{(b)}
If~$U$ is an indifferent component, then~$U$ is a region of the form
\[
  U=\PBerk\smallsetminus \big( V_1 \cup \cdots \cup V_N \big),
\]
where~$V_1,\ldots,V_N$ are disjoint rational closed disks.
The~$N$ boundary points of~$U$ are all periodic under~$f^n$.
\end{parts}
\end{theorem}

In the same paper, Rivera-Letelier proved that a map~$f\in\CC_v(z)$ of
degree~$d\geq 2$ can have at most~$d-1$ cycles of attracting periodic
components that are of Cantor type.  However, if there are wild
points~$\zeta\in\PBerk$, i.e., points at which the local degree of~$f$
is divisible by the residue characteristic, then there may be
infinitely many attracting components that are disks. For example,
the Gauss point is fixed and wild for
the map~$f(z)=z^p\in\CC_p(z)$, and
every periodic residue class~$D(a,1)$ is an attracting component.  By
contrast, if there are no wild points in~$\PBerk$, for example
if~$d<p$ or~$p=0$, then there are at most~$2d-2$ cycles of attracting
components; see \cite[Corollary~1.6]{MR3265554}.

There can also be infinitely many indifferent periodic disks. For
example, if~$f(z)=z^d$, and if the residue characteristic~$p$ does not
divide~$d$, then every periodic residue class~$D(a,1)$ besides the
ones at~$0$ and~$\infty$ is indifferent. Rivera-Letelier showed
that~$f$ can have many cycles of non-disk indifferent
components~\cite[Proposition~6.4]{MR2040006}, but his construction
increases the degree of~$f$ as it increases the number of such
cycles, suggesting the following question.

\begin{question}
Let~$f\in\CC_v(z)$ with~$\deg f = d\geq 2$. Can~$f$ have infinitely many periodic cycles
of indifferent components that are not disks? If not, is there a bound in terms of~$d$
for the number of such cycles? Do the answers depend on the field~$\CC_v$? 
\end{question}

A component of~$\Fcal_f$ that is not preperiodic is said to be a
\emph{wandering domain}.  A famous theorem of
Sullivan~\cite{sullivan:nowanderingdomains} in complex dynamics says
that the Fatou set of a rational function~$f\in\CC(z)$ has no
wandering domains in~$\PP^1(\CC)$. The nonarchimedean situation is
more complicated.  The following result first appeared for~$p$-adic
fields, but its proof generalizes to function fields of positive
characteristic under stronger hypotheses.  To state it, we say that a 
point~$b$ is \emph{recurrent} under~$f$ if~$b$ is itself an
accumulation point of its forward orbit~$\Orbit_f(b)$.

\begin{theorem}[Benedetto, 2000, \cite{MR1781331}]
\label{thm:RobNWD}
Let~$K$ be a locally compact subfield of~$\CC_v$, and suppose that the
algebraic closure~$\bar{K}$ of~$K$ is dense in~$\CC_v$.  Let~$p\geq 2$
be the residue characteristic of~$\CC_v$, and let~$f\in K(z)$ be a
rational function of degree at least~$2$.  Suppose that one of the
following two conditions holds\textup:
\begin{parts}
\Part{(a)}
$\operatorname{char} K=0$, and there are
no wild recurrent critical points of Type~I in~$\Jcal_f$.
\Part{(b)}  
$\operatorname{char} K=p$, and there are
no wild critical points of Type~I in~$\Jcal_f$.
\end{parts}
Then~$f$ has no wandering domains.
\end{theorem}

The residue field of a locally compact nonarchimedean field is
necessarily finite. Thus Theorem~\ref{thm:RobNWD}(a) applies
to~$p$-adic fields, and Theorem~\ref{thm:RobNWD}(b) applies to
function fields over finite fields. However, if~$f\in\CC_v(z)$ is
\emph{not} required to be defined over a locally compact subfield
of~$\CC_v$, then~$f$ may have wandering domains, as in the following
explicit result.  The construction relies heavily on both wild
behavior, i.e., maps that are locally~$m$-to-$1$ with~$p\mid m$, and the
fact that~$\CC_v$ is not locally compact.

\begin{theorem}[Benedetto, 2002 \&\ 2006, \cite{MR1941304,arxiv0312029}]
\label{thm:RobWD}
Suppose that~$\CC_v$ has residue characteristic~$p\geq 2$.  Then for
any~$b\in\CC_v$ with~$|b|_v>1$ and any~$\epsilon>0$, there
exists~$a\in\CC_v$ with~$|a-b|_v<\epsilon$ such that the Fatou set of
the polynomial
\[
  f(z) = (1-a) z^{p+1} + az^p
\]
has a wandering domain.
\end{theorem}

The contrasting conclusions of Theorems~\ref{thm:RobNWD}
and~\ref{thm:RobWD} suggest the following conjecture.

\begin{conjecture}
[\cite{MR1781331}]
Let~$K$ be a locally compact subfield of~$\CC_v$, and suppose that the
algebraic closure~$\bar{K}$ of~$K$ is dense in~$\CC_v$.
Let~$f\in{K}(z)$ be a rational function of degree at least~$2$.
Then~$f$ has no wandering domains.
\end{conjecture}

In contrast to the locally compact case,
if the residue field of~$\CC_v$
is \emph{not} an algebraic extension of a finite field, then it turns
out that there are wandering domains attached to \emph{every}
periodic point of Type~II. For example, the map
$f(z)=z^d$ fixes the Gauss point,
and every open disk of the form~$U=D(b,1)$ with~$|b|=1$ satisfies
$\phi^n(U) = D(b^{d^n},1)$. If~$\CC_v$ is of residue characteristic~$0$,
then some such disks~$U$ contain no roots of unity, e.g., for~$b=2$,
and in fact,~$U$ is a wandering domain of~$\Fcal_f$. Such wandering
domains are said to lie in the \emph{attracting basin} of the associated
(Julia) periodic point of Type~II.
However, it appears that no other wandering domains are possible
in residue characteristic zero.

\begin{theorem}[(a) Benedetto, 2005, {\cite[Theorem~B]{arxiv0312034}}; (b) Trucco,
  2014, \cite{MR3295719}]
\label{thm:RobTrucco}
Let~$f\in\CC_v(z)$ be a rational function of degree~$d\geq 2$,
and let~$p$ be the residue characteristic of~$\CC_v$.
Suppose that at least one of the following conditions holds\textup:%
\begin{parts}
\Part{(a)}
The function $f$ is defined over a discretely valued subfield~$K$ of~$\CC_v$
whose algebraic closure~$\bar{K}$ is dense in~$\CC_v$,
and either~$p=0$ or~$p>d$.
\Part{(b)}
The function $f$ is a polynomial, i.e., $f\in\CC_v[z]$, and~$p=0$.
\end{parts}
Then every wandering domain of~$f$ lies in the attracting basin
of a  Type~II Julia periodic point.
\end{theorem}

The condition in~Theorem~\ref{thm:RobTrucco}(a) that~$\bar{K}$ is dense in~$\CC_v$ cannot be
removed, because if~$|\bar{K}^{\times}|_v \subsetneq |\CC_v^\times|_v$,
then any Latt\`{e}s map arising from an elliptic curve having
multiplicative reduction will have wandering domains defined
over~$\CC_v$, but which correspond to Type~III points over the
completion of~$\bar{K}$.

The conditions on the residue characteristic in
Theorem~\ref{thm:RobTrucco} are used mainly to eliminate the
possibility of wild behavior. Thus we are led to the following
conjecture, which applies to maps defined over the~$p$-adic
field~$\CC_p$ and maps defined over nonarchimedean function fields of
arbitrary characteristic.

\begin{conjecture}
\label{conj:char0NWD}
Let~$K\subset\CC_v$ be a discretely valued subfield whose algebraic
closure~$\bar{K}$ is dense in~$\CC_v$.  Let~$p\geq 0$ be the residue
characteristic of~$\CC_v$, and let~$f\in\CC_v(z)$ be a rational
function of degree at least~$2$.  Suppose that there are no
points~$\zeta\in\PBerk$ at which~$f$ has local degree divisible
by~$p$.  Then every wandering domain of~$f$ lies in the attracting
basin of a Type~II Julia periodic point.
\end{conjecture}

\begin{remark}
Another important tool used in studying nonarchimedean dynamics is
measure theory.  Given~$f\in\CC_v(z)$, there is a naturally associated
Borel probability measure~$\mu_f$ on~$\PBerk$ called the
\emph{canonical measure} or \emph{invariant measure} or
\emph{equilibrium measure} for~$f$. Constructed independently by Baker
and Rumely~\cite{arxiv0407426}, Chambert-Loir~\cite{arxiv0304023}, and
Favre and Rivera-Letelier~\cite{MR2092012,arxiv0407471}, the
measure~$\mu_f$ has many useful properties. In particular, its support
is precisely the Julia set~$\Jcal_f$, and it is invariant under~$f$,
i.e.,~$\mu_f(f^{-1}(E))=\mu_f(E)$ for every Borel subset~$E$
of~$\PBerk$.  The existence of~$\mu_f$ permits a meaningful study of
ergodic theory, Lyapunov exponents, and more, in the context of
nonarchimedean dynamics.  In addition, if~$f$ is defined over a global
field~$K$, then the measures~$\mu_{f,v}$ attached to the different
absolute values~$v$ on~$K$ are closely related to the local canonical
height functions~\eqref{eqn:hatlambdafv} discussed in
Section~\ref{section:ch}.  A deeper discussion of the equilibrium
measure and open questions surrounding it would unfortunately take us
too far afield in this survey; we refer the interested reader
to~\cite[Section~10.1]{bakerrumelyberkbook}
or~\cite[Chapter~13]{benedettobook} for more information.
\end{remark}

\begin{remark}
As mentioned in the introduction, there is a thriving theory of
ergodic theory and entropy over $p$-adic fields dating back to at
least a 1975 paper of Oselies and Zieschang~\cite{MR0369301}.  For
recent work on $p$-adic ergodic theory for rational maps, see the
paper~\cite{MR2578470} of Favre and Rivera-Letelier, and for a
discussion of $p$-adic entropy, see the paper~\cite{lindward:solenoids}
of Lind and Ward.
\end{remark}

%%%%%%%%%%%%%%%%%%%%%%%%%%%%%%%%%%%%%%%%%%%%%%%%%%%%%%%%%%%%%%%%%%%%%%
\section{Dynamics over Finite Fields\WhoDone{Joe}}
\label{section:dff}
\WhoWrite{Joe}
%%%%%%%%%%%%%%%%%%%%%%%%%%%%%%%%%%%%%%%%%%%%%%%%%%%%%%%%%%%%%%%%%%%%%%

Projective space~$\PP^N(\FF_q)$ over a finite field~$\FF_q$ contains
only finitely many points, precisely
$\#\PP^N(\FF_q)=(q^{N+1}-1)/(q-1)$, so a morphism $f:\PP^N\to\PP^N$
defined over~$\FF_q$ is a self-map of the finite set~$\PP^N(\FF_q)$.
Thus every point is preperiodic, but it is interesting to study the
subset of points that are strictly periodic.

\begin{definition}
Let~$S$ be a (finite) set. We write
\begin{align*}
  \End_\Set(S)&=\text{the collection of set maps $S\to S$,} \\*
  \Aut_\Set(S)&=\text{the collection of set bijections $S\to S$,} \\
  \Per(f,S)&=\text{the set of strictly periodic points for~$f\in\End_\Set(S)$} \\*
  &= \{P\in S : f^n(P)=P~\text{for some $n\ge1$}\}.
\end{align*}
We note that a map $f\in\End_\Set(S)$ gives~$S$ the structure of a
directed graph~$\G_f$ whose vertices are the element of~$S$ and such
that there is an arrow from~$P$ to~$f(P)$ for each~$P\in S$.
\end{definition}

If we let~$f$ range over $\End_\Set(S)$ or $\Aut_\Set(S)$ and let
$\#S\to\infty$, then there are known statistics for $\#\Per(f,S)$ and
for various natural quantities associated to the graph~$\G_f$ such
as its number of components, its largest component, its largest cycle,
and many others. See for example~\cite{MR1083961,MR0119227,MR0062973}.
In particular, Flajolet and Odlyzko~\cite{MR1083961} analyze a large
number of properties for average random maps, and one may ask whether
similar statistics hold if one restricts to maps defined by
polynomials.

\begin{question}[Vague Motivating Question]
To what extent do polynomial maps behave like random set maps?
\end{question}

In this section we briefly discuss some potential answers to this
question, where we note that there are many sorts of averages that one
might consider. For example, we can average over maps defined over a
given field~$\FF_q$, or we can fix a map defined over~$\FF_q$ and ask
what happens to the points in~$\PP^N(\FF_{q^n})$ as $n\to\infty$, or we can fix a map
defined over~$\QQ$ and average over its mod~$p$ reductions.

It is important to note that the collection of degree~$d$ morphisms
$f:\PP^N(\FF_q)\to\PP^N(\FF_q)$ is generally only a small part of the
collection of set maps $\End_\Set\bigl(\PP^N(\FF_q)\bigr)$. For
example, the inverse image of a point $P\in\PP^N(\FF_q)$ satisfies $\#
f^{-1}(P)\le\deg(f)$, and further, for a given~$1\le k\le\deg(f)$, the
``probability'' that $\#f^{-1}(P)=k$ may depend on~$k$ in subtle ways;
see~\cite{MR3681355}.

We have already discussed in Section~\ref{section:goodred} the
situation where we start with a morphism $f:\PP^N\to\PP^N$ defined
over~$\QQ$ and a point~$P\in\PP^N(\QQ)$, and we reduce modulo various
primes~$p$ and study the orbit of the reduced point~$\tilde P$ for the
reduced map~$\tilde{f}$. In particular,
Propostion~\ref{proposition:goodredeval} says that iteration and
reduction modulo~$p$ interact nicely with one another for all but
finitely many primes~$p$. If~$P$ is $f$-wandering in~$\PP^N(\QQ)$, it
is interesting to ask how~$\Orbit_{\tilde f}(\tilde P\bmod p)$ varies
as we vary the prime~$p$. For example, is it usually large, as
predicted by the following conjecture,\footnote{This conjecture should
  probably be viewed as a provocation. A safer formulation would be to
  phrase it as the problem of classfiying all maps for which the
  conjecture fails to be true!}  where we refer the reader to
Definition~\ref{definition:density} for the definition of the
(natural) density of a set of primes.

\begin{conjecture}
\label{conjecture:sizeorbitmodp}
\begin{parts}
\Part{(a)}
Let $f:\PP^N\to\PP^N$ be a morphism of degree $d\ge2$ defined over
$\QQ$, and let $P\in\PP^N(\QQ)$ be a point whose $f$-orbit
$\Orbit_f(P)$ is Zariski dense in~$\PP^N$. Then for every~$\e>0$, the
set of primes
\[
  \bigl\{ p : \#\Orbit_{\tilde f}(\tilde P\bmod p) \ge  p^{(1/2-\e)N} \bigr\}
\]
is a set of density~$1$.
\Part{(b)}
Let $f:\AA^N\to\AA^N$ be a polynomial automorphism defined over~$\QQ$,
i.e., the map~$f$ is defined by polynomials, and it has an inverse
defined by polynomials.  Let $P\in\AA^N(\QQ)$ be a point whose
$f$-orbit $\Orbit_f(P)$ is Zariski dense in~$\AA^N$.  Then for
every~$\e>0$, the set of primes
\[
  \bigl\{ p : \#\Orbit_{\tilde f}(\tilde P\bmod p) \ge  p^{(1-\e)N} \bigr\}
\]
is a set of density~$1$.
\end{parts}
\end{conjecture}

Only much weaker estimates are known. Here is a typical result; see
also~\cite{arXiv:0707.1505,MR3767353} for generalizations involving
parameterized families of collections of maps.

\begin{theorem}
\cite[Akbary--Ghioca (2009)]{MR2549537}
Let $\e>0$.   
\begin{parts}
\Part{(a)}
Let $f:\PP^N\to\PP^N$ be an endomorphism defined over~$\QQ$.  Then the
set of primes
\[
  \left\{ p : \#\Orbit_{\tilde f}(\tilde P\bmod p)
  \ge \frac{\log(p)}{\deg(f)+\e}\right\}
\]
is a set of density~$1$. 
\Part{(b)}
Let $f:\AA^N\to\AA^N$ be a polynomial automorphism defined
over~$\QQ$. Then the set of primes
\[
  \left\{ p : \#\Orbit_{\tilde f}(\tilde P\bmod p)
  \ge p^{1/2-\e}\right\}
\]
is a set of density~$1$.
\end{parts}
\end{theorem}

In the context of Conjecture~\ref{conjecture:sizeorbitmodp}(a), one
might ask instead for a density~$1$ conclusion for primes whose orbits
are smaller than~$p^{(1/2+\e)N}$.  As explained in~\cite{MR3731313},
it seems unlikely that such a result holds in general, so we pose the
problem of characterizing those maps for which there is such a result.

\begin{question}
Let $f:\PP^N\to\PP^N$ be a morphism of degree $d\ge2$ defined over
$\QQ$, and let $P\in\PP^N(\QQ)$ be a point whose $f$-orbit
$\Orbit_f(P)$ is Zariski dense in~$\PP^N$.
Under what conditions on~$f$ is it true that for every~$\e>0$, the
set of primes
\[
  \bigl\{ p : \#\Orbit_{\tilde f}(\tilde P\bmod p) \le  p^{(1/2+\e)N} \bigr\}
\]
is a set of density~$1$?
\end{question}

There are a handful of published results giving upper bounds for the size of
the orbit \text{$\Orbit_{\tilde{f}}(\tilde P\bmod p)$},
including~\cite{MR3731313,MR3544625,MR3437764}.  To illustrate, we
quote an easily stated corollary from one of these papers.

\begin{theorem}
[Heath-Brown, 2017, \cite{MR3731313}] Let $f(x)=ax^2+b$ with
$a,b\in\ZZ$ positive integers. There is a constant $C(a,b)$ such that
for every $\a\in\ZZ$ and every odd prime~$p$,
\[
  \#\Orbit_{\tilde f}(\tilde\a\bmod p) \le C(a,b)\frac{p}{\log\log p}.
\]
\end{theorem}

We next consider the graph~$\G_f$ induced by the action of~$f$ on the
points of~$\PP^N$.  For any graph~$\G$, we write
\[
  \Comp(\G) = \text{\# of components in the graph $\G$}.
\]
Kruskal~\cite{MR0062973} proved that the average value of~$\Comp(\G)$
for random maps of a set~$S$ with~$n$ elements is given by
\[
  \Compavg(n) := \frac{1}{n^n} \sum_{f\in\End_\Set(S)} \Comp(\G_f)
  = \frac12\log n + C + o(1),
\]
with an explicit value for~$C$. This suggests a natural question.

\begin{question}
\label{question:avgcomp}
Fix an integer $N\ge1$ and a prime power $q$. Is it true that
\[
  \frac{1}{\#\End_d^N(\FF_q)} \sum_{f\in\End_d^N(\FF_q)} \Comp(\G_f) = \Compavg\bigl(\#\PP^N(\FF_q)\bigr) + O(1)
  \quad\text{as $d\to\infty$?}
\]
\end{question}

One could, in Question~\ref{question:avgcomp}, ask for a stronger
result with $o(1)$, or for a weaker result with~$\sim$ or~$\asymp$ in
place of~$=$. And rather than taking endomorphisms of~$\PP^N$, one
could instead average over endomorphisms of~$\AA^N$, i.e., over
polynomial self-maps of~$\FF_q^N$.

Flynn and Garton~\cite{MR3190008} prove that if~$N=1$
and~$d\ge\sqrt{q}$ , then the average in
Question~\ref{question:avgcomp} is greater than $\frac12q-4$. See
also~\cite{MR3681355} for an affirmative answer to
Question~\ref{question:avgcomp} for~$N=1$ assuming an
independence hypothesis.  For material on the connectivity of the
graphs~$\G_{x^2+c}$ as~$c$ ranges over~$\FF_q$, see~\cite{MR3425238};
and for results on the number of non-isomorphic graphs~$\G_f$ as~$f$
ranges over polynomials in~$\FF_q[x]$ of fixed degree,
see~\cite{doi:10.1080/10586458.2017.1391725}.

Rather than averaging over all maps in$\End_d^N(\FF_q)$, we might
instead fix a map defined over a number field and look at the
components of the graphs that we obtain by reducing modulo different
primes. This leads to the following question.

\begin{question}
Let $f:\PP^N\to\PP^N$ be a morphism of degree $d\ge2$ defined
over~$\QQ$. Is it true that for every $\e>0$, the set of primes
\[
  \left\{ p : \Comp(\G_{\tilde f\bmod p}) \le \Compavg(p^{(1-\e)N})\right\}
\]
is a set of density~$0$?  Similarly with~$\AA^N$ in place of~$\PP^N$? For
which maps~$f$ can one obtain stronger results, e.g.,
$\Comp(\G_{\tilde{f}\bmod p})\asymp\Compavg(p^N)$ as $p\to\infty$?
\end{question}

For our next problem, we start with an (irreducible) polynomial
$f(x)\in\FF_q[x]$ and ask how the iterates~$f^n(x)$ factor
in~$\FF_q[x]$. It turns out to be useful to introduce a supplementary
polynomial~$g(x)\in\FF_q[x]$ and study the factorization
of~$g\bigl(f^n(x)\bigr)$. In this generality, very little is known, or
even specifically conjectured, so our next question is necessarily
somewhat vague.

\begin{question}
\label{question:factorfnoverFq}
Let $f(x)\in\FF_q[x]$ be an \textup(irreducible\textup) polynomial,
let~$g(x)\in\FF_q[x]$, and for each~$n\ge1$,
factor~$g\bigl(f^n(x)\bigr)$ in~$\FF_q[x]$ as
\[
  g\bigl(f^n(x)\bigr) = p_1(x)^{e_1} p_2(x)^{e_2} \cdots p_r(x)^{e_r}
\]
with $p_i(x)$ irreducible in~$\FF_q[x]$.
\begin{parts}
\Part{(a)}
Describe the behavior, as $n\to\infty$, of quantities such as
\[
  r,\quad \sum_{i=1}^r e_i,\quad \max_{1\le i\le r} e_i,\quad \max_{1\le i\le r}\deg p_i(x).
\]
\Part{(b)}
Describe or predict the factorization of $g\bigl(f^{n+1}(x)\bigr)$
from a given factorization of $g\bigl(f^n(x)\bigr)$.  
\Part{(c)}
Answer the analogous questions for the iterates of a rational function
$f(x)\in\FF_q(x)$. Note that the exponents $e_1,\ldots,e_r$ may now be
both positive and negative, and hence it may be more natural to study $\sum |e_i|$
instead of $\sum e_i$, and $\max |e_i|$ instead of $\max e_i$.
\end{parts}
\end{question}

Under the assumptions that~$q$ is odd, that~$f(x)$ is quadratic,
that~$g(x)$ is irreducible of even degree, and that every
iterate~$f^n(x)$ is separable, Boston, Goksel, Jones, and
Xia~\cite{MR3359218,MR2888174} describe a conjectural multi-stage
Markov model to explain the factorization of successive values
of~$g\bigl(f^n(x)\bigr)$.  See also~\cite{MR3293428} for estimates for
the number and degree of the irreducible factors of iterates
of~$f(x)\in\FF_q[x]$, and see Section~\ref{section:isi} for a further
discussion of factorization of iterates over both finite fields and
number fields.

We might also look at orbits of dominant rational maps
$f:\PP^N\dashrightarrow\PP^N$ over finite fields, although just as in
Section~\ref{section:orm}, we need to avoid points~$P\in\PP^N(\FF_q)$
whose orbit lands in the indeterminacy locus~$I(f)$ of~$f$. We recall
our general notation~$\PP^N(K)_f$ for the set of points~$P\in\PP^N(K)$
for which~$f^n(P)\notin I(f)$ for all $n\ge1$.

\begin{conjecture}
\label{conjecture:IfoverFqk}
Let $f:\PP^N\dashrightarrow\PP^N$ be a dominant rational map defined
over~$\FF_q$.
\begin{parts}
\Part{(a)}
$\PP^N(\bar\FF_q)_f$ is Zariski dense in~$\PP^N$.
\Part{(b)}
$\displaystyle
  \lim_{k\to\infty} \#\PP^N(\FF_{q^k})_f/\#\PP^N(\FF_{q^k}) = 1.
$
\end{parts}
\end{conjecture}

\begin{remark}
We note that~(b) implies~(a) in
Conjecture~\ref{conjecture:IfoverFqk}. Assuming that the conjecture is
true, we may ask for a more accurate estimate for the ratio in~(b).
Thus let let~$m=\codim I(f)$. Then a rough probabilistic model
suggests that for most (all?) maps~$f$ and all $\e>0$, there should be
a constant~$C(f,\e)>0$ so that for all $k\ge1$,
\[
  1 - Cq^{-(m-\frac12+\e)k} \le
  \#\PP^N(\FF_{q^k})_f/\#\PP^N(\FF_{q^k}) \le 1 -
  Cq^{-(m-\frac12-\e)k}.
\]
\end{remark}

\begin{remark}
There is also the subject of \emph{permutation polynomials}, which are
polynomials $f(x)\in\FF_q[x]$ having the property that the map
$f:\FF_q\to\FF_q$ is a bijection. Since they act as permutations
of~$\FF_q$, one might expect that their dynamics mirrors that of random
automorphisms, rather than random endomorphisms. In any case, the
study of permutation polynomials is a field in its own right, with
many interesting open questions; see for example,
\cite{MR1541277,MR1542258} and \cite[Chapter~7]{MR1429394}.
\end{remark}

%%%%%%%%%%%%%%%%%%%%%%%%%%%%%%%%%%%%%%%%%%%%%%%%%%%%%%%%%%%%%%%%%%%%%%
\section{Irreducibilty and Stability of Iterates\WhoDone{Joe}}
\label{section:isi}
\WhoWrite{Joe}
%%%%%%%%%%%%%%%%%%%%%%%%%%%%%%%%%%%%%%%%%%%%%%%%%%%%%%%%%%%%%%%%%%%%%%

Let~$K$ be a number field, and let $f(x)\in K(x)$ be a rational
function of degree at least~$2$.  Write the~$n$th iterate of~$f(x)$
as $f^n(x)=F_n(x)/G_n(x)$ with relatively prime
polynomials~$F_n,G_n\in K[x]$.  Then for~$\a\in K$, we see that the polynomial
\[
  \text{$F_n(x)-\a G_n(x)$ is the numerator of the rational function
    $f^n(x)-\a$.}
\]

\begin{definition}
The pair~$(f,\a)$ is said to be \emph{stable over~$K$} if for
every~$n\ge1$, the polynomial $F_n(x)-\a G_n(x)$ is irreducible
in~$K[x]$.
\end{definition}

We observe that stability is not preserved under field
extension,\footnote{One might say that ``stability is not very
  stable!'' More seriously, there are two notions of stability in
  play, one being stability for varying~$n$, the other for varying base
  field.}  since for example~$(f,\a)$ is not stable over the
field~$K\bigl(f^{-1}(\a)\bigr)$. This prompts the following more
resilient version of stability, as described
in~\cite{MR2439638,MR3704363}.

\begin{definition}
\label{definition:eventualstab}
The pair~$(f,\a)$ is \emph{eventually stable over~$K$} if there is
an~$m\ge0$ such that for every~$\b\in f^{-m}(\a)$, the pair~$(f,\b)$
is stable over~$K(\b)$.
This is equivalent~\cite[Proposition~2.1(e)]{MR3704363} to the existence of 
a bound~$C(f,\a)$ such that the following holds:
\begin{equation}
  \label{eqn:FnalphaGn}
  \parbox{.65\hsize}{\noindent For every $n\ge1$, the polynomial
  $F_n(x)-\a G_n(x)$ is a product of at most $C(f,\a)$ irreducible
  factors in $K[x]$, counted with multiplicity.}
\end{equation}
%% *** This is correct, I think, but probably should be checked.
\end{definition}

\begin{example}
\label{example:evenstab4}
Let $f(x)=x^2+\frac{1}{48}$.
Then initially the number of factors grows,
\begin{align*}
  f^2(x) &=
  \left(x^2 - \tfrac{1}{2} x + \tfrac{7}{48}\right)  \left(x^2 + \tfrac{1}{2} x + \tfrac{7}{48}\right), \\
  f^3(x) &=
\left( x^2 - \tfrac{1}{2} x + \tfrac{19}{48} \right)
\left( x^2 + \tfrac{1}{2} x + \tfrac{19}{48} \right)
\left( x^4 - \tfrac{11}{24} x^2 + \tfrac{313}{2304} \right).
\end{align*}
Further, for all~$n\ge3$ one can show that~$f^n(x)$ has exactly three
irreducible factors in~$\QQ[x]$, so in~\eqref{eqn:FnalphaGn} we can
take $C(x^2+\frac{1}{48},0)=3$.  In particular,~$(x^2+\frac{1}{48},0)$
is eventually stable over~$\QQ$.
%% Let $f(x)=x^2-\frac{16}{9}$. Then initially the number of factors grows,
%% \begin{align*}
%%   f(x) &= \left(x-\tfrac43\right)\left(x+\tfrac43\right), \\
%%   f^2(x) &= \left(x-\tfrac23\right)\left(x+\tfrac23\right)\left(x^2-\tfrac{28}{9}\right), \\
%%   f^3(x) &= \left( x^2-2x+\tfrac29 \right)\left( x^2+2x+\tfrac29
%%   \right)\left( x^2-\tfrac{22}{9} \right)\left( x^2-\tfrac{10}{9}
%%   \right).
%% \end{align*}
%% However, for all~$n\ge3$ one can show~\cite{MR3335237} that~$f^n(x)$
%% has exactly four irreducible factors in~$\QQ[x]$,
%% so in~\eqref{eqn:FnalphaGn} we can take $C(x^2-\frac{16}{9},0)=4$.  In
%% particular,~$(x^2-\frac{16}{9},0)$ is eventually stable over~$\QQ$.
\end{example}

\begin{conjecture}
\label{conjecture:evenstab}
\textup{(Jones--Levy, 2017, \cite[Conjecture~1.2]{MR3704363})} Let $K$
be a number field, let $f(x)\in K(x)$ be a rational function of degree
at least~$2$, and let $\a\in K$ be a point that is not periodic
for~$f$.  Then~$(f,\a)$ is eventually stable over~$K$.
\end{conjecture}

One might be bold and ask if there is a uniform bound for the eventual
stability predicted by Conjecture~\ref{conjecture:evenstab}, but some
care is needed, since it is not hard to see that for any~$k\ge1$, the
constant~$C\bigl(f(x),f^k(\alpha)\bigr)$ in
Definition~$\ref{definition:eventualstab}$ must be at least~$k$.

\begin{question}
Let $f\in\End_d^1(K)$, and let~$\a\in K$ be a wandering point for~$f$
with the property that $f^{-1}(\a)\cap\PP^1(K)=\emptyset$ .  Is there
a constant $\kappa(d,K)$ so that the statement~$\eqref{eqn:FnalphaGn}$
in Definition~$\ref{definition:eventualstab}$ is true
with~$C(f,\a)=\kappa(d,K)$?  In particular, may we
take~$\kappa(2,\QQ)=3$; cf.\ Example~$\ref{example:evenstab4}$?
Further, if~$K$ is a number field, might the same statement be true
for a constant~$\kappa$ depending only on~$d$ and~$[K:\QQ]$?
\end{question}

\begin{remark}
Sookdeo~\cite{MR2782838} has shown that
Conjecture~\ref{conjecture:evenstab} for~$\a$ wandering follows from
the dynamical Lehmer conjecture
(Conjecture~\ref{conjecture:dynlehmer}).
\end{remark}

\begin{remark}
We note that stability is closely related to the arboreal
representation of the pair~$(f,\a)$;
cf.\ Section~\ref{section:ar}. Thus~$(f,\a)$ is stable over~$K$ if and
only if for every $n\ge1$ we have
\[
  \#f^{-n}(\a)=(\deg f)^n\quad\text{and}\quad
  \text{$\Gal(\Kbar/K)$ acts transitively on~$f^{-n}(\a)$.}
\]
Eventual stability has a similar interpretation since it is defined in
terms of stability of the pairs~$(f,\b)$ with~$\b\in f^{-m}(\a)$.
\end{remark}

There are various equivalent ways to extend the notion of stability to
finite endomorphisms of higher dimensional varieties such
as~$\PP^N$. We give one such definition.  Let~$X\subsetneq\PP^N$ be a proper
subvariety of dimension~$r$ that is defined over~$K$, e.g.,~$X$ could
be a point or a hypersurface. For any endomorphism~$f\in\End_d^N(K)$, the
pull-back of~$X$ has the form
\[
  f^*X
  = \sum_{\substack{Y/K\subset\PP^N,\;\dim(Y)=r\\\text{$Y$ irreducible over $K$}\\}}
  m(f,X;Y)Y,
\]
i.e.,~$f^*X$ is a formal sum of $K$-irreducible $r$-dimensional
subvarieties of~$\PP^N$ with multiplicities $m(f,X;Y)\in\ZZ$.

\begin{definition}
\label{definition:evenstabXinPN}
Let~$X\subsetneq\PP^N$ be an $r$-dimensional subvariety defined over~$K$, and let
$f:\PP^N\to\PP^N$ be an endomorphism defined over~$K$. The pair $(f,X)$
is \emph{eventually stable over~$K$} if there is a constant~$C(f,X)$
such that
\begin{equation}
  \label{eqn:evenstabXinPN}
  \sup_{n\ge1}  \sum_{\substack{Y/K\subset\PP^N,\;\dim(Y)=r\\\text{$Y$ irreducible over $K$}\\}}
  m(f^n,X;Y) \le C(f,X).
\end{equation}

The pair~$(f,X)$ is said to be \emph{geometrically eventually stable}
if it is eventually stable over~$\Kbar$, an algebraic closure of~$K$.
Geometrical eventual stability essentially never happens
if~$\dim(X)=0$, but should be common for \text{$\dim(X)\ge1$}.
\end{definition}

\begin{conjecture}
\label{conjecture:fPevenstab}
Let $K$ be a number field, let $f:\PP^N\to\PP^N$ be an endomorphism of
degree at least~$2$ defined over~$K$, and let $P\in\PP^N(K)$ be a
point that is not periodic for~$f$.  Then~$(f,P)$ is eventually stable
over~$K$.
\end{conjecture}

\begin{question}
\label{question:fXevenstab}
Let $K$ be a number field, let $f:\PP^N\to\PP^N$ be an endomorphism of
degree at least~$2$ defined over~$K$, and let $X\subsetneq\PP^N$ be a
subvariety that is defined over~$K$ and irreducible over~$K$.
\begin{parts}
\Part{(a)}
If the iterated forward orbit $\bigcup_{n\ge1} f^n(X)$ of~$X$ is
Zariski dense in~$\PP^N$, is it true that~$(f,X)$ is eventually
stable over~$K$?
\Part{(b)}
If the iterated forward orbit is Zariski dense and $\dim(X)\ge1$, is it
true that~$(f,X)$ is geometrically eventually stable?
\Part{(c)}
If $\dim(X)\ge1$, is eventual stability over~$K$ equivalent to
geometric eventual stability?
\end{parts}
\end{question}

The question of (eventual) stability, or lack thereof, for maps over
finite fields has also attracted attention.  We restrict to~$\PP^1$,
and via an affine change of coordinates $f(x+\a)-\a$, we note that it
suffices to discuss stability with base point~$0$.  An initial
question is to determine the number of stable rational maps or
polynomials.

\begin{question}
\label{question:ratstabFq}
We recall that~$\End_d^1(K)$ denotes the set of degree~$d$ rational
functions~$f\in K(x)$, and we write~$\Poly_d^1(K)$ for the set of
degree~$d$ polynomials~$f\in K[x]$. We define
\begin{align*}
\Poly_d^1(\FF_q)^\stab &:= \bigl\{ f\in\Poly_d^1(\FF_q) : \text{$(f,0)$ is stable over $\FF_q$} \bigr\}, \\  
\Rat_d^1(\FF_q)^\stab &:= \bigl\{ f\in\Rat_d^1(\FF_q) : \text{$(f,0)$ is stable over $\FF_q$} \bigr\},
\end{align*}
to be the indicated sets of $\FF_q$-stable maps.  For fixed~$d\ge2$,
how does the size of the sets $\Poly_d^1(\FF_q)^\stab$ and
$\Rat_d^1(\FF_q)^\stab$ grow as~$q\to\infty$?
Analogous questions apply for the sets of eventually stable maps,
and for maps of~$\AA^N$ or~$\PP^N$ with $N\ge2$.
\end{question}

Very little is known about Question~\ref{question:ratstabFq}.  For
polynomials of degree~$2$, Gomez and Nicol{\'a}s~\cite{MR2727344}
prove that there is an absolute constant~$C$ so that
\[
  \frac14(q-1)^2 \le \#\Poly_2^1(\FF_q)^\stab \le Cq^{5/2}\log q.
\]
Note that the upper bound is non-trivial, since
\[
  \Poly_2^1(\FF_q)=\{ax^2+bx+c:a\in\FF_q^*,\;b,c\in\FF_q\}
\]
has roughly~$q^3$ elements.

In a somewhat different vein, for a given finite set of
polynomials~$F\subset\FF_q[x]$, let~$\Ccal(F)$ denote the set of all
compositions of the elements of~$F$, where each~$f(x)\in F$ may be
used arbitrarily many times.  If~$F$ is a set of monic quadratic
polynomials, Ferraguti, Micheli, and Schnyder~\cite{MR3649090} give
necessary and sufficient conditions for~$\Ccal(F)$ to consist entirely
of irreducible polynomials, and they show~\cite{arxiv1701.06040} that
the set of irreducible polynomials in~$\Ccal(F)$ has the structure of
a regular language. Further, Heath-Brown and
Micheli~\cite{arxiv1701.05031} prove that for every~$k\ge1$ there 
exists a set~$F$ containing~$k$ distinct quadratic polynomials such
that~$\Ccal(F)$ consists entirely of irreducible polynomials.

In those cases that~$f(x)$ is not eventually stable, it is an
interesting problem to describe the limiting distribution of the
decomposition of~$f^n(x)$ into irreducible factors.  See
Question~\ref{question:factorfnoverFq} and the subsequent material in
Section~\ref{section:dff} for a brief discussion.

We might instead fix a map in~$\QQ(x)$, reduce it modulo~$p$, and
consider its stability over~$\FF_p$.

\begin{question}
Let $f(x)\in\QQ(x)$ be a rational function of degree~$d\ge2$, and for
each prime~$p$, let~$\tilde f_p\in\FF_p(x)$ denote its reduction modulo~$p$;
see Section~$\ref{section:goodred}$.
\begin{parts}
\Part{(a)}
Is it true that the set
\begin{equation}
  \label{eqn:stableprimes}
  \bigl\{ p : \text{$\tilde f_p(x)$ is stable over $\FF_p$} \bigr\}
\end{equation}
is a finite set?
\Part{(b)}
Is it true that the set
\begin{equation}
  \label{eqn:evenstableprimes}
  \bigl\{ p : \text{$\tilde f_p(x)$ is eventually stable over $\FF_p$} \bigr\}
\end{equation}
is a finite set?
\end{parts}
\end{question}

As a specific example, Jones~\cite[Conjecture~6.3]{arXiv:1012.2857}
has conjectured that~$x^2+1$ is stable over~$\FF_p$ if and only
if~$p=3$, i.e., for $f(x)=x^2+1$, the set~\eqref{eqn:stableprimes}
is~$\{3\}$.  Ferraguti~\cite{MR3787342} has shown that under some
Galois-theoretic assumptions, the set of stable
primes~\eqref{eqn:stableprimes} has density zero.

%%%%%%%%%%%%%%%%%%%%%%%%%%%%%%%%%%%%%%%%%%%%%%%%%%%%%%%%%%%%%%%%%%%%%%
\section{Primes, Prime Divisors, and Primitive Divisors in Orbits\WhoDone{Joe+Rafe+Patrick}}
\label{section:ppdds}
\WhoWrite{Joe+Rafe+Patrick}

The study of primes appearing in sequences is a recurring theme in
number theory, ranging from Dirichlet's theorem on primes in
arithmetic progressions to the conjectural existence of
infinitely many primes in the sequences~\text{$n^2+1$} and~\text{$2^n-1$}.
It is thus natural to look for primes appearing in orbits. A classical
example is provided by the orbit of~$3$ under the map~$f(x)=x^2-2x+2$,
which leads to the sequence of Fermat numbers
\begin{equation}
  \label{eqn:fermatnumber}
  f^n(3) = 2^{2^n}+1.
\end{equation}
In this case it is known that~$f^n(3)$ is prime for $0\le n\le 4$
and that~$f^n(3)$ is composite for $5\le n\le 32$.

\begin{remark}
We give a heuristic calculation to estimate the number of primes in an
orbit.  Let $f(x)\in\ZZ[x]$ have degree~$d\ge2$, and let~$\a\in\ZZ$ be
a wandering point for~$f$. Then as described in
Section~\ref{section:ch}, the canonical height~$\hhat_f(\a)$ is
strictly positive, which shows that
$\log\bigl|f^n(\a)\bigr|=d^n\hhat_f(\a)+O(1)$ grows quite rapidly.  If
we assume that the integers~$f^n(\a)$ in the orbit of~$\a$ are
``random and independent'' in an appropriate sense, then the prime
number theorem says that the probability that~$\bigl|f^n(\a)\bigr|$ is prime is
roughly $1/\log\bigl|f^n(\a)\bigr|$. This leads to the estimate
\begin{align*}
\#\bigl\{ n\ge0 : \text{$\bigl|f^n(\a)\bigr|$ is prime} \bigr\}
&\approx \sum_{n\ge 0} \Prob\big(\text{$\bigl|f^n(\a)\bigr|$ is prime} \bigr) \\*
&\approx \sum_{n\ge0} \frac{1}{\log\bigl|f^n(\a)\bigr|} \\*
&\approx \sum_{n\ge0} \frac{1}{d^n\hhat_f(\a)+O(1)} < \infty.
\end{align*}
The convergence of this series suggests that~$\Orbit_f(\a)$ should
contain only finitely many primes. However, and this is a major
caveat, our ``random and independent'' assumption has some definite
flaws. For example, the Fermat numbers~\eqref{eqn:fermatnumber} are
pairwise relatively prime, and their prime divisors have congruence
constraints.
\end{remark}

\begin{conjecture}
\label{conjecture:primesinorbit}
Let $f(x)\in\ZZ[x]$ be a polynomial of degree at least~$2$, and
let~$\a\in\ZZ$ be a wandering point for~$f$. Then~$\Orbit_f(\a)$
contains only finitely many primes.
\end{conjecture}

More generally, let~$K$ be a number field. For each prime
ideal~$\gp$ of the ring of integers of~$K$, let~$\ord_\gp:K^*\onto\ZZ$
be the associated normalized valuation, and
let~$\ord_\gp^{\scriptscriptstyle+}(\a)=\max\bigl\{\ord_\gp(\a),0\bigr\}$. Then we
have the following natural generalization of
Conjecture~\ref{conjecture:primesinorbit}.

\begin{conjecture}
\label{conjecture:nufnabc}
Let~$K$ be a number field, let $f(x)\in K(x)$ be a rational function
of degree at least~$2$, let~$\a\in K$ be a wandering point for~$f$,
let~$\b\in K$, and let~$c$ be a constant. Then the set
\[
  \left\{ n\ge0 : \sum_{\gp} \ord_\gp^{\scriptscriptstyle+}\bigl(f^n(\a)-\b\bigr) \le c \right\}
  \quad\text{is finite.}
\]
\end{conjecture}

Conjectures~\ref{conjecture:primesinorbit}
and~\ref{conjecture:nufnabc} suggest that orbits contain very few
primes. Further, even for sequences such as~$n^2+1$ and~$2^n-1$, the
existence of infinitely many primes is still conjectural.  For many
applications it suffices to know that there are a lot of different primes
dividing some term in the sequence, or that each successive term in
the sequence is divisible by a new prime. These considerations prompt
the following definitions.

\begin{definition}
 Let $\Acal=(a_1,a_2,a_3,\ldots)$ be a sequence of integers. The
 \emph{support\footnote{The fancy reason for the name ``Support'' is
     that a non-zero integer~$a$ may be viewed as the function on
     $\Spec(\ZZ)$ that sends the prime ideal~$(p)\in\Spec(\ZZ)$ to the
     value $\bar a\in\ZZ/p\ZZ$. The \emph{support of~$a$} is then the set of
     points of~$\Spec(\ZZ)$ at which~$a$ vanishes.} of~$\Acal$} is the
 set of primes
\[
  \Support(\Acal) := \{\text{primes $p$} : \text{$p$ divides at least one non-zero $a_n$} \}.
\]
\end{definition}

\begin{definition}
A prime~$p$ is said to be a \emph{primitive prime divisor of~$a_n$}
if~$p$ divides~$a_n$ and~$p$ does not divide~$a_k$ for all~$k<n$
with~$a_k\ne0$. The \emph{Zsigmondy set of~$\Acal$} is the set of
indices whose terms do not have a primitive prime divisor,
\[
  \Zsig(\Acal) := \{ n : \text{$a_n\ne0$ and $a_n$ has no primitive prime divisors} \}.
\]
\end{definition}

More generally, if~$\Acal$ is a sequence of rational numbers,
then~$\Support(\Acal)$ and $\Zsig(\Acal)$ are defined, respectively, to
be the support and the Zsigmondy set of the sequence of numerators
of~$a_1,a_2,\ldots$.

\begin{example}
Let $\a\in\QQ$ be a rational number with $\a\notin\{0,\pm1\}$.
In~1892 Zsigmondy~\cite{MR1546236}, generalizing an earlier weaker result of
Bang, proved that
\[
\Zsig\bigl((\a^n-1)_{n\ge1}\bigr)\subseteq\{1,2,3,6\}.
\]
Carmichael~\cite{MR1502458} proved in~1913 that the Zsigmondy set of
the Fibonacci sequence is $\{1,2,6,12\}$, and recently, Bilu, Hanrot,
and Voutier~\cite{MR1863855} generalized all of these old results by
showing that if~$\Acal$ is a Lucas sequence or a Lehmer sequence, then
$\Zsig(\Acal)\subseteq\{1,2,\ldots,30\}$.
\end{example}

Let $f(x)\in\QQ(x)$ be a rational function of degree at least~$2$, and
let~$\a\in\QQ$ be an $f$-wandering point. Then for any~$\b\in\QQ$, we
would like to understand the support and the Zsigmondy set of the
dynamically defined sequence
\[
  \Acal(f,\a,\b) := \bigl( f^n(\a) - \b \bigr)_{n\ge1}.
\]

\begin{theorem}
\label{theorem:zsigwanpre}
\textup{(Ingram--Silverman, 2009, \cite{MR2475968})} Let
$f(x)\in\QQ(x)$ be a rational function of degree at least~$2$,
let~$\a\in\QQ$ be an $f$-wandering point, and let~$\b\in\QQ$ be an
$f$-preperiodic point. Assume further that for all $k\ge1$, the
equation $f^k(x)=\b$ has at least one complex root other than
$x=\b$.
Then the Zsigmondy set of~$\Acal(f,\a,\b)$ is finite.
\end{theorem}

If the point~$\b$ in Theorem~\ref{theorem:zsigwanpre} is allowed to be
an~$f$-wandering point, then the problem becomes much more
difficult. Krieger~\cite[(2013)]{MR3142262} has shown that for all
$d\ge2$ and all $c\in\QQ$ with $c\notin\{0,-1,-2\}$,
\[
  \# \Zsig \bigl( \Acal(x^d+c,0,0) \bigr) \le 23.
\]
Notice that the point~$0$ is a critical wandering point
of~$x^d+c$.  There are stronger conditional results if one assumes
the~$ABC$-conjecture~\cite{MR3138487,arxiv1711.01507} or Vojta's
conjecture~\cite{MR3057058}.  See also~\cite{MR2863906} for a proof
that the Zsigmondy set of a sequence of the
form~$\bigl(f^n(\a)-f^{n-1}(\a)\bigr)_{n\ge1}$ is finite. 

\begin{conjecture}
\label{conjecture:zsigfinite}
Let $f(x)\in\QQ(x)$ be a rational function of degree at least~$2$, and
let~$\a,\b\in\QQ$ be $f$-wandering points.
Then the Zsigmondy set of~$\Acal(f,\a,\b)$ is finite.
\end{conjecture}

Conjecture~\ref{conjecture:zsigfinite} implies, in particular, that
there are infinitely many primes dividing at least one term in the
sequence~$\Acal(f,\a,\b)$.  Could it happen that every prime divides
at least one term? Or is the set of such primes relatively sparse?

There is a large literature studying the density of the support of
sequences such as $a^n-1$ and the Fibonacci sequence. In the dynamical
setting, it is interesting to study the density of the support of
orbits, as in the next theorem, which generalizes an earlier
result of Odoni~\cite{MR813379} for~$f(x)=x^2-x+1$. (We refer the
reader to Definition~\ref{definition:density}  for the definition of
the density of a set of primes.)

\begin{theorem}
\label{theorem:densupx2A}
\textup{(Jones, 2008, \cite{MR2439638})}
Let $\a\in\ZZ$. Then
\begin{align*}
\textup{(a)}&&   \Density\bigl( \Support\bigl( \Acal(x^2 + c,\a,0) \bigr) \bigr) &=0 &&\text{for all $c\in\ZZ\setminus\{-1\}$,} \\
\textup{(b)}&&   \Density\bigl( \Support\bigl( \Acal(x^2 - cx + c,\a,0) \bigr) \bigr) &=0 &&\text{for all $c\in\ZZ$.} 
\end{align*}
\end{theorem}

See~\cite{MR3335237} for a generalization of
Theorem~\ref{theorem:densupx2A} to polynomials of the form $x^d+c$
over number fields and function fields.  But before one guesses that
dynamical densities are always zero, we note~\cite{MR2679455,MR2439638}
\begin{align*}
  \Density\bigl( \Support\bigl( \Acal(x^2,q,q) \bigr) \bigr) & =1/3
  \quad\text{for all primes $q\ge3$,} \\
  \Density\bigl( \Support\bigl( \Acal(x^2-2,16,16) \bigr) \bigr) & = 7/24.
\end{align*}
As the following conjecture shows, it is no coincidence that these
examples use~$x^2$ and~$x^2-2$, since they are PCF polynomials. 

\begin{conjecture}
\label{conjecture:denssupdyn}
\textup{(Jones, 2008, \cite[Conjecture~4.7]{MR2439638})} Let
$f(x)\in\ZZ[x]$ be a monic polynomial of degree~$2$, and assume
that~$f(x)$ is stable\footnote{We recall from
  Section~\ref{section:isi} that~$f(x)$ is stable if every
  iterate~$f^n(x)$ is irreducible in~$\QQ[x]$.}  and that its critical
point is wandering.  Then
\[
  \Density\bigl( \Support\bigl( \Acal(f,\a,0) \bigr) \bigr) =0\quad\text{for all $\a\in\ZZ$.} 
\]
\end{conjecture}

\begin{question}
Give a concise characterization, analogous to that given in
Conjecture~$\ref{conjecture:denssupdyn}$ for maps
satisfying~$\deg(f)=2$, which ensures that~$f(x)\in\QQ(x)$
and~$\a,\b\in\QQ$ satisfy
\[
  \Density\bigl( \Support\bigl( \Acal(f,\a,\b) \bigr) \bigr)=0.
\]
\end{question}

%%%%%%%%%%%%%%%%%%%%%%%%%%%%%%%%%%%%%%%%%%%%%%%%%%%%%%%%%%%%%%%%%%%%%%
\section{Integral Points in Orbits\WhoDone{Joe}}
\label{section:ipo}
\WhoWrite{Joe}
%%%%%%%%%%%%%%%%%%%%%%%%%%%%%%%%%%%%%%%%%%%%%%%%%%%%%%%%%%%%%%%%%%%%%%

Let $f(X)\in\QQ(X)$ be a rational function of one variable, and
let~$\a\in\QQ$ be a rational starting point. A natural number
theoretic question is to ask whether the orbit~$\Orbit_f(\a)$ 
contains infinitely many integers. In some cases, for
example \text{$f(X)\in\ZZ[X]$} and~$\a\in\ZZ$, the entire orbit consists of
integers.  Similarly, for the map~\text{$f(X)=1/X^2$} and $\a\ge2$ an
integer, every other entry in the orbit is an integer, which simply
reflects the fact that~$f^2(X)=X^4$ is a polynomial. The following
proposition is well-known to dynamicists.  For a generalization to
endomorphisms of~$\PP^N$, see~\cite[Proposition~4.2]{MR1285389}
and~\cite[Exercise~7.16]{MR2316407}, and for a characteristic~$p$
version, see~\cite{silverman:polynomialiterate}.

\begin{proposition}
\label{proposition:fnpolyf2poly}
Let $f(X)\in\CC(X)$ be a rational function of degree~$d\ge2$ such that there
is an $n\ge1$ with $f^n(X)\in\CC[X]$. Then already~$f^2(X)\in\CC[X]$.
\end{proposition}

We start with an arithmetic analogue of
Proposition~\ref{proposition:fnpolyf2poly}, which is at the same
time a dynamical analogue of Siegel's theorem that curves of genus at
least~$1$ have only finitely many integral points.

\begin{theorem}
\label{theorem:intptsorbitP1}
\textup{(Silverman, 1993, \cite{silverman:dynamicalintegerpoints})}
Let $f(X)\in\QQ(X)$ be a rational function of degree at least~$2$,
and assume that~$f^2(X)$ is not a polynomial. Then for all~$\a\in\QQ$,
the orbit~$\Orbit_f(\a)$ contains only finitely many integers.
\end{theorem}

A strong version of Siegel's theorem, alluded to earlier, puts a limit
on the allowable non-integrality of points on curves. Here is a
dynamical analogue.

\begin{theorem}[Dynamical Siegel Theorem for $\PP^1$]
\label{theorem:dynsiegelP1}
\textup{(Silverman, 1993, \cite{silverman:dynamicalintegerpoints})}
Let $f(X)\in\QQ(X)$ be a rational function of degree at least~$2$,
and assume that both~$f^2(X)$ and $f^2(X^{-1})^{-1}$
are not polynomials. Let~$\a\in\QQ$
be a wandering point for~$f$, and write~$f^n(\a)$ as a fraction 
\[
  f^n(\a) = A_n/B_n\quad\text{with $A_n,B_n\in\ZZ$ and $\gcd(A_n,B_n)=1$.}
\]
Then
\[
  \lim_{n\to\infty} \frac{\log|A_n|}{\log|B_n|} = 1.
\]
\end{theorem}

In words, Theorem~\ref{theorem:dynsiegelP1} says that the numerator
and the denominator of~$f^n(\a)$ have roughly the same number of
digits when~$n$ is large.

Lang formulated a conjectural generalization of Siegel's theorem to
integral points on affine subvarieties of abelian varieties, and
Faltings~\cite{faltings:mordelllangconj,MR1307396} proved Lang's
conjecture. We give a conjectural dynamical analogue for~$\PP^N$. We
say that a point~$P\in\PP^N(\QQ)$ is an \emph{integral point}
(\emph{relative to the hyperplane~$x_0=0$}) if it has homogeneous
coordinates of the form
\[
  P = [1,a_1,a_2,\ldots,a_N]\quad\text{with $a_1,\ldots,a_N\in\ZZ$.}
\]

\begin{conjecture}[Dynamical Lang--Siegel Conjecture for $\PP^N$]
\label{conjecture:intptsorbitPN}
Let $f:\PP^N\to\PP^N$ be a morphism of degree~$d\ge2$ defined over~$\QQ$,
and assume that for all $1\le n\le N+1$, the first coordinate function of $f^n$ is not a constant
multiple of~$X_0^{d^n}$. Let $P\in\PP^N(\QQ)$. Then~$\Orbit_f(P)$ contains
only finitely points that are integeral relative to the hyperplane~$x_0=0$.
\par
More precisely, assume that~$P$ is wandering, and for each $n\ge1$, write
\begin{multline*}
  f^n(P) = \bigl[A_0(n),A_1(n),\ldots,A_N(n)\bigr] \\*
  \quad\text{with
    $A_0(n),\ldots,A_N(n)\in\ZZ$ and
    $\gcd\bigl(A_0(n),\ldots,A_N(n)\bigr)=1$.}
\end{multline*}
Then
\[
\lim_{n\to\infty} \frac{\log\bigl|A_0(n)\bigr|}
{\log\max\bigl\{ \bigl|A_0(n)\bigr|,\ldots,\bigl|A_N(n)\bigr| \bigr\} } = 1.
\]
\end{conjecture}

\begin{remark}
All of the theorems/conjectures in this sections can be proven/stated
for number fields~$K$ and for points that are integral in rings
of~$S$-integers of~$K$, but we have restricted attention to~$\QQ$ in
order to keep the notation and the technicalities to a minimum.
\end{remark}

%%%%%%%%%%%%%%%%%%%%%%%%%%%%%%%%%%%%%%%%%%%%%%%%%%%%%%%%%%%%%%%%%%%%%%
\section{Local--Global Questions in Dynamics\WhoDone{Joe}}
\label{section:lgqd}
\WhoWrite{Joe}

An enduring theme in mathematics is the use of local data to
characterize global phenomena. In number theory, local data can mean
information modulo~$m$ for all, or at least most, integers~$m$.  The
Chinese remainder theorem reduces studying solutions modulo~$m$ to
studying solutions modulo prime powers, and the latter are best
described in terms of solutions in the field of~$p$-adic
numbers~$\QQ_p$.

For example, a classical theorem of Legendre says that
for~$A,B,C\in\ZZ$, the equation
\[
AX^2 + BY^2 + CZ^2 = 0
\]
has a non-trivial solution in integers if and only if it has a
non-trivial solution in~$\QQ_p$ for every prime~$p$ and a non-trivial
solution in~$\RR$.  A modern example of a local-global principle is
the conjecture of Birch and Swinnerton-Dyer, which describes the
$\QQ$-rational points on an elliptic curve in terms of its~$\RR$
points and its~$\QQ_p$ points as~$p$ ranges over all primes.

There are various natural ways in which one might formulate
local-global problems in arithmetic dynamics. For example, to what
extent does the existence of $\QQ_p$ and~$\RR$ periodic points force
there to be $\QQ$ periodic points? The following result,
although covering only a few cases, gives a flavor for what we have in
mind.  (See also \cite{MR3654074,MR3062909}.)

\begin{theorem}  
\textup{(Krumm, 2016, \cite{MR3562026})} Let $f(x)=ax^2+bx+c\in\QQ[x]$
be a quadratic polynomial.  Let $n\in\{1,2,3\}$. Suppose that for
every prime~$p$, the polynomial~$f(x)$ has a $\QQ_p$-rational periodic
point of exact period~$n$. Then~$f(x)$ has a $\QQ$-rational periodic
point of exact period~$n$. Further, for~$n=4$ and, under an additional
assumption, for~$n=5$, there  exist infinitely many primes~$p$
such that~$f(x)$ does not have a $\QQ_p$-rational periodic point of
exact period~$n$. \textup(This may be compared with
Theorem~\textup{\ref{theorem:x2cperpt}(b)}.\textup)
\end{theorem}

We formulate the next question over~$\QQ$, but it is easily
generalized to number fields.

\begin{question}
Let $f:\PP^N\to\PP^N$ be a morphism of degree $d$ defined over~$\QQ$.
\begin{parts}
  \Part{(a)}
  For which $n\ge1$ is it true that~$f$ has a point of exact
  period~$n$ in~$\PP^N(\QQ)$ if and only if~$f$ has a point of exact
  period~$n$ in~$\PP^N(\RR)$ and in~$\PP^N(\QQ_p)$ for all primes~$p\/$\textup{?}
  \Part{(b)}
  For which $n\ge1$ is it true that there are infinitely many
  primes~$p$ such that the map~$f$ has no points of exact period~$n$
   in~$\PP^N(\QQ_p)$\textup{?} Is this true for all
  $n\ge n_0(N,d)$\textup{?}
\end{parts}
More generally, we may ask the same questions for preperiodic points
of type~$(m,n)$, i.e., points with tails of exact length~$m$ and
periods of exact length~$n$.
\end{question}

One can also formulate local-global questions related to the Dynamical
Mordell--Lang Conjecture (Conjecture~\ref{conjecture:dynmordelllang}).
We state the next conjecture using the topology of the adele ring of a
number field~$K$.  This ring amalgamates all of the archimedean and
$p$-adic data of~$K$. The conjecture says, roughly, that $p$-adic data
associated to the orbit of a point suffices to pin down the orbit.

\begin{conjecture}[Dynamical Brauer--Manin Conjecture]
\label{conjecture:dynBM}
\textup{(Hsia--Silverman, 2009, \cite{MR2537714})}
Let $K$ be a number field, let $X/K$ be a variety, let $f:X\to X$ be
a morphism defined over~$K$, let~$Y\subset X$ be a subvariety
defined over~$K$, and let~$P\in X(K)$.
Let~$\Acal_K$ be the ring of adeles of~$K$, and for any set~$S$,
let~$\Ccal(S)$ denote the adelic closure of~$S$. Then
\[
  \Ccal\bigl(Y\cap\Orbit_f(P)\bigr) = Y(\Acal_K) \cap \Ccal\bigl(\Orbit_f(P)\bigr).
\]
\end{conjecture}

Classically, a
Diophantine equation
\[
F(x_1,\ldots,x_n)=0, \quad\text{where $F\in\ZZ[x_1,\ldots,x_n]$,}
\]
is said to have a \emph{Brauer obstruction} to $\QQ$-rational
solutions if either it has no solutions in~$\RR$ or if there exists a
prime~$p$ such that the equation has no solutions in~$\QQ_p$. When
there is no Brauer obstruction, Manin devised a fancier test for the
non-existence of solutions. This has led to a considerable amount of
research studying for which equations the lack of a Brauer--Manin
obstruction suffices to ensure the existence of a $\QQ$-rational
solution. Conjecture~\ref{conjecture:dynBM} is a dynamical analogue of
Scharaschkin's~\cite{MR2700328} geometric formulation of the
Brauer--Manin obstruction for
curves. Conjecture~\ref{conjecture:dynBM} is known for $X=\PP^1$ and
for various other special cases; see~\cite{MR3424832,MR2496466}.

Here is another local-global formulation of dynamical Mordell--Lang.
We state it for polynomials, but it would be interesting to formulate
it for rational functions, or even for morphisms on~$\PP^N$.

\begin{question}
Let~$K$ be a number field, let $f(x),g(x)\in K[x]$ be polynomials of degree at least~$2$,
and let~$a,b\in K$.
\begin{parts}
\Part{(a)}
Suppose that for all but finitely many primes~$\gp$,
the mod~$\gp$ orbits of~$a$ and~$b$ satisfy
\begin{equation}
  \label{eqn:OfaOfbmodp}
  \Orbit_{\tilde f}(\tilde a \bmod \gp) \cap   \Orbit_{\tilde g}(\tilde b \bmod \gp)
  \ne \emptyset.
\end{equation}
Is it true that true that~$f$ and~$g$ have the same complex Julia
sets?\/\footnote{Aside from power maps and Chebyshev polynomials, having
  the same Julia set is more-or-less equivalent to~$f$ and~$g$ having
  a common iterate, up to symmetries of the Julia set.} If the answer
is affirmative, does it suffice to assume that~\eqref{eqn:OfaOfbmodp}
holds for a set of primes of density~$1$?
\Part{(b)}
Fix some $\e>0$, and suppose that the
intersection~\eqref{eqn:OfaOfbmodp} has at least~$\Norm\gp^\e$
elements for infinitely many primes~$\gp$. Is it true that true
that~$f$ and~$g$ have the same complex Julia sets?
\end{parts}
\end{question}

%% %%%%%%%%%%%%%%%%%%%%%%%%%%%%%%%%%%%%%%%%%%%%%%%%%%%%%%%%%%%%%%%%%%%%%%
%% \section{XXX}
%% \label{section:***}
%% \WhoWrite{}
%% %%%%%%%%%%%%%%%%%%%%%%%%%%%%%%%%%%%%%%%%%%%%%%%%%%%%%%%%%%%%%%%%%%%%%%

%%%%%%%%%%%%%%%%%%%%%%%%%%%%%%%%%%%%%%%%%%%%%%%%%%%%%%%%%%%%%%%%%%%%%%%%
% Acknowldegements and Bibliography
%%%%%%%%%%%%%%%%%%%%%%%%%%%%%%%%%%%%%%%%%%%%%%%%%%%%%%%%%%%%%%%%%%%%%%%%

\begin{acknowledgement}
The authors would like to thank Bryna Kra and Steven J. Miller for
suggesting that this article be written, and Eric Bedford, Laura
DeMarco, Mattias Jonsson, Juan Rivera-Letelier, Karou Sano, and Igor
Shparlinski for helpful advice.
\end{acknowledgement}

%% \begin{thebibliography}{99}
%% \itemsep=\smallskipamount
%% \end{thebibliography}

%% \bibliographystyle{plain}
%% \bibliography{/Users/jhs/Dropbox/AAJHS/Book/ADS/ArithDyn}

\def\cprime{$'$}

%%%%%%%%%%%%%%%%%%%%%%%%%%%%%%%%%%%%%%%%%%%%%%%%%%%%%%%%%%%%%%%%%%%%%%
%% \appendix
%%%%%%%%%%%%%%%%%%%%%%%%%%%%%%%%%%%%%%%%%%%%%%%%%%%%%%%%%%%%%%%%%%%%%%

%%%%%%%%%%%%%%%%%%%%%%%%%%%%%%%%%%%%%%%%%%%%%%%%%%%%%%%%%%%%%%%%%%%%%%
%% \section{Appendix Section}
%% \label{section:***}
%%%%%%%%%%%%%%%%%%%%%%%%%%%%%%%%%%%%%%%%%%%%%%%%%%%%%%%%%%%%%%%%%%%%%%

\end{document}